\documentclass[a4paper,10pt]{article}
\usepackage{amsthm,amsmath,amssymb,bbm,geometry,epsfig,listings,hyperref,enumerate,comment,nicefrac,verbatim,multirow,titlesec}
\newtheorem{lemma}{Lemma}[section]

\newtheorem{prop}[lemma]{Proposition}

\newcommand{\1}{\ensuremath{\mathbbm{1}}}

\providecommand{\N}{{\ensuremath{\mathbbm{N}}}}
\providecommand{\Z}{{\ensuremath{\mathbbm{Z}}}}
\providecommand{\R}{{\ensuremath{\mathbbm{R}}}}

\providecommand{\E}{{\ensuremath{\mathbb{E}}}}
\renewcommand{\P}{{\ensuremath{\mathbb{P}}}}
\newcommand{\eps}{{\ensuremath{\varepsilon}}}

\newcommand{\Var}{{\ensuremath{\operatorname{Var}}}}
\newcommand{\Lip}{{\ensuremath{\operatorname{Lip}}}}

\title{On  
multilevel Picard numerical approximations\\
for high-dimensional nonlinear parabolic partial\\
differential equations and high-dimensional nonlinear\\
backward stochastic differential equations}

\author{Weinan E$^{1,2} $, Martin Hutzenthaler$^{3}$, Arnulf Jentzen$^{4}$, \& Thomas Kruse$^{5}$
\bigskip
\\
\small{$^1$ Department of Mathematics and 
Program in Applied and Computational Mathematics,}
\\
\small{Princeton University, Princeton, NJ 08544-1000, USA,
e-mail: weinan@math.princeton.edu}
\smallskip
\\
\small{$^2$ BICMR and School of Mathematical Sciences,  Peking University, Beijing, China, 100870}
\smallskip
\\
\small{$^3$ Faculty of Mathematics, University of Duisburg-Essen,}
\\
\small{45117 Essen, Germany, e-mail: martin.hutzenthaler@uni-due.de}
\smallskip
\\
\small{$^4$ Seminar f\"ur Angewandte Mathematik,
ETH Zurich,}
\\
\small{8092 Z\"urich, Switzerland, e-mail: arnulf.jentzen@sam.math.ethz.ch}
\smallskip
\\
\small{$^5$ Faculty of Mathematics, University of Duisburg-Essen,}
\\
\small{45117 Essen, Germany, e-mail: thomas.kruse@uni-due.de}
}

\begin{document}

\maketitle
\makeatletter
\let\@makefnmark\relax
\let\@thefnmark\relax
\@footnotetext{\emph{AMS 2010 subject classification:} 65M75}
\@footnotetext{\emph{Key words and phrases:}
curse of dimensionality, high-dimensional PDEs, high-dimensional nonlinear BSDEs, multilevel Picard approximations, multilevel Monte Carlo method
  }
\makeatother

\begin{abstract}
Parabolic partial differential equations (PDEs) and backward stochastic differential equations (BSDEs)
are key ingredients in a number of models in physics and financial engineering. In particular, parabolic PDEs and BSDEs  
are fundamental tools in the state-of-the-art pricing and hedging of financial derivatives. 
The PDEs and BSDEs appearing in such 
applications are often high-dimensional and nonlinear. 
Since explicit solutions of such PDEs and BSDEs 
are typically not available, it is a very active topic of research to 
solve such PDEs and BSDEs approximately. 
In the recent article [E, W., Hutzenthaler, M., Jentzen, A.,  \& Kruse, T. Linear scaling algorithms for solving high-dimensional
nonlinear parabolic differential equations.\ {\it arXiv:1607.03295} (2017)] 
we proposed a family of approximation methods based on Picard approximations and multilevel Monte Carlo methods 
and showed under suitable regularity assumptions on the exact solution for semilinear heat equations
that the computational complexity is bounded by $O( d \, \eps^{-(4+\delta)})$ for any $\delta\in(0,\infty)$,
where $d$ is the dimensionality of the problem and $\eps\in(0,\infty)$ is the prescribed accuracy.
In this paper, we test the applicability of this algorithm on a variety of $100$-dimensional nonlinear PDEs that arise in physics and finance
by means of
numerical simulations presenting 
approximation accuracy against runtime.
The simulation results for these 100-dimensional example PDEs are very satisfactory in terms of accuracy and speed. In addition, we also provide a review of other approximation methods for nonlinear PDEs and BSDEs from the literature.
\end{abstract}

\tableofcontents

\section{Introduction and main results}

Parabolic partial differential equations (PDEs) and backward stochastic
differential equations (BSDEs) have a wide range of applications.
To give specific examples we focus now on a number of applications in finance.
There are several fundamental assumptions incorporated
in the Black-Scholes model that are not met in the real-life trading 
of financial derivatives.
A number of derivative pricing models have been developed in about the last four decades
to relax these assumptions;
see, e.g., 
\cite{Bergman1995,ElKarouiPengQuenez1997,BenderSchweizerZhuo2014,
LemorGobetWarin2006,GobetLemorWarin2005,LarentAmzelekBonnaud2014,
Crepeyetal2013}
for models taking into account
the fact that the ``risk-free'' bank account has higher interest 
rates for borrowing than for lending, particularly, 
due to the default risk of the trader,
see, e.g., 
\cite{Henry-Labordere2012,Crepeyetal2013}
for models incorporating the default risk 
of the issuer of the financial derivative,
see, e.g., 
\cite{WindcliffWangForsythVetzal2007,BayraktarYoung2008,BayraktarMilevskyPromislowYoung2009}
for models for the pricing of financial derivatives 
on underlyings which are not tradeable such as financial derivatives on the temperature
or mortality-dependent financial derivatives,
see, e.g., 
\cite{Amadori2003}
for models incorporating that the hedging strategy influences 
the price processes through demand and supply 
(so-called large investor effects),
see, e.g.,
\cite{ForsythVetzal2001,Leland1985,GuyonLabordere2011}
for models taking the transaction costs
in the hedging portfolio into account,
and see, e.g., 
\cite{AvellanedaLevyParas1995,GuyonLabordere2011}
for models incorporating 
uncertainties in the model parameters 
for the underlying.
In each of the above references the value function 
$ u $, describing the price of the financial derivative, 
solves a \emph{nonlinear} parabolic PDE. 
Moreover, 
the PDEs for the value functions emerging from the above models 
are often high-dimensional as the financial derivative depends in several cases 
on a whole basket of underlyings and as a portfolio containing several 
financial derivatives must often be treated as a whole 
in the case where the above nonlinear effects are taken into account
(cf., e.g.,
\cite{Crepeyetal2013,
ForsythVetzal2001,
BenderSchweizerZhuo2014}).
These high-dimensional nonlinear PDEs can typically not be solved explicitly
and, in particular, there is a strong demand from the financial engineering industry
to approximately compute the solutions of such high-dimensional nonlinear parabolic PDEs.

The numerical analysis literature contains a number of 
deterministic approximation methods for nonlinear parabolic PDEs
such as finite element methods, finite difference methods,
spectral Galerkin approximation methods,
or sparse grid methods
(cf., e.g., \cite[Chapter~14]{t97}, \cite[Section~3]{Tadmor2012}, \cite{Smolyak1963}, and \cite{PetersdorffSchwab2004}).
Some of these methods achieve high convergence rates 
with respect to the computational effort and, in particular, 
provide efficient approximations in low or moderate dimensions.
However, these approximation methods can not be used in high dimensions as 
the computational effort grows \emph{exponentially} in the dimension
$ d \in \N = \{ 1, 2, \dots \} $
of the considered nonlinear parabolic PDE
and then the method fails to terminate within years even for low  accuracies.

In the case of linear parabolic PDEs the \emph{Feynman-Kac formula} 
establishes an explicit representation of the solution of the PDE 
as the expectation of the solution of an appropriate 
stochastic differential equation (SDE).
(Multilevel) Monte Carlo methods 
together with suitable discretizations of the SDE (see, e.g.,
\cite{m55,kp92,hjk12,HutzenthalerJentzen2014})
then result in a numerical approximation method 
with a computational effort that grows at most polynomially in the dimension $ d \in \N $ of the PDE
and that grows up to an arbitrarily small order quadratically in the reciprocal of the approximation precision
(cf., e.g., \cite{GrahamTalay2013Stochastic_simulation_and_Monte_Carlo_methods,g08b,h98,heinrich01}).
These multilevel Monte Carlo approximations are, however, limited
to \emph{linear} PDEs as the classical Feynman-Kac formula provides only in the 
case of a linear PDE an explicit representation of the solution of the PDE.
For lower error bounds in the literature on random and deterministic numerical approximation methods for high-dimensional linear PDEs the reader is, e.g., referred to Heinrich~\cite[Theorem 1]{Heinrich2006}.

In the seminal papers~\cite{PardouxPeng1990,Peng1991,PardouxPeng1992},
Pardoux \& Peng 
developed 
the theory of nonlinear backward
stochastic differential equations 
and, in particular, established a considerably generalized nonlinear Feynman-Kac formula to 
obtain an explicit representation of the solution of a nonlinear parabolic PDE
by means of the solution of an appropriate BSDE;
see also Cheridito et al.~\cite{CheriditoSonerTouziVictoir2007}
for second-order BSDEs.
Discretizations of BSDEs, however, require suitable discretizations of nested
conditional expectations 
(see, e.g., \cite{BouchardTouzi2004,Zhang2004,FahimTouziWarin2011,GuoZhangZhuo2015,CheriditoSonerTouziVictoir2007}).
Discretization methods for these nested conditional expectations proposed in the literature include
the 'straight forward' Monte Carlo method,
the quantization tree method (see~\cite{BallyPages2003a}),
the regression method based on Malliavin calculus or based on kernel estimation (see~\cite{BouchardTouzi2004}),
the projection on function spaces method (see~\cite{GobetLemorWarin2005}),
the cubature on Wiener space method (see~\cite{CrisanManolarakis2012}),
and the Wiener chaos decomposition method (see~\cite{BriandLabart2014}).
None of these discretization methods has the property
that the computational effort of the method grows at most polynomially
both in the dimension and in the reciprocal of the prescribed accuracy
(see Subsections~\ref{ssec:MC}--\ref{ssec:WienerChaos} below for a detailed discussion).
We note that solving high-dimensional semilinear parabolic PDEs at single space-time points and 
solving high-dimensional nonlinear BSDEs
at single time points
is essentially equivalent due to the generalized nonlinear Feynman-Kac formula 
established by Pardoux \& Peng.
In recent years the concept of fractional smoothness in the sense of function spaces has been used for
studying variational properties of BSDEs.
This concept of fractional smoothness quantifies the propagation of singularities in time
and shows that certain non-uniform time grids are more suitable in the presence of singularities;
see, e.g.,
Geiss \& Geiss~\cite{GeissGeiss2004},
Gobet \& Makhlouf~\cite{GobetMakhlouf2010}
or
Geiss, Geiss, \& Gobet~\cite{GeissGeissGobet2012} for details.
Also these temporal discretization methods require suitable discretizations of nested conditional expectations
resulting in the same problems as in the case of uniform time grids.

Another probabilistic representation for solutions of some nonlinear parabolic PDEs with polynomial nonlinearities
has been established in Skorohod \cite{Skorohod1964} by means 
of \emph{branching diffusion processes}.
Recently, this classical representation
has been extended under suitable assumptions
in Henry-Labord{\`e}re \cite{Henry-Labordere2012} to more general analytic nonlinearities
and
in Henry-Labord{\`e}re et al.~\cite{Henry-LabordereEtAl2016}
to polynomial nonlinearities in the pair 
$ ( u(t,x), (\nabla_x u)(t,x) ) \in \R^{ 1 + d } $, $ (t,x) \in [0,T] \times \R^d $, 
where $ u $ is the solution of the PDE, $d\in\N$ is the dimension, and $T\in(0,\infty)$ is the time horizon.
This probabilistic representation has been successfully used 
in Henry-Labord{\`e}re~\cite{Henry-Labordere2012} (see also Henry-Labord{\`e}re, Tan, \& Touzi~\cite{Henry-LabordereTanTouzi2014})
and
in Henry-Labord{\`e}re et al.~\cite{Henry-LabordereEtAl2016}
to obtain a Monte Carlo approximation method for semilinear parabolic PDEs
with a computational 
complexity which is bounded by $O(d\,\eps^{-2})$
where $d$ is the dimensionality of the problem and $\eps\in(0,\infty)$ is the prescribed accuracy.
The major drawback of the branching diffusion method is 
its insufficient applicability,
namely
it requires
the terminal/initial condition of the parabolic PDE 
to be quite small (see Subsection~\ref{ssec:branchingdiffusions}
below
for a detailed discussion).

In the recent article \cite{EHutzenthalerJentzenKruse2016} we proposed
a family of approximation methods which we denote as multilevel Picard approximations
(see~\eqref{eq:def_scheme} for its definition and Section~\ref{sec:semilinear}
for its derivation).
Corollary 3.18 in \cite{EHutzenthalerJentzenKruse2016}
shows
under suitable regularity assumptions (including smoothness and Lipschitz continuity)
on the exact solution
that the computational complexity of this algorithm
is bounded by $ O( d \, \eps^{-(4+\delta)} ) $ for any $ \delta\in(0,\infty)$,
where $d$ is the dimensionality of the problem and $\eps\in(0,\infty)$ is the prescribed accuracy.
In this paper we complement the theoretical complexity analysis
of \cite{EHutzenthalerJentzenKruse2016}
with a simulation study.
Our simulations in Section~\ref{sec:numerics} indicate
that the computational complexity
grows at most linearly in the dimension and quartically in the reciprocal of the prescribed accuracy
also for several $ 100 $-dimensional nonlinear PDEs from physics and finance
with non-smooth and/or non-Lipschitz nonlinearities and terminal condition functions. 
The simulation results for these $ 100 $-dimensional example PDEs are very satisfactory 
in terms of accuracy and speed.

\subsection{Notation}
Throughout this article we frequently use the following notation.
We denote by
$
  \langle \cdot, \cdot \rangle \colon
  \left(
    \cup_{ n \in \N }
    (
    \R^n \times \R^n
    )
  \right)
  \to
  [0,\infty)
$
the function that satisfies
for all $ n \in \N $, $ v = ( v_1, \dots, v_n ) $, $ w = ( w_1, \dots, w_n ) \in \R^n $
that
$
  \langle
    v, w
  \rangle
    =
   \sum_{i=1}^n v_i w_i
$.
For every
topological space $(E,\mathcal E)$ we denote by $\mathcal{B}(E)$
the Borel-sigma-algebra on  $(E,\mathcal E)$.
For all measurable spaces $(A,\mathcal{A})$ and $(B,\mathcal{B})$
we denote by $\mathcal{M}(\mathcal A,\mathcal B)$ the set of $\mathcal{A}$/$\mathcal{B}$-measurable
functions from $A$ to $B$.
For every probability space $(\Omega,\mathcal{A},\P)$ we denote by
$\left\|\cdot\right\|_{L^2(\P;\R)}\colon\mathcal{M}(\mathcal{A},\mathcal{B}(\R))\to[0,\infty]$ the function that
 satisfies for all
$X\in\mathcal{M}(\mathcal{A},\mathcal{B}(\R))$ that
$\|X\|_{L^2(\P;\R)}=\sqrt{\E\!\left[|X|^2\right]}$.
For all metric spaces $(E,d_E)$ and $(F,d_F)$ we denote by
$\Lip(E,F)$ the set of all globally Lipschitz continuous functions from $E$ to $F$.
For every $d\in\N$ we denote by $\operatorname{I}_{\R^{d\times d}}$ the identity matrix in $\R^{d\times d}$
and
we denote by $\R^{d\times d}_{ \operatorname{Inv} }$ the set of invertible matrices in $\R^{d\times d}$.
For every $d\in\N$ and every $A\in\R^{d\times d}$ we denote by $A^{*}\in\R^{d\times d}$ 
the transpose of $A$.
For every $d\in\N$ and every $x=(x_1,\ldots,x_d)\in\R^d$ we denote by $\operatorname{diag}(x)\in\R^{d\times d}$
 the diagonal matrix with diagonal entries $x_1,\ldots,x_d$.
For every $T\in(0,\infty)$ we denote by $\mathcal{Q}_T$ the set given by
$\mathcal{Q}_T=\{w\colon[0,T]\to\R\colon w^{-1}(\R\backslash\{0\})\text{ is a finite set}\}$.
We denote by $\lfloor\cdot\rfloor\colon\R\to\Z$
and $[\cdot]^{+}\colon\R\to[0,\infty)$
 the functions that satisfy for all $x\in\R$ that
$\lfloor x\rfloor =\max(\Z\cap(-\infty,x])$
and $[x]^+=\max\{x,0\}$.

\section{Multilevel Picard approximations for high-dimensional semilinear PDEs}
\label{sec:semilinear}

In Subsection~\ref{sec:algorithm_semilinear}
below we define 
multilevel Picard
approximations
(see \eqref{eq:def_scheme} below)
in the case of semilinear PDEs
(cf.\ \eqref{eq:PDE_quasilinear_0} in Subsection~\ref{sec:derivation.semilinear} below).
These approximations have been introduced in \cite{EHutzenthalerJentzenKruse2016}.
In Subsection~\ref{sec:abstract_picture} we explain the abstract idea behind multilevel Picard approximations.
In Subsection~\ref{sec:derivation.semilinear} we derive a fixed-point equation for semilinear PDEs
which is based on the Feynman-Kac and Bismut-Elworthy-Li formulas.

\subsection{Abstract picture for multilevel Picard approximations}
\label{sec:abstract_picture}
Roughly speaking, a key idea in our approach to solve high-dimensional nonlinear
PDEs/BSDEs is to formulate the solution of the considered PDE/BSDE as the solution 
of a suitable fixed-point equation
and then to approximate the fixed-point by suitable 
\emph{multilevel Picard approximations}
(see \eqref{eq:abs_scheme} below).
We now first outline this idea in an abstract form
(see \eqref{eq:multilevel_fixed_point}--\eqref{eq:abs_scheme} below)
and, thereafter, we demonstrate (see  Subsections~\ref{sec:derivation.semilinear}--\ref{sec:algorithm_semilinear})
how this general idea 
is applied to high-dimensional semilinear PDEs 
and high-dimensional BSDEs.

Let $ ( \mathcal{V}, \left\| \cdot \right\|_{ \mathcal{V} } ) $ be an $ \R $-Banach space
and let $ {\bf \Phi} \colon \mathcal{V} \to \mathcal{V} $ be a contraction on 
$ ( \mathcal{V}, \left\| \cdot \right\|_{ \mathcal{V} } ) $.
The \emph{Banach fixed-point theorem} then ensures that there exists 
a unique fixed-point of $ {\bf \Phi} $,
that is, the Banach fixed-point theorem establishes the existence 
of a unique element $ {\bf u}_{ \infty } \in \mathcal{V} $
with the property that
\begin{equation}
\label{eq:abs_fixed_point}
  {\bf u}_{ \infty } = {\bf \Phi}( {\bf u}_{ \infty } ) .
\end{equation}
We think of $ \mathcal{V} $ as a suitable space of functions,
we think of \eqref{eq:abs_fixed_point} as a suitable reformulation 
of the considered deterministic PDE, 
and we think of $ {\bf u}_{ \infty } \in \mathcal{V} $
as the pair consisting of the solution of the considered PDE and of its spatial derivative.
Next let 
$ ({\bf u}_k)_{k\in\N_0} \subseteq \mathcal{V} $
be a sequence of fixed-point iterations associated to \eqref{eq:abs_fixed_point},
i.e., let 
$ ( {\bf u}_k )_{ k \in \N_0 } \subseteq \mathcal{V} $
satisfy for all $ k \in \N $
that
$
  {\bf u}_k = {\bf \Phi}( {\bf u}_{ k - 1 } )
$.
In examples we choose 
$ ( {\bf u}_k )_{ k \in \N_0 } \subseteq \mathcal{V} $
such that $ {\bf u}_0 $ is known explicitly or easily computable.
The Banach fixed-point theorem also ensures that
$ \lim_{ k \to \infty } {\bf u}_k = {\bf u}_{ \infty } $
and that this convergence happens \emph{exponentially fast}.
To introduce our 
multilevel Picard approximations,
we need suitable computable approximations of the typically nonlinear function
$ {\bf \Phi} \colon \mathcal{V} \to \mathcal{V} $,
that is,
we consider a family
$ \Psi^{ k, n }_{ l, \rho } \colon \mathcal{V} \to \mathcal{V} $,
$ l, k \in \N $,
$ n \in \{ 0, 1 \} $,
$ \rho \in (0,\infty) $,
of functions on $ \mathcal{V} $ 
and we think of 
$ \Psi^{ k, n }_{ l, \rho } $
for 
$ l, k \in \N $,
$ n \in \{ 0, 1 \} $,
$ \rho \in (0,\infty) $
as appropriate computable approximations of $ {\bf \Phi} $ in the sense that
$ \Psi^{ k, 0 }_{ l, \rho } $
and
$ \Psi^{ k, 1 }_{ l, \rho } $
get closer to $ {\bf \Phi} $ in a suitable sense as $ l \in \N $
(and $ k \in \N $ and $ \rho \in (0,\infty) $ respectively) gets larger.
This,
the fact that
$
  \lim_{ k \to \infty } {\bf u}_k = {\bf u}_{ \infty }
$,
and the fact that
for all
$
  k \in \N
$
it holds that
$
  {\bf u}_k = {\bf \Phi}( {\bf u}_{ k - 1 } )
$
then ensure for all sufficiently large $ k \in \N $ that
\begin{equation}
\label{eq:multilevel_fixed_point}
\begin{split}
  {\bf u}_{ \infty }
  \approx
  {\bf u}_k 
  =
  {\bf u}_1
  +
  \sum_{ l = 1 }^{ k - 1 }
  \left[
    {\bf u}_{ l + 1 }
    -
    {\bf u}_l
  \right]
& =
  {\bf \Phi}( {\bf u}_0 )
  +
  \sum_{ l = 1 }^{ k - 1 }
  \Big[
    {\bf \Phi}( {\bf u}_l )
    -
    {\bf \Phi}( {\bf u}_{ l - 1 } )
  \Big]
\\ &
\approx
  \Psi_{ k, \rho }^{ k, 0 }( {\bf u}_0 )
  +
  \sum_{ l = 1 }^{ k - 1 }
  \Big[
    \Psi^{ k, 0 }_{ k - l, \rho }( {\bf u}_l )
    -
    \Psi^{ k, 1 }_{ k - l, \rho }( {\bf u}_{ l - 1 } )
  \Big]
  .
\end{split}
\end{equation}
The first $ \approx $ in \eqref{eq:multilevel_fixed_point}
exploits the fact that 
$
  \lim_{ k \to \infty } {\bf u}_k = {\bf u}_{ \infty }
$
and the second 
$ \approx $
in \eqref{eq:multilevel_fixed_point}
uses our approximation assumption 
on 
$ \Psi^{ k, n }_{ l, \rho } $, 
$ l, k \in \N $, 
$ n \in \{ 0, 1 \} $,
$ \rho \in (0,\infty) $.
Display~\eqref{eq:multilevel_fixed_point}
suggests to introduce suitable numerical approximations 
of $ {\bf u}_k $, $ k \in \N $,
and $ {\bf u}_{ \infty } $, respectively.
More formally, let
$
  \psi_{ k, \rho } \colon \mathcal{V}^k \to \mathcal{V}  
$,
$ k \in \N $,
$ \rho \in (0,\infty) $,
be the functions which satisfy for all 
$ k \in \N $, $\rho\in(0,\infty)$, $ v_0, v_1, \dots, v_{ k - 1 } \in \mathcal{V} $
that
\begin{equation}
\label{eq:abs_psi_definition}
  \psi_{ k, \rho }( v_0, v_1, \dots, v_{ k - 1 } )
  =
  \sum_{ l = 0 }^{ k - 1 }
  \Big[
    \Psi_{ k - l, \rho }^{ k, 0 }\big( v_l \big)
    -
    \mathbbm{1}_{ \N }( l )
    \,
    \Psi_{ k - l, \rho }^{ k, 1 }\big( 
      v_{ [ l - 1]^{+} } 
    \big)
  \Big]
\end{equation}
and let $ {\bf U}_{ k, \rho } \in \mathcal{V} $,
$ k \in \N_0 $,
$ \rho \in (0,\infty) $,
satisfy
for all $ k \in \N $,
$ \rho \in (0,\infty) $
that
$
  {\bf U}_{ 0, \rho } = {\bf u}_0
$
and
\begin{equation}
\label{eq:abs_scheme}
  {\bf U}_{ k, \rho }
  =
  \psi_{ k, \rho }\big(
    {\bf U}_{ 0, \rho }
    ,
    {\bf U}_{ 1, \rho }
    ,
    \dots 
    ,
    {\bf U}_{ k - 1, \rho }
  \big)
  =
  \sum_{ l = 0 }^{ k - 1 }
  \Big[
    \Psi^{ k, 0 }_{ k - l, \rho }\big( 
      {\bf U}_{ l, \rho } 
    \big)
    -
    \mathbbm{1}_{ \N }( l )
    \,
    \Psi_{ k - l, \rho }^{ k, 1 }\big( 
      {\bf U}_{[l - 1]^{+}, \rho } 
    \big)
  \Big]
  .
\end{equation}
We would like to point out that the approximations 
in \eqref{eq:abs_scheme} are 
\emph{full history recursive} 
in the sense that for every $k\in\N$ and every $\rho\in(0,\infty)$ the ``full history'' 
$ {\bf U}_{ 0, \rho } $,
$ {\bf U}_{ 1, \rho } $,
$ \dots $,
$ {\bf U}_{ k - 1, \rho } $
needs to be computed recursively 
in order to compute 
$
  {\bf U}_{ k, \rho }
$.
Moreover, we note that the approximations in \eqref{eq:abs_scheme}
exploit multilevel/multigrid ideas (cf., e.g., \cite{h98,hs99,heinrich01,g08b}).
Typically multilevel ideas appear where the different levels correspond to 
approximations with different time or space step sizes 
while here different levels correspond to different stages 
of the fixed-point iteration.
This, in turn, results in numerical approximations which require 
the full history of the approximations.
A key feature of the approximations in \eqref{eq:abs_scheme}
is that -- depending on the choice of the space 
$ ( \mathcal{V} , \left\| \cdot \right\|_{\mathcal{V}} ) $,
the function $ {\bf \Phi} $,
and 
the functions 
$ 
  ( \Psi^{ k, n }_{ l, \rho } )_{ l, k \in \N, n \in \{ 0, 1 \}, \rho \in (0,\infty) } 
$
--
the approximations~\eqref{eq:abs_scheme} often preserve the 
frequently exceedingly high convergence speed of the $ ( {\bf u}_k )_{ k \in \N } $
to $ {\bf u}_{ \infty } $
while keeping the computational cost moderate compared
to the desired approximation precision.

\subsection{A fixed-point equation for semilinear PDEs}
\label{sec:derivation.semilinear}

To get a better understanding of the approximation scheme
introduced in Subsection~\ref{sec:algorithm_semilinear},
we present in this subsection a rough derivation
of a fixed-point equation on which
the scheme~\eqref{eq:def_scheme} is based on.
For this fixed-point equation we impose for simplicity of presentation
appropriate additional hypotheses that are not needed for the definition
of the scheme~\eqref{eq:def_scheme} 
(cf.\ \eqref{eq:PDE_quasilinear_0}--\eqref{eq:def_Phi_semilinear} in this subsection with Subsection~\ref{sec:algorithm_semilinear}).

Let 
$ T \in (0,\infty) $, 
$ d \in \N $, 
let 
$ 
  g \colon \R^d \to \R
$, 
$
  f \colon [0,T] \times \R^d \times \R \times \R^{d} \to \R  
$,
$
  u \colon [0,T] \times \R^d \to \R
$,
$
  \eta \colon \R^d \to \R^d
$,
$
  \mu \colon [0,T] \times \R^d \to \R^d
$,
and
$ 
  \sigma = ( \sigma_1,\ldots,\sigma_d )
  \colon [0,T] \times \R^d \to \R^{ d \times d }_{ \operatorname{Inv} } 
$
be sufficiently regular functions, 
 assume 
for all 
$ t \in [0,T) $, $ x \in \R^d $
that
$
  u(T,x) = g(x) 
$
and
\begin{multline}  
\label{eq:PDE_quasilinear_0}
 ( \tfrac{ \partial }{ \partial t } u )( t, x )
  + 
  f\big( 
    t, x, 
    u( t, \eta( x ))    , 
   [ \sigma(t,\eta(x))]^{*}(\nabla_x u) ( 
      t, \eta( x ) 
    ) 
  \big)
  +
  \langle 
    \mu( t, x ) 
    ,
    ( \nabla_x u )( t, x ) 
  \rangle
\\
  +
  \tfrac{ 1 }{ 2 }
  \operatorname{Trace}\!\big(
    \sigma(t,x) [ \sigma(t,x) ]^*
    ( \operatorname{Hess}_x u)( t, x )
  \big) 
  = 0
  ,   
\end{multline}
let
$
  ( 
    \Omega, \mathcal{F}, \P, ( \mathbb{F}_t )_{ t \in [0,T] } 
  )
$
be a stochastic basis 
(cf., e.g., \cite[Appendix~E]{PrevotRoeckner2007}), let
$
  W = 
$
$
  ( W^{ 1 }, \dots, 
$
$  
  W^{ d } ) 
  \colon 
  [0,T] \times \Omega \to \R^d
$
be a standard $ ( \mathbb{F}_t )_{ t \in [0,T] } $-Brownian motion,
and for every $ s \in [0,T] $, $ x \in \R^d $
let 
$ 
  X^{ s, x }
  \colon [s,T] \times \Omega \to \R^d
$
and
$
  D^{ s, x } 
  \colon [s,T] \times \Omega \to \R^{ d \times d }
$
be $ ( \mathbb{F}_t )_{ t \in [s,T] } $-adapted stochastic processes 
with continuous sample paths which satisfy
that for all $ t \in [s,T] $ 
it holds $ \P $-a.s.\ that
\begin{equation}  
\label{eq:def_XD}
\begin{split} 
  X^{ s, x }_t
&
  = 
  x + \int_s^t \mu(r, X^{ s, x }_r ) \, dr
  + 
  \sum_{ j = 1 }^d
  \int_s^t \sigma_j(r, X^{ s, x }_r ) \, d W^j_r,
\\
  D^{ s, x }_t 
  &
  =
  \operatorname{I}_{\R^{d\times d}} 
  +
  \int_s^t (\tfrac{\partial}{\partial x}\mu)( r,X^{ s, x }_r ) \, D^{ s, x }_r \, dr
  +
  \sum_{ j = 1 }^d
  \int_s^t (\tfrac{\partial}{\partial x}\sigma_j)( r,X^{ s, x }_r ) \, D^{ s, x }_r \, d W^{ j }_r
  .
\end{split}     
\end{equation}
For every $s\in[0,T]$ the processes $ D^{ s, x } $, $x\in \R^d$, are in a suitable sense the 
\emph{derivative processes} of $ X^{ s, x } $, $x\in \R^d$ with respect to $ x \in \R^d $.
The function $ \eta $ in \eqref{eq:PDE_quasilinear_0} allows 
to include a possible space shift in the PDE.
Typically we are interested in the case where $ \eta $ is the identity, that is,
for all $x\in\R^d$ it holds that
$
\eta(x) = x
$.
Our approximation scheme in \eqref{eq:def_scheme} below is based on a suitable  
\emph{fixed-point formulation} of the solution of the PDE~\eqref{eq:PDE_quasilinear_0}.
To obtain such a fixed-point formulation,
we apply the \emph{Feynman-Kac formula} and 
the {\it Bismut-Elworthy-Li formula} 
(see, e.g., Elworthy \& Li~\cite[Theorem~2.1]{elworthy1994formulae} or Da Prato \& Zabczyk~\cite[Theorem~2.1]{DaPratoZabczyk1997}).
More precisely, 
let
$ 
  {\bf u}^{ \infty }
  \in
  \Lip( [0,T] \times \R^d, \R^{ 1 + d } )
$
satisfy
for all $ (t,x) \in [0,T) \times \R^d $
that
$
  {\bf u}^{ \infty }( t, x ) 
  = 
  \big( u(t,x), [\sigma(t,x)]^{*}(\nabla_x u )( t, x ) \big)
$
and let 
$ 
  \Phi \colon 
  \Lip( [0,T] \times \R^d, \R^{ 1 + d } )
  \to
  \Lip( [0,T] \times \R^d, \R^{ 1 + d } )
$
satisfy
for all 
$ 
  {\bf v} \in 
  \Lip( [0,T] \times \R^d, \R^{ 1 + d } )
$,
$ (s,x) \in [0,T) \times \R^d $
that
\begin{equation}
\label{eq:def_Phi_semilinear}
\begin{split}
  \big( \Phi( {\bf v} ) \big)( s, x )
& =
\textstyle
  \E\!\left[ 
    \left( 
      g( X^{ s, x }_T ) - g(x) 
    \right)
    \left( 
      1, 
      \tfrac{ [\sigma(s,x)]^{*} }{ T - s } 
      \smallint\nolimits_s^T
      \big[ 
        \sigma( r, X_r^{ s, x } )^{ - 1 } 
        D_r^{ s, x } 
      \big]^{ * }
      d W_r
    \right)
  \right]
  +
    \left( g(x), 0 \right)
\\ 
&
\quad
\textstyle
  +
    \int_s^T 
    \E\!\left[
      f\Big(
        t, X_t^{ s, x }, {\bf v}\big(t, \eta( X_t^{ s, x } ) \big)
      \Big) 
      \,
      \big(
      1 ,
      \tfrac{ [\sigma(s,x)]^{*} }{ t - s }
      \smallint\nolimits_s^t
      \big[ 
        \sigma( r, X_r^{ s, x } )^{ - 1 } 
        D_r^{ s, x } 
      \big]^{ * } 
      \, d W_r 
      \big)
    \right]
  dt
\end{split}
\end{equation}
Combining \eqref{eq:def_Phi_semilinear} with 
the \emph{Feynman-Kac formula}
and the {\it Bismut-Elworthy-Li formula} 
ensures that 
$
  {\bf u}^{ \infty } = \Phi( {\bf u}^{ \infty } ) 
$.
Note that we have incorporated a zero expectation term in \eqref{eq:def_Phi_semilinear}.
The purpose of this term is to slightly reduce the variance 
when approximating the right-hand side of \eqref{eq:def_Phi_semilinear}
by Monte Carlo approximations.
Now we approximate
the non-discrete quantities in~\eqref{eq:def_Phi_semilinear}
(expectation and time integral) by discrete quantities (Monte Carlo averages and quadrature formulas)
with different degrees of discretization on different levels
(cf.\ the remarks in Subsection \ref{sec:remarks.on.scheme} below).
This yields a family of approximations 
of $\Phi$.
With these approximations of $\Phi$ we finally define multilevel Picard approximations of ${\bf u}^{\infty}$ through~\eqref{eq:abs_scheme}
which results in the approximations~\eqref{eq:def_scheme}.

\subsection{The approximation scheme}
\label{sec:algorithm_semilinear}

In this subsection we introduce 
multilevel Picard approximations in the case 
of semilinear PDEs
(see \eqref{eq:def_scheme} below).
To this end we consider the following setting.

Let 
$ T \in (0,\infty) $, 
$ d \in \N $, 
$
  \Theta = \cup_{ n \in \N } \R^n
$,
let 
$ 
  g \colon \R^d \to \R
$, 
$
  f \colon [0,T] \times \R^d \times \R^{d+1} \to \R  
$,
$
  \eta \colon \R^d \to \R^d
$,
$
  \mu \colon [0,T] \times \R^d \to \R^d
$,
$ 
  \sigma
  \colon [0,T] \times \R^d \to \R^{ d \times d }_{ \operatorname{Inv} } 
$
be measurable functions,
let 
$ 
  ( q^{ k, l, \rho }_s )_{ k, l \in \N_0, \rho \in (0,\infty), s \in [0,T) } 
  \subseteq \mathcal{Q}_T 
$,
$
  ( m^g_{ k, l, \rho } )_{ k, l \in \N_0, \rho \in (0,\infty) }$,
$
  ( m^f_{ k, l, \rho } )_{ k, l \in \N_0, \rho \in (0,\infty) } \subseteq \N
$,
let
$
  ( 
    \Omega, \mathcal{F}, \P, ( \mathbb{F}_t )_{ t \in [0,T] } 
  )
$
be a stochastic basis,
let 
$
  W^{ \theta } \colon [0,T] \times \Omega \to \R^d 
$,
$ \theta \in \Theta $,
be independent
standard $ ( \mathbb{F}_t )_{ t \in [0,T] } $-Brownian motions
with continuous sample paths,
for every 
$ l \in \Z $,
$ \rho \in (0,\infty) $,
$ \theta \in \Theta $,
$ x \in \R^d $,
$ s \in [0,T) $, $ t \in [s,T] $
let 
$
  \mathcal{X}^{ l, \rho, \theta }_{ x, s, t } 
  \colon  \Omega \to \R^d
$,
$
  \mathcal{D}^{ l, \rho, \theta }_{ x, s, t }
  \colon  \Omega \to \R^{ d \times d }
$,
and
$
  \mathcal{I}^{ l, \rho, \theta }_{ x, s, t }
  \colon  \Omega \to \R^{ 1+d  }
$
be functions,
and 
for every
$ \theta \in \Theta $,
$
  \rho \in (0,\infty)
$
let
$ 
  {\bf U}^{ \theta }_{ k, \rho } 
  \colon [0,T]\times\R^d \times \Omega \to \R^{ d + 1 } 
$,
$
  k \in \N_0 
$,
be 
functions
which satisfy
for all 
$k\in\N$,
$ (s,x) \in [0,T)\times\R^d $
that
\begin{equation}  \begin{split}
\label{eq:def_scheme}
  {\bf U}^{ \theta }_{ k, \rho }( s, x )
&= 
  \sum_{ l = 0 }^{ k - 1} 
  \sum_{ i = 1 }^{ m^g_{ k, l , \rho } }
  \frac{ 1 }{
    m^g_{ k, l, \rho } 
  }
  \,
  \big[
    g(
      \mathcal{X}_{ x, s, T }^{ l, \rho , (\theta, l, -i) }
    )
    -
    \1_{ \N }( l ) \,
    g(
      \mathcal{X}_{ x, s, T }^{ l - 1, \rho, (\theta, l, -i)}
    )
    -
    \1_{ \{ 0 \} }( l ) \,
    g(x)
  \big]
  \,
  \mathcal{I}^{ l, \rho, ( \theta, l, - i ) }_{ x, s, T }
\\
&
  +
  \big( g(x), 0 \big)
  +
  \sum_{ l = 0 }^{ k - 1 }
  \sum_{ i = 1 }^{
    m^f_{ k, l, \rho } 
  }
  \sum_{ t \in [s,T] }
  \frac{ 
    q^{ k, l, \rho }_s( t )
  }{
    m^f_{ k, l, \rho } 
  }
  \,
  \Big[
    f\Big(
      t, 
      \mathcal{X}^{ k - l ,
        \rho,  
        ( \theta, l, i ) 
      }_{ x, s, t }
      , 
      {\bf U}^{ ( \theta, l, i , t) }_{ l, \rho }\big( 
        t,
        \eta(
          \mathcal{X}_{ x, s, t }^{ 
            k - l , \rho , 
            ( \theta, l, i ) 
          }
        )
      \big)
    \Big)
\\ 
& 
   -
   \1_{ \N }( l )
   \,
    f\Big(
      t, 
      \mathcal{X}^{ 
        k - l , \rho ,
        ( \theta, l, i ) 
      }_{ x, s, t }
      , 
      {\bf U}^{ ( \theta, -l, i, t ) }_{ [ l - 1 ]^{+} , \rho }\big( 
        t,
        \eta(
          \mathcal{X}_{ x, s, t }^{ 
            k - l , \rho ,
            ( \theta, l, i ) 
          }
        )
      \big)
    \Big)
  \Big]
  \,
  \mathcal{I}^{ k - l, \rho, (\theta, l, i) }_{ x, s, t }
  .
\end{split}     \end{equation}
\subsection{Remarks on the approximation scheme}
\label{sec:remarks.on.scheme}
In this subsection we add a few comments on the numerical 
approximations~\eqref{eq:def_scheme}.
For this we assume the setting in Section~\ref{sec:algorithm_semilinear}.
The set $ \Theta $ allows to index families
of independent random variables which we need 
for Monte Carlo approximations.
The natural numbers 
$ ( m^g_{ k, l, \rho } )_{ k, l \in \N, \rho \in (0,\infty) }, ( m^f_{ k, l, \rho } )_{ k, l \in \N_0, \rho \in (0,\infty) } \subseteq \N $
specify the number of Monte Carlo samples
in the corresponding levels 
for approximating the expectations involving $g$ and $f$
appearing on the right-hand side 
of \eqref{eq:def_Phi_semilinear}.
The family
$ 
  ( q^{ k, l, \rho }_s )_{ k, l \in \N_0, \rho \in (0,\infty), s \in [0,T) } \subseteq \mathcal{Q}_T 
$
provides the quadrature formulas that we employ to approximate
the time integrals
$ \int_s^T \dots dt $, $s\in[0,T]$, appearing on the right-hand
side of \eqref{eq:def_Phi_semilinear}.
In Subsections~\ref{subsec:counterparty}--\ref{subsec:borrowlend}
these parameters satisfy that
for every $k,l\in\N_0$, $\rho\in\N$ it holds
that $m_{k,l,\rho}^g=\rho^{k-l}$, $m_{k,l,\rho}^f=\operatorname{round}(\sqrt{\rho}^{k-l})$
and that for every $k,l\in\N_0$, $\rho\in\N$
it holds that $q^{k,l,\rho}$ is a Gau\ss-Legendre quadrature rule with 
$\operatorname{round}(\Gamma^{-1}(\rho^{(k-l)/2}))$
nodes.
In Subsections~\ref{subsec:allencahn}--\ref{subsec:exp}
these parameters satisfy that
for every $k,l\in\N_0$, $\rho\in\N$ it holds 
that $m_{k,l,\rho}^g=m_{k,l,\rho}^f=\rho^{k-l}$
and that for every $k,l\in\N_0$, $\rho\in\N$
it holds that $q^{k,l,\rho}$ is a Gau\ss-Legendre quadrature rule with 
$\operatorname{round}(\Gamma^{-1}(\rho^{(k-l)/2}))$
nodes.
For every $l\in\N$, $\rho\in(0,\infty)$, $\theta\in\Theta$, $(s,x)\in[0,T]\times\R^d$, $v\in (s,T]$
we think of the processes
$
  ( \mathcal{X}^{ l, \rho, \theta }_{ x, s, t } )_{ t \in [s,T] }
$
and 
$
  ( \mathcal{D}^{ l, \rho, \theta }_{ x, s, t } )_{ t \in [s,T] }
$
as $
  (\mathbb{F}_t)_{t\in[s,T]}
$-optional measurable computable approximations with
$\P \big( \int_{s}^T\big\|
    \sigma( r,
      \mathcal{X}_{ x, s, r }^{ l, \rho, \theta } 
    )^{ - 1 }
    \,
    \mathcal{D}_{ x, s, r }^{ l, \rho, \theta) }
  \big\|_{L(\R^d,\R^d)}^2\,dr<\infty \big )=1$
  (e.g., piecewise constant 
c\`{a}dl\`{a}g Euler-Maruyama approximations)
of the processes
$ 
  ( X^{ s, x }_t )_{ t \in [s,T] }
$
and 
$ 
  ( D^{ s, x }_t )_{ t \in [s,T] }
$
in \eqref{eq:def_XD}
and 
we think of
$ \mathcal{I}_{ x, s, v }^{ l, \rho, \theta }$
as a random variable that
satisfies $\P$-a.s.\ that
\begin{equation}\label{eq:I}
  \mathcal{I}^{ l, \rho, \theta }_{ x, s, v }
  =
  \left(
  1 ,
  \tfrac{ [\sigma(s,x)]^{*} }{ v - s }
  \smallint\nolimits_s^v
  \big[
    \sigma( r,
      \mathcal{X}_{ x, s, r }^{ l, \rho, \theta } 
    )^{ - 1 }
    \,
    \mathcal{D}_{ x, s, r }^{ l, \rho, \theta }
  \big]^{ * }
  dW_r^{ \theta }
  \right)
  .
\end{equation}
Note that if $\mathcal{X}_{x,s,\cdot}^{k,\rho,\theta}$ and
$\mathcal{D}_{x,s,\cdot}^{k,\rho,\theta}$ are piecewise constant
then
 the stochastic integral on the right-hand side of~\eqref{eq:I} reduces to a stochastic Riemann-type sum
which is not difficult to compute.
Observe that our approximation scheme~\eqref{eq:def_scheme} 
employs Picard fixed-point approximations (cf., e.g., \cite{BenderDenk2007}),
multilevel/multigrid techniques (see, e.g., \cite{g08b,h98,heinrich01,cdmr09}),
discretizations of the SDE system~\eqref{eq:def_XD},
as well as quadrature approximations for the time integrals.
Roughly speaking, the numerical approximations~\eqref{eq:def_scheme} are 
full history recursive in the sense that for every $(k,\rho)\in\N\times(0,\infty)$ 
the full history 
$ {\bf U}^{ ( \cdot ) }_{ 0, \rho } $,
$ {\bf U}^{ ( \cdot ) }_{ 1, \rho } $,
$ \dots $,
$ {\bf U}^{ ( \cdot ) }_{ k - 1, \rho } $
needs to be computed recursively 
in order to compute 
$
  {\bf U}^{ ( \cdot ) }_{ k, \rho }
$.
In this sense the numerical approximations~\eqref{eq:def_scheme}
are full history recursive multilevel Picard approximations.
\subsection{Special case: semilinear heat equations}
In this subsection
we specialize
the numerical scheme~\eqref{eq:def_scheme}
to the case of semilinear heat equations.
\begin{prop}\label{p:specialcase.abm}
Assume the setting in Section~\ref{sec:algorithm_semilinear},
assume for all
$
  k \in \N_0
$,
$ \rho \in (0,\infty) $,
$ \theta \in \Theta $,
$ x \in \R^d $,
$ s \in [0,T) $,
$ 
  t \in [s,T]
$,
$ 
  u \in (s,T]
$
that
$
  \eta( x ) = x
$,
$
  \mathcal{X}^{ k, \rho, \theta }_{ x, s, t } =
  x + W^{ \theta }_t - W^{ \theta }_s 
$,
$
  \mathcal{D}_{ x, s, t }^{ k, \rho, \theta }
  =
  \sigma( s, x )
  =
  \operatorname{I}_{\R^{d\times d}}
$,
$
  \mathcal{I}_{x,s,s}^{k,\rho,\theta}=0
$,
$
  \mathcal{I}_{x,s,u}^{k,\rho,\theta}=(1,\tfrac{W_u^{\theta}-W_s^{\theta}}{u-s})
$,
and for every 
$ \theta \in \Theta $,
$ a \in [0,T] $, 
$ b \in [a,T] $
let 
$ 
  \Delta W^{ \theta }_{ a, b } 
  \colon \Omega \to \R^d
$
be the function given by
$ 
  \Delta W^{ \theta }_{ a, b } 
  = W_b^{ \theta } - W^{ \theta }_a 
$.
Then it holds 
for all 
$ \theta \in \Theta $,
$ k\in\N$, $\rho \in (0,\infty) $,
$ (s,x) \in [0,T)\times\R^d $
that
\begin{equation}  
\label{eq:scheme_Laplace}
\begin{split}
  {\bf U}^{ \theta }_{ k, \rho }( s, x )
&= 
  \big(
    g(x)
    , 0
  \big)
  +
  \sum_{ i = 1 }^{ m^g_{ k, 0, \rho } }
  \frac{ 1 }{
    m^g_{ k, 0, \rho }
  }
  \Big[
    g\big(
      x 
      + 
      \Delta W^{ ( \theta, 0, -i) }_{ s, T } 
    \big)
    -
    g(x)
  \Big]
  \Big(
  1 ,
  \tfrac{ 1 }{ T - s }
  \Delta W^{ ( \theta, 0, -i) }_{ s, T } 
  \Big)
\\ &\quad
  +
  \sum_{ l = 0 }^{ k - 1 }
  \sum_{ i = 1 }^{
    m^f_{ k, l , \rho }
  }
  \sum_{ t \in (s,T] }
  \frac{ 
    q^{ k, l , \rho }_s( t )
  }{
    m^f_{ k, l, \rho }
  }
  \Big[
    f\big(
      t, 
      x 
      + 
      \Delta W^{ ( \theta, l, i) }_{ s, t } 
      , 
      {\bf U}^{ ( \theta, l, i, t ) }_{ l, \rho }( 
        t,
        x + 
        \Delta W^{ ( \theta, l, i) }_{ s, t } 
      )
    \big)
\\ &\quad
   -
   \1_{ \N }( l )
    f\big(
      t, 
      x + 
      \Delta W^{ ( \theta, l, i) }_{ s, t } 
      , 
      {\bf U}^{ ( \theta, - l, i, t ) }_{ [ l - 1 ]^{ + } , \rho }( 
        t,
        x + 
        \Delta W^{ ( \theta, l, i) }_{ s, t } 
      )
    \big)
  \Big]
  \big(
    1 ,
    \tfrac{ 1 }{ t - s }
    \Delta W^{ ( \theta, l, i) }_{ s, t } 
  \big)
  .
\end{split}     
\end{equation}
\end{prop}
\noindent
The proof of Proposition~\ref{p:specialcase.abm} is clear and therefore omitted.

\subsection{Special case: geometric Brownian motion}
In this subsection
we specialize
the numerical scheme~\eqref{eq:def_scheme}
to the case of the forward diffusion being a geometric Brownian motion.
This case often appears in the financial engineering literature.
\begin{prop}\label{p:specialcase.gbm}
Assume the setting in Section~\ref{sec:algorithm_semilinear},
let $\bar{\mu}\in\R$, $\bar{\sigma}\in(0,\infty)$,
for every 
$ \theta \in \Theta $,
$ a \in [0,T] $, 
$ b \in [a,T] $
let 
$ 
  \Delta W^{ \theta }_{ a, b } 
  \colon \Omega \to \R^d
$
be the function given by
$ 
  \Delta W^{ \theta }_{ a, b } 
  = W_b^{ \theta } - W^{ \theta }_a 
$,
and
assume for all
$
  k \in \N_0
$,
$ \rho \in (0,\infty) $,
$ \theta \in \Theta $,
$ x \in (0,\infty)^d $,
$ s \in [0,T) $,
$ 
  t \in [s,T]
$,
$ 
  u \in (s,T]
$
that
$
  \eta( x ) = x
$,
$
  \mathcal{D}_{ x, s, t }^{ k, \rho, \theta }
  =
  \exp((\bar{\mu}-\tfrac{\bar{\sigma}^2}{2})(t-s))\exp(\bar{\sigma}\operatorname{diag}(\Delta W^{\theta}_{s,t}))
$,
$
  \mathcal{X}^{ k, \rho, \theta }_{ x, s, t } =
  \mathcal{D}_{ x, s, t }^{ k, \rho, \theta } x
$,
$
  \sigma(s,x)=\bar{\sigma}\operatorname{diag}(x)
$,
$
  \mathcal{I}_{x,s,s}^{k,\rho,\theta}=0
$,
$
  \mathcal{I}_{x,s,u}^{k,\rho,\theta}=(1,\tfrac{1}{u-s}\Delta W_{s,u}^{\theta})
$.
Then it holds 
for all 
$ \theta \in \Theta $,
$ k\in\N$ , $\rho \in (0,\infty) $,
$ (s,x) \in [0,T)\times(0,\infty)^d $
that
\begin{equation}  
\label{eq:scheme_Laplace2}
\begin{split}
&
  {\bf U}^{ \theta }_{ k, \rho }( s, x )
= 
  \big(
    g(x)
    , 0
  \big)
  +
  \sum_{ i = 1 }^{ m^g_{ k, 0, \rho } }
  \frac{ 1 }{
    m^g_{ k, 0, \rho }
  }
  \Big[
    g\big(
    \mathcal{X}^{ 0, \rho, (\theta,0,-i) }_{ x, s, T }
    \big)
    -
    g(x)
  \Big]
  \Big(
  1 ,
  \tfrac{ 1 }{ T-s }
  \Delta W^{ ( \theta, 0, -i) }_{ s, T } 
  \Big)
\\ &
  +
  \sum_{ l = 0 }^{ k - 1 }
  \sum_{ i = 1 }^{
    m^f_{ k, l , \rho }
  }
  \sum_{ t \in (s,T] }
  \frac{ 
    q^{ k, l , \rho }_s( t )
  }{
    m^f_{ k, l, \rho }
  }
  \Big[
    f\big(
      t, 
      \mathcal{X}_{ x, s, t }^{ 
            k - l , \rho , 
            ( \theta, l, i ) 
          }
      , 
      {\bf U}^{ ( \theta, l, i, t ) }_{ l, \rho }( 
        t,
        \mathcal{X}_{ x, s, t }^{ 
            k - l , \rho , 
            ( \theta, l, i ) 
          }
      )
    \big)
\\ &
   -
   \1_{ \N }( l )
    f\big(
      t, 
      \mathcal{X}_{ x, s, t }^{ 
            k - l , \rho , 
            ( \theta, l, i ) 
          }
      , 
      {\bf U}^{ ( \theta, - l, i, t ) }_{ [ l - 1 ]^{ + } , \rho }( 
        t,
       \mathcal{X}_{ x, s, t }^{ 
            k - l , \rho , 
            ( \theta, l, i ) 
          }
      )
    \big)
  \Big]
  \big(
    1 ,
    \tfrac{ 1 }{ t-s }
    \Delta W^{ ( \theta, l, i) }_{ s, t } 
  \big)
  .
\end{split}     
\end{equation}
\end{prop}
\noindent
The proof of Proposition~\ref{p:specialcase.gbm} is clear and therefore omitted.
In the setting of Proposition~\ref{p:specialcase.gbm} we note that 
for all 
$k\in\N$, $\rho\in(0,\infty)$, $\theta\in\Theta$, $(s,x)\in[0,T)\times\R^d$, $t\in(s,T]$
it holds $\P$-a.s.\ that
\begin{equation}\label{eq:xdgbm}  \begin{split}
  \mathcal{X}_{x,s,t}^{k,\rho,\theta}
  &=x+\int_s^t\bar{\mu}\mathcal{X}_{x,s,r}^{k,\rho,\theta}\,dr
  +\int_s^t\bar{\sigma}\operatorname{diag}(\mathcal{X}_{x,s,r}^{k,\rho,\theta})\,dW_r^{\theta}
\\
  \mathcal{D}_{x,s,t}^{k,\rho,\theta}
  &=\operatorname{I}_{\R^{d\times d}}+\int_s^t\bar{\mu}\mathcal{D}_{x,s,r}^{k,\rho,\theta}\,dr
  +\int_s^t\bar{\sigma}\operatorname{diag}(\mathcal{D}_{x,s,r}^{k,\rho,\theta})\,dW_r^{\theta}
\\
  \mathcal{I}_{x,s,t}^{k,\rho,\theta}
  &=
  \left(
  1 ,
  \tfrac{ [\sigma(s,x)]^{*} }{ t - s }
  \int_s^t
  \big[
    \sigma( r,
      \mathcal{X}_{ x, s, r }^{ k, \rho, \theta } 
    )^{ - 1 }
    \,
    \mathcal{D}_{ x, s, r }^{ k, \rho, \theta }
  \big]^{ * }
  dW_r^{ \theta }
  \right).
\end{split}     \end{equation}

\section{Numerical simulations of high-dimensional nonlinear PDEs}
\label{sec:numerics}

In this section we apply the algorithm \eqref{eq:def_scheme} to 
approximate the solutions of several nonlinear PDEs; see 
Subsections~\ref{subsec:counterparty}--\ref{subsec:exp} below.
The solutions of the PDEs in
Subsections~\ref{subsec:counterparty}--\ref{subsec:allencahn} are not known explicitly. 
The solution of the PDE in Subsection~\ref{subsec:exp} is known explicitly. 
In Subsections~\ref{subsec:counterparty}--\ref{subsec:allencahn} 
the algorithm is tested for a one-dimensional and a one hundred-dimensional version
of a PDE. In the one-dimensional cases in Subsections~\ref{subsec:counterparty}--\ref{subsec:allencahn} we present
the error of our algorithm relative to a high-precision approximation of
the exact solution of the PDE provided by a finite difference approximation scheme 
(see the left-hand sides of Figures~\ref{fig:counter},~\ref{fig:cva},~\ref{fig:borrowlend},
 and~\ref{fig:allencahn} and Tables~\ref{tab:counterd1},~\ref{tab:cvad1},~\ref{tab:borrowlendd1},
and~\ref{tab:allencahnd1} below).
In the one hundred-dimensional cases in Subsections~\ref{subsec:counterparty}--\ref{subsec:allencahn}
we present the approximation increments of our scheme to analyze 
the performance of our scheme
in the case of high-dimensional PDEs
(see the right-hand sides of Figures~\ref{fig:counter},~\ref{fig:cva},~\ref{fig:borrowlend},
 and~\ref{fig:allencahn} and Tables~\ref{tab:counterd100},~\ref{tab:cvad100},~\ref{tab:borrowlendd100},
and~\ref{tab:allencahnd100} below).
In Subsection~\ref{subsec:exp} we employ the explicit formula for the solution
of the considered one hundred-dimensional PDE (see \eqref{eq:exp_sol} below) to present 
the error of our scheme
relative to the explicitly known exact solution (see the left-hand side of Figure~\ref{fig:exp} and Table~\ref{tab:expd100}).
Moreover, for each of the PDEs in Subsections~\ref{subsec:counterparty}--\ref{subsec:exp}
 we illustrate the growth of the computational effort with respect
to the dimension
by running the algorithm for each PDE for every dimension $d\in \{5,6,\ldots, 100\}$
 and recording the associated runtimes (see Figures~\ref{fig:countermultidim} and \ref{fig:borrowlendmultidim} and the 
 right-hand side of Figure~\ref{fig:exp}). All simulations are performed with 
 {\sc Matlab} on a 2.8 GHz Intel i7 processor with 16 GB RAM.

Throughout this section assume the setting in Subsection \ref{sec:algorithm_semilinear}, let $x_0\in \R^d$,
 let $u\in C^{1,2}([0,T]\times \R^d, \R)$ be a function which satisfies for all $(t,x)\in [0,T)\times \R^d$
that
$
  u(T,x) = g(x) 
$
and
\begin{multline}  
\label{eq:PDE_quasilinear_1}
  ( \tfrac{ \partial }{ \partial t } u )( t, x )
  + 
  f\big( 
    t, x, 
    u( t, \eta( x ))    , 
    [\sigma(t,\eta(x))]^{*}(\nabla_x u) ( 
      t, \eta( x ) 
    ) 
  \big)
  +
  \langle 
    \mu( t, x ) 
    ,
    ( \nabla_x u )( t, x ) 
  \rangle
\\
  +
  \tfrac{ 1 }{ 2 }
  \operatorname{Trace}\!\big(
    \sigma(t,x) [ \sigma(t,x) ]^*
    ( \operatorname{Hess}_x u)( t, x )
  \big) 
  = 0
  ,   
\end{multline}
and
assume for all $\theta \in \Theta$, $\rho \in (0,\infty)$, $(s,x)\in [0,T)\times \R^d$ that
\begin{equation}
{\bf U}^{ \theta }_{ 0, \rho }( s, x )
=  \big( g(x), 0 \big)+
  \sum_{ i = 1 }^{ m^g_{ 0, 0 , \rho } }
  \frac{ 1 }{
    m^g_{ 0, 0, \rho } 
  }
  \,
  \big[
    g(
      \mathcal{X}_{ x, s, T }^{ 0, \rho , (\theta, 0, -i) }
    )
    -
    g(x)
  \big]
  \,
  \mathcal{I}^{ 0, \rho, ( \theta, 0, - i ) }_{ x, s, T }.
\end{equation} 
To obtain smoother results we average over 10 independent simulation runs. More precisely, 
 for the numerical results in Subsections~\ref{subsec:counterparty}--\ref{subsec:borrowlend},
 for every $d\in \{1,100\}$ we run {\sc Matlab} code \ref{code:testrun} twice to produce one 
 realization of
 \begin{equation}
 \{ 1, 2, \dots, 7 \}\times \{ 1, 2, \dots,10 \} \ni (\rho,i) \mapsto {\bf U}^i_{ \rho,  \rho }(0, x_0 )=({\bf U}^{i,[1]}_{ \rho ,  \rho   }( 0,x_0 ),{\bf U}^{i,[2]}_{ \rho ,  \rho }( 0,x_0 ),\ldots,{\bf U}^{i,[d+1]}_{ \rho , \rho }(0, x_0 ))\in \R^{d+1},
 \end{equation}
 where in the 
second run, line 2 of  {\sc Matlab} code \ref{code:testrun} is replaced by \texttt{rng(2017)} 
to initiate the pseudorandom number generator with a different 
seed.
\sloppy{Moreover, for the numerical results in Subsections \ref{subsec:allencahn} and \ref{subsec:exp}, we run {\sc Matlab} code \ref{code:testrun} once, where
lines 4, 5, and 14 are replaced by \texttt{average=10;}, 
\texttt{rhomax=5;}, and \texttt{[a,b]=approximateUZabm(n(rho),rho,zeros(dim,1),0);}, respectively.}
\begin{lstlisting}[
	caption={{\sc Matlab} code to perform a testrun.},
	label=code:testrun]
global Mf Mg Q c w T dim f g mu sigma eta;
rng(2016)
format long
average=5;
rhomax=7;
rhomin=1;
[T,dim,f,g,eta,mu,sigma]=modelparameters();
[Mf,Mg,Q,c,w,n] = approxparameters(rhomax);
value=zeros(average,rhomax);
time=value;
for rho=rhomin:rhomax
  for k=1:average
    tic
    [a, b]=approximateUZgbm(n(rho),rho,100*ones(dim,1),0);
    value(k,rho)=a;
    time(k,rho)=toc;
  end
end
name = [datestr(now, 'yymmddTHHMMSS') '.mat']; 
save(name,'n','Q','Mf','Mg','value','time');
\end{lstlisting}
\sloppy{{\sc Matlab} code \ref{code:testrun} calls the {\sc Matlab} functions
 \texttt{approximateUZgbm} (respectively \texttt{approximateUZabm}), 
 \texttt{modelparameters}, and \texttt{approxparameters}.
The {\sc Matlab} functions \texttt{approximateUZgbm} and
 \texttt{approximateUZabm}
 are presented 
in  {\sc Matlab} codes~\ref{code:approximategbm} and \ref{code:approximateabm} and 
implement the schemes~\eqref{eq:scheme_Laplace2} and \eqref{eq:scheme_Laplace}, respectively.} More precisely,
up to rounding errors and the fact that random numbers are replaced by pseudo random numbers, it holds 
for all $\theta\in \Theta $, $n\in \N_0$, $\rho \in \N$, $x\in \R^d$, $s\in [0,T)$ that
$\texttt{approximateUZgbm}(n,\rho,x,s)$ returns one realization of ${\bf U}^\theta_{ n,  \rho }(s, x)$ satisfying \eqref{eq:scheme_Laplace2}.
Moreover, up to rounding errors and the fact that random numbers are replaced by pseudo random numbers, it holds 
for all $\theta\in \Theta $, $n\in \N_0$, $\rho \in \N$, $x\in \R^d$, $s\in [0,T)$ that
$\texttt{approximateUZabm}(n,\rho,x,s)$ returns one realization of ${\bf U}^\theta_{ n,  \rho }(s, x)$ satisfying \eqref{eq:scheme_Laplace}.
\lstinputlisting[
	caption={A {\sc Matlab} function with input $\theta\in \Theta$, $n\in \N_0$, $\rho \in \mathbb N $, $ x \in \mathbb R^d $, 
	$ t \in [0,T) $ and output one 
	realization of ${\bf U}^\theta_{n,\rho}(t,x)$ satisfying \eqref{eq:scheme_Laplace2}.},
	label=code:approximategbm]
	{Simulation/approximateUZgbm.m}
	
\lstinputlisting[
	caption={A {\sc Matlab} function with input $\theta\in \Theta$, $n\in \N_0$, $\rho \in \mathbb N $, $ x \in \mathbb R^d $, 
	$ t \in [0,T) $ and output one 
	realization of ${\bf U}^\theta_{n,\rho}(t,x)$ satisfying \eqref{eq:scheme_Laplace}.},
	label=code:approximateabm]
	{Simulation/approximateUZabm.m}

The {\sc Matlab} function \texttt{modelparameters} 
called in line 7 of {\sc Matlab} code \ref{code:testrun}
returns 
 the parameters $T\in (0,\infty)$, $d\in \N$, $f\colon [0,T]\times \R^d \times \R \times \R^{d} \to \R$, 
 $g\colon \R^d \to \R$, $\eta\colon \R^d \to \R$,
  $\bar \mu\in \R$, and $\bar \sigma \in \R$
 for each example considered in Subsections~\ref{subsec:counterparty}--\ref{subsec:exp}.  
 For the implementations of the  {\sc Matlab} function \texttt{modelparameters} we refer
 to {\sc Matlab} code \ref{code:datacounterd100} in Subsection \ref{subsec:counterparty}, 
 to {\sc Matlab} code \ref{code:datacvad100} in Subsection \ref{subsec:cva},
 to {\sc Matlab} code \ref{code:databorrowlendd100} in Subsection \ref{subsec:borrowlend},
 to {\sc Matlab} code \ref{code:dataallencahnd100} in Subsection \ref{subsec:allencahn},
 and to {\sc Matlab} code \ref{code:dataexpd100} in Subsection \ref{subsec:exp}.

The  {\sc Matlab} function \texttt{approxparameters}
 called in line 8 of {\sc Matlab} code \ref{code:testrun}
 provides for every example considered in Subsections~\ref{subsec:counterparty}--\ref{subsec:borrowlend}
 (respectively Subsections~\ref{subsec:allencahn}--\ref{subsec:exp})
 and every $\rho \in \{1,2,\ldots,  7\}$ (respectively $\rho \in \{1,2,\ldots, 5\}$)
 the numbers of Monte-Carlo samples $(m^g_{k,l,\rho})_{k,l \in \N_0}$ and 
 $(m^f_{k,l,\rho})_{k,l \in \N_0}$ and the quadrature formulas $(q^{k,l,\rho}_s)_{k,l \in \N_0, s\in [0,T)}$.
More precisely, throughout this section assume
for every
$
  s \in [0,T], 
  k,l \in \N_0, 
  \rho \in \N
$
with
$k\ge l$
that
$
  q^{ k, l, \rho }_s
$
is the Gauss-Legendre
quadrature formula on $ (s,T) $
with $\operatorname{round}(\varphi(\rho^{(k-l)/2}))$ nodes, where
$\varphi\colon [1,\infty) \to [2,\infty)$ is an approximation of the inverse gamma function which is 
given by {\sc Matlab} code \ref{code:invgamma}. To compute the Gauss-Legendre nodes
 and weights we use the {\sc Matlab} function \texttt{lgwt} that was written by Greg von Winckel and that can be downloaded from
\url{www.mathworks.com}. In addition, for every
$
  k,l \in \N_0, 
  \rho \in \N
$
we choose in Subsections~\ref{subsec:counterparty}--\ref{subsec:borrowlend}
that
$
  m^f_{ k, l, \rho }
  =
  \operatorname{round}(\rho^{(k-l)/2})
$
and 
$
  m^g_{ k, l, \rho }
  =
  \rho^{k-l}
$
and in Subsections~\ref{subsec:allencahn}--\ref{subsec:exp}
that 
$
  m^f_{ k,l, \rho }
  =
  \rho^{k-l}
$
and 
$
  m^g_{ k, l, \rho }
  =
  \rho^{k-l}
$.
For the numerical results in Subsections~\ref{subsec:counterparty}--\ref{subsec:borrowlend} {\sc Matlab} code \ref{code:approxparameterscounter} presents
the implementation of \texttt{approxparameters}. For the numerical results in 
Subsections~\ref{subsec:allencahn}--\ref{subsec:exp}
 line 10 in {\sc Matlab} code \ref{code:approxparameterscounter} is replaced by \texttt{Mf(rho,k)=rho\^{}k;}. The reason for choosing in Subsections \ref{subsec:counterparty}--\ref{subsec:borrowlend} fewer Monte-Carlo samples $(m^f_{k,l,\rho})_{k,l \in \N_0, \rho \in \N}$ than in 
 Subsections~\ref{subsec:allencahn}--\ref{subsec:exp} is that in the former cases 
 for every $s \in [0,T)$ the variance $ \Var(f(s,X^{0,x_0}_s, \E[g(X^{s,x}_T)(1,\frac{W_T-W_s}{T-s})]\big|_{x=X^{0,x_0}_s}))$ of the nonlinearity is of smaller magnitude than the variance
$\Var(g(X^{0,x_0}_T))$ of the terminal condition. Therefore, the nonlinearity requires fewer Monte-Carlo samples to obtain a Monte-Carlo error of 
the same magnitude as the terminal condition. Averaging the nonlinearity less saves 
 computational effort and allows to employ a higher maximal number of Picard iterations ($7$ in Subsections~\ref{subsec:counterparty}--\ref{subsec:borrowlend} compared to $5$ in 
 Subsections~\ref{subsec:allencahn}--\ref{subsec:exp}).

\begin{lstlisting}[
	caption={A {\sc Matlab} function that returns the approximation parameters.},
	label=code:approxparameterscounter]
function [Mf,Mg,Q,c,w,n] = approxparameters(rhomax)
    global T;
    n=1:1:rhomax;
    Q=zeros(rhomax);
    Mf=Q;
    Mg=zeros(rhomax,rhomax+1);
    for rho=1:rhomax
        for k=1:n(rho)
            Q(rho,k)=round(inverse_gamma(rho^(k/2)));
            Mf(rho,k)=round((rho)^((k)/2));
            Mg(rho,k)=rho^(k-1);
        end
        Mg(rho,rho+1)=rho^rho;
    end
    qmax=max(max(Q));
    c=zeros(qmax);
    w=c;
    for k=1:qmax
        [ctemp,wtemp] = lgwt(k,0,T);
        c(:,k)=[flip(ctemp);zeros(qmax-k,1)];
        w(:,k)=[flip(wtemp);zeros(qmax-k,1)];
    end
end
\end{lstlisting}
\begin{lstlisting}[
	caption={{\sc Matlab} function to approximate the inverse Gamma function.},
	label=code:invgamma]
function y=inverse_gamma(x)
    c=0.036534;
    L= log((x+c)/sqrt(2*pi));
    y=L/lambertw(L/exp(1))+0.5;
end
\end{lstlisting}

Solutions of one-dimensional PDEs can be efficiently approximated by finite difference approximation 
schemes.
{\sc Matlab} code~\ref{code:counterfinitediffgbm} 
implements such an approximation scheme in the setting of Proposition~\ref{p:specialcase.gbm} and
{\sc Matlab} code~\ref{code:counterfinitediffabm} 
implements such an approximation scheme in the setting of Proposition~\ref{p:specialcase.abm}.

\lstinputlisting[
	caption={A {\sc Matlab} code to approximate the solution $u$ of \eqref{eq:PDE_quasilinear_1} at $(t_0,x_0)\in [0,T)\times \R$
with a finite difference approximation scheme in the setting of Proposition \ref{p:specialcase.gbm} with $d=1$.},
	label=code:counterfinitediffgbm]
	{Simulation/approximateUfinitediffgbm.m}
	
\lstinputlisting[
	caption={A {\sc Matlab} code to approximate the solution $u$ of \eqref{eq:PDE_quasilinear_1} at $(t_0,x_0)\in [0,T)\times \R$
with a finite difference approximation scheme in the setting of Proposition \ref{p:specialcase.abm} with $d=1$.},
	label=code:counterfinitediffabm]
	{Simulation/approximateUfinitediffabm.m}

 Figures \ref{fig:counter}, \ref{fig:cva}, \ref{fig:borrowlend}, \ref{fig:allencahn}, and the left-hand side
 of Figure \ref{fig:exp} illustrate the
empirical convergence of our scheme.
 In Figures \ref{fig:counter}, \ref{fig:cva}, and \ref{fig:borrowlend}
(respectively \ref{fig:allencahn}) the left-hand side depicts 
for the settings of Subsections \ref{subsec:counterparty}--\ref{subsec:borrowlend} 
(respectively \ref{subsec:allencahn})
in the one-dimensional case 
the relative 
approximation errors 
\begin{equation}\label{eq:relerr}
  \frac{
    \frac{ 1 }{ 10 }
    \sum_{ i = 1 }^{ 10 }
    |
      {\bf U}^{i,[1]}_{ \rho ,  \rho  }(0, x_0 )
      -
      \texttt v
    |
  }
  { 
  |\texttt v|
   }
\end{equation}
against the average 
runtime 
needed to compute the realizations $({\bf U}^{i,[1]}_{ \rho,  \rho  }(0, x_0 ))_{i\in \{ 1, 2, \dots,10\}}$
for 
$ \rho \in \{ 1, 2, \dots, 7 \} $ (respectively $ \rho \in \{ 1, 2, \dots, 5 \} $),
\sloppy{where  $\texttt{v}\in \R$ is the approximation obtained through the finite difference approximation scheme, i.e.,
\texttt{v=approximateUfinitediffgbm(0,x0,2\^{}11)} (respectively \texttt{v=approximateUfinitediffabm(0,x0,2\^{}11)}).}
The left-hand side
 of Figure \ref{fig:exp} shows the relative approximation errors \eqref{eq:relerr} for $ \rho \in \{ 1, 2, \dots, 5 \} $ 
 in the setting of Subsection \ref{subsec:exp},
 where $\texttt{v}=u(0,x_0)$ denotes the value of the exact solution of the PDE. 
These figures are obtained by executing the command
\texttt{ploterrorvsruntime(v,value,time)} (the matrices \texttt{value} and \texttt{time}
are produced in {\sc Matlab} code \ref{code:testrun}). The
 {\sc Matlab} function \texttt{ploterrorvsruntime} is presented in 
{\sc Matlab} code \ref{code:ploterror}.

The right-hand side of the Figures \ref{fig:counter}, \ref{fig:cva}, and \ref{fig:borrowlend}
(respectively \ \ref{fig:allencahn}) depicts for the settings of 
Subsections~\ref{subsec:counterparty}--\ref{subsec:borrowlend}
 (respectively \ref{subsec:allencahn}) in the one hundred-dimensional case with $\rho_{\max}=7$  (respectively $\rho_{\max}=5$) the relative 
approximation increments
\begin{equation}\label{eq:relinc}
  \frac{
    \frac{ 1 }{ 10 }
    \sum_{ i = 1 }^{ 10 }
    |
      {\bf U}^{i,[1]}_{ \rho + 1 , \rho + 1 }(0, x_0 )
      -
      {\bf U}^{i,[1]}_{ \rho, \rho }(0, x_0 ) 
    |
  }
  { 
  \frac{1}{10}
    |\sum_{i=1}^{10}{\bf U}^{i,[1]}_{\rho_{\max},\rho_{\max}  }(0,x_0)|
   }
\end{equation}
\sloppy{against the average 
runtime 
needed to compute the realizations $({\bf U}^{i,[1]}_{ \rho,  \rho}(0, x_0 ))_{i\in \{ 1, 2, \dots,10\}}$
for
$ \rho \in \{ 1, 2, \dots, \rho_{\max}-1  \} $.}
They are obtained by executing the command \texttt{plotincrementvsruntime(value,time)},
where the {\sc Matlab} function \texttt{plotincrementvsruntime} is presented in 
{\sc Matlab} code \ref{code:plotincrements}.

\begin{lstlisting}[
	caption={{\sc Matlab} function to plot relative approximation errors against runtime.},
	label=code:ploterror]
function ploterrorvsruntime(v,value,time)
    merror=mean(abs(value-v))/abs(v);
    mtime=mean(time);
    loglog(mtime,merror,'black-o');
    hold on
    loglog(mtime,1./(mtime).^(1/3)*mtime(1)^(1/3)*merror(1),'black');
    ylabel('relative approximation error');
    xlabel('runtime (seconds)');
    legend('relative approximation error','slope -1/3');
end
\end{lstlisting}
	
\begin{lstlisting}[
	caption={{\sc Matlab} function to plot relative approximation increments against runtime.},
	label=code:plotincrements]
function plotincrementvsruntime(value,time)
    [~,rhomax]=size(value);
    merror=mean(abs(value(:,2:rhomax)-value(:,1:rhomax-1)))/abs(mean(value(:,rhomax)));
    mtime=mean(time(:,1:rhomax-1));
    loglog(mtime,merror,'black-o');
    hold on
    loglog(mtime,1./(mtime).^(1/3)*mtime(1)^(1/3)*merror(1),'black');
    ylabel('relative approximation increments');
    xlabel('runtime (seconds)');
    legend('relative approximation increments','slope -1/3');
end
\end{lstlisting}	

\sloppy{Tables \ref{tab:counterd1}--\ref{tab:expd100} present several statistics for the simulations. More precisely,
Tables \ref{tab:counterd1}--\ref{tab:borrowlendd100} (respectively \ref{tab:allencahnd1}--\ref{tab:expd100}) show for the settings of
Subsections \ref{subsec:counterparty}--\ref{subsec:borrowlend} (respectively \ref{subsec:allencahn}) for all $d\in\{1,100\}$,
$ \rho \in \{ 1, 2, \dots, 7 \} $ (respectively $d\in\{1,100\}$, $ \rho \in \{ 1, 2, \dots, 5 \} $) the average runtime needed to compute 
$({\bf U}^{i,[1]}_{ \rho,  \rho}(0, x_0 ))_{i\in \{ 1, 2, \dots,10\}}$, the empirical mean $ \overline {{\bf U}}^{[1]}_{ \rho ,  \rho  }(0, x_0 )=\frac{ 1 }{ 10 }\sum_{ i = 1 }^{ 10 }{\bf U}^{i,[1]}_{ \rho ,  \rho  }(0, x_0) $, and
the empirical standard deviation 
$
\sqrt{\frac{1}{9}\sum_{ i = 1 }^{ 10 }|{\bf U}^{i,[1]}_{ \rho ,  \rho  }(0, x_0 )-\overline {{\bf U}}^{[1]}_{ \rho ,  \rho  }(0, x_0 )|^2}$.}
  Tables \ref{tab:counterd1}, \ref{tab:cvad1}, \ref{tab:borrowlendd1}, \ref{tab:allencahnd1}, and \ref{tab:expd100} show additionally the relative approximation error \eqref{eq:relerr}. Furthermore, Tables \ref{tab:counterd100}, \ref{tab:cvad100}, \ref{tab:borrowlendd100}, and \ref{tab:allencahnd100} present the relative approximation increments~\eqref{eq:relinc}.

Figures \ref{fig:countermultidim}, \ref{fig:borrowlendmultidim}, and the right-hand side of 
Figure \ref{fig:exp} show the growth of the runtime of our algorithm 
with respect to the dimension
 for each of the example PDEs.
More precisely, Figure \ref{fig:countermultidim} and the left-hand side of Figure \ref{fig:borrowlendmultidim}
show for the settings in Subsections~\ref{subsec:counterparty}--\ref{subsec:borrowlend} the runtime needed to compute 
one realization of ${\bf U}^1_{6,6}(0,x_0)$ against the dimension $d\in \{5,6,\ldots,100\}$. 
The left-hand side of Figure \ref{fig:countermultidim} is obtained by running {\sc Matlab} code
\ref{code:testrunmultidim} in combination with {\sc Matlab} codes \ref{code:approxparameterscounter},
\ref{code:plotruntimevsdim}, and \ref{code:datacounterd100}.
The right-hand side of Figure \ref{fig:countermultidim} is obtained by running {\sc Matlab} code
\ref{code:testrunmultidim} in combination with {\sc Matlab} codes \ref{code:approxparameterscounter},
\ref{code:plotruntimevsdim}, and \ref{code:datacvad100}.
The left-hand side of Figure \ref{fig:borrowlendmultidim} is obtained by running {\sc Matlab} code
\ref{code:testrunmultidim} in combination with {\sc Matlab} codes \ref{code:approxparameterscounter},
\ref{code:plotruntimevsdim}, and \ref{code:databorrowlendd100}.
The right-hand sides of Figures~\ref{fig:borrowlendmultidim} and \ref{fig:exp} show for the the settings in 
Subsections~\ref{subsec:allencahn}--\ref{subsec:exp} the average runtime needed to compute 
20 realizations of ${\bf U}^1_{4,4}(0,x_0)$ against the dimension $d\in \{5,6,\ldots,100\}$. We average over 20 runs here
to obtain smoother results. The right-hand side of Figure \ref{fig:borrowlendmultidim} 
is obtained by running {\sc Matlab} code \ref{code:testrunmultidim} (with line 4 in {\sc Matlab} code \ref{code:testrunmultidim} replaced by \texttt{average=20;} and line 5 in {\sc Matlab} code~\ref{code:testrunmultidim} replaced by \texttt{rhomax=4;})
in combination with {\sc Matlab} codes \ref{code:approxparameterscounter}, 
\ref{code:plotruntimevsdim}, and \ref{code:dataallencahnd100} (with line 10 in {\sc Matlab} code~\ref{code:approxparameterscounter} replaced by \texttt{Mf(rho,k)=rho\^{}k;}).
The right-hand side of Figure \ref{fig:exp} is obtained by running {\sc Matlab} code \ref{code:testrunmultidim} (with line 4 in {\sc Matlab} code \ref{code:testrunmultidim} replaced by \texttt{average=20;} and line 5 in {\sc Matlab} code \ref{code:testrunmultidim} replaced by \texttt{rhomax=4;})
in combination with {\sc Matlab} codes \ref{code:approxparameterscounter}, 
\ref{code:plotruntimevsdim}, and \ref{code:dataexpd100} 
(with line 10 in {\sc Matlab} code \ref{code:approxparameterscounter} replaced by \texttt{Mf(rho,k)=rho\^{}k;}).
\begin{lstlisting}[
	caption={A {\sc Matlab} code to compute 
one realization of ${\bf U}^1_{6,6}(0,x_0)$ for all $d\in \{5,6,\ldots, 100\}$.},
	label=code:testrunmultidim]
global Mf Mg Q c w T dim f g mu sigma eta;
rng(2016)
format long
average=1;
rhomax=6;
dmax=100;
dmin=5;
[T,dim,f,g,eta,mu,sigma]=modelparameters();
[Mf,Mg,Q,c,w,n] = approxparameters(rhomax);
value=zeros(average,dmax-dmin+1);
time=value;
for d=dmin:dmax
  for k=1:average
    dim=d
    tic
    [a, b]=approximateUZgbm(n(rhomax),rhomax,100*ones(d,1),0);
    value(k,d-dmin+1)=a;
    time(k,d-dmin+1)=toc;
  end
end
name = [datestr(now, 'yymmddTHHMMSS') '.mat']; 
save(name,'n','Q','Mf','Mg','value','time')
plotruntimevsdim(dmin, dmax, time)
\end{lstlisting}	
	
\begin{lstlisting}[
	caption={A {\sc Matlab} code to plot the runtime against the dimension. },
	label=code:plotruntimevsdim]
function plotruntimevsdim(dmin, dmax, time)
    mtime=mean(time,1);
    plot((dmin:dmax),mtime,'black')
    hold on
    ylabel('runtime (seconds)')
    xlabel('dimension')
end
\end{lstlisting}

\subsection{Recursive pricing with default risk}\label{subsec:counterparty}
In this subsection we discuss an example which is a special case of the recursive pricing model with default risk due to Duffie, Schroder, \& Skiadas~\cite{DuffieSchroderSkiadas1996}. 
The five-dimensional version of this example has also 
been used as a test example in the literature on numerical approximations of BSDEs
(see, e.g., Bender, Schweizer, \& Zhuo~\cite{BenderSchweizerZhuo2014}).

\begin{figure}[htb]
\begin{center}
\includegraphics[width=0.5\textwidth]{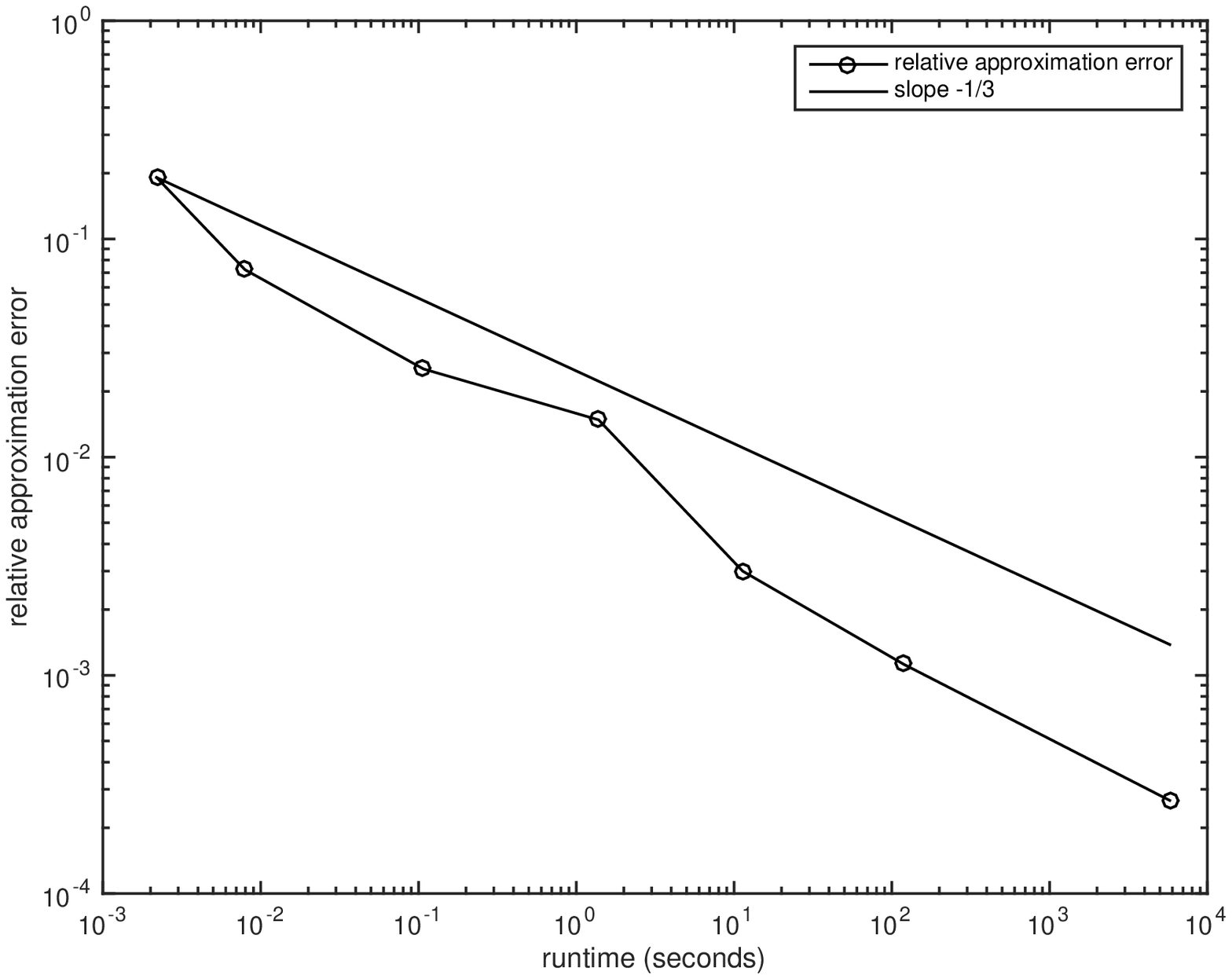}\includegraphics[width=0.5\textwidth]{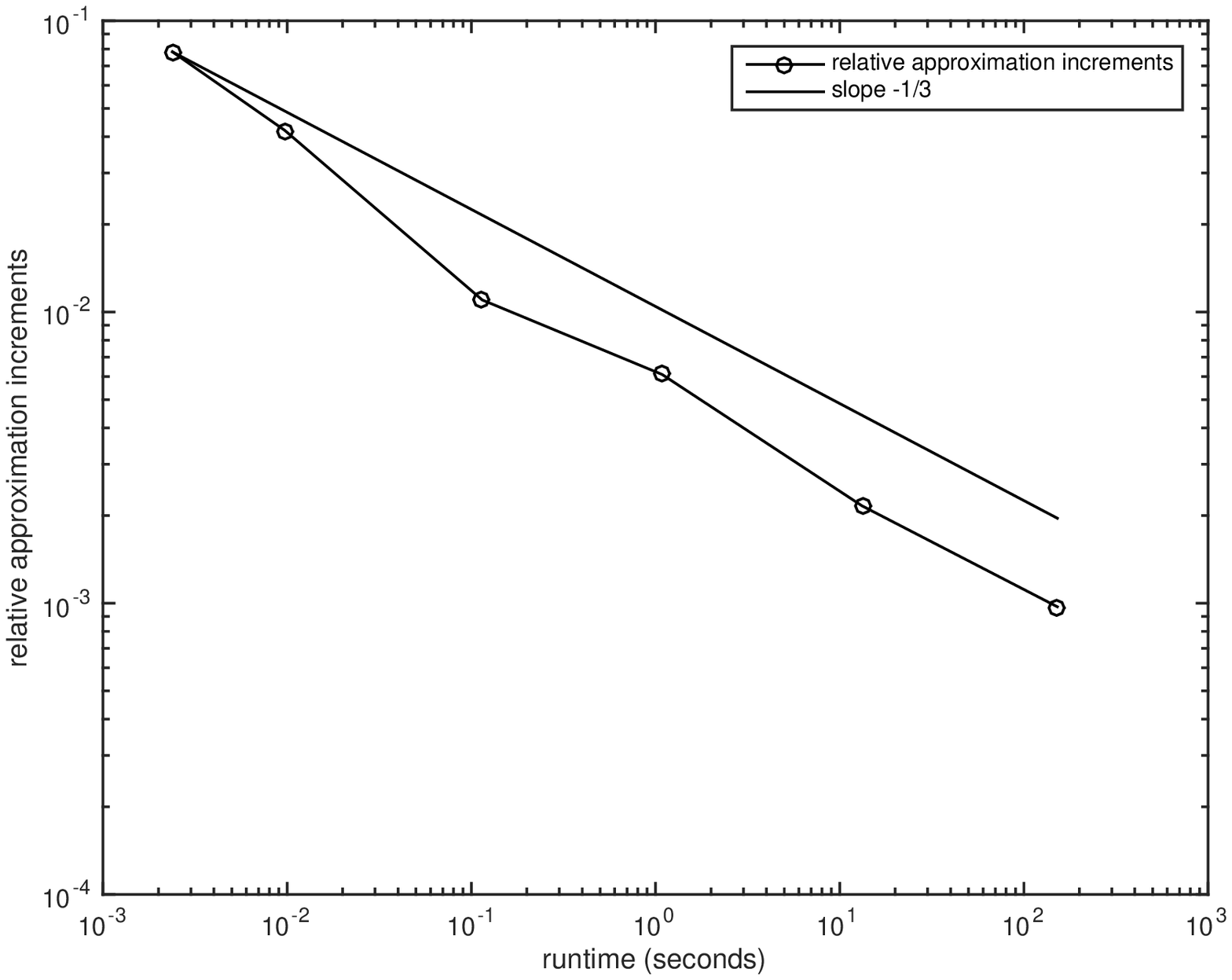}
\end{center}
\caption{Empirical convergence of the scheme \eqref{eq:scheme_Laplace2} for the recursive pricing with default risk example in Subsection \ref{subsec:counterparty}. 
  Left: Relative 
approximation errors 
$
\scriptstyle{
    \frac{ 1 }{ 10  |\texttt v |}
    \sum_{ i = 1 }^{ 10 }
    |
      {\bf U}^{i,[1]}_{ \rho ,  \rho  }(0, x_0 )
      -
      \texttt v
    |
   }
$
for $\rho \in \{1,2,\ldots, 7\}$
against the average 
runtime in the case $d=1$.  Right: Relative 
approximation increments
$\scriptstyle{
  \left(
    \frac{ 1 }{ 10 }
    \sum_{ i = 1 }^{ 10 }
    |
      {\bf U}^{i,[1]}_{ \rho + 1 , \rho + 1 }(0, x_0 )
      -
      {\bf U}^{i,[1]}_{ \rho, \rho }(0, x_0 ) 
    |
  \right)
  \big /
  \left(
  \frac{1}{10}
    |\sum_{i=1}^{10}{\bf U}^{i,[1]}_{7,7}(0,x_0)|
   \right)}
$
for $\rho \in \{1,2,\ldots,6\}$
against the average 
runtime in the case $d=100$.}
\label{fig:counter}
\end{figure}

Throughout this subsection assume the setting in the beginning of Section~\ref{sec:numerics}, let
$\delta=\nicefrac{2}{3}$,
$R=0.02$,
$\gamma^h=0.2$, $\gamma^l=0.02$,
$
  \bar \mu = 0.02
$,
$
  \bar \sigma =0.2
$,
$v^h,v^l \in (0,\infty)$ satisfy $v^h<v^l$,
and assume
for all
$
  s \in [0,T]$, 
   $t \in [s,T]$, $x = ( x_1, \dots, x_d ) \in \R^d$,  
  $y \in \R$, $z \in \R^d$, 
  $k \in \N_0$,  
  $\rho \in \N$, 
  $\theta \in \Theta
$
that
$ T =1 $, 
$
  \eta( x ) = x
$, 
$
  \mu( s, x ) = \bar \mu x
$,
$
  \sigma(s, x) = \bar \sigma \operatorname{diag}(x)
$, 
$
  x_0 = 100 \cdot ( 1, 1, \dots, 1 ) \in \R^{d }
$, 
$
  \mathcal{D}_{ x, s, t }^{ k, \rho, \theta }
  =
  \exp((\bar{\mu}-\tfrac{\bar{\sigma}^2}{2})(t-s))\exp\!\left(\bar{\sigma}\operatorname{diag}(\Delta W^{\theta}_{s,t})\right)
$, 
$
  \mathcal{X}^{ k, \rho, \theta }_{ x, s, t } =
  \mathcal{D}_{ x, s, t }^{ k, \rho, \theta } x
$, 
$
  g(x)
  =
  \min_{
    j \in \{ 1, 2, \dots, d \}
  }
    x_j
$, and 
\begin{equation}
\begin{split}
  f(s,x,y,z) & =
  -
  \left( 1 - \delta \right)
  y 
  \bigg[ 
    \mathbbm{1}_{ ( - \infty, v^h) }( y)
    \,
    \gamma^h
    +
    \mathbbm{1}_{
      [ v^l, \infty )
    }( y )
    \,
    \gamma^l
    +
    \mathbbm{1}_{
      [ v^h, v^l )
    }(y )
    \left[
      \tfrac{ ( \gamma^h - \gamma^l ) }{ ( v^h - v^l ) }
      \left( y - v^h \right)
      +
      \gamma^h
    \right]
  \bigg]
  -
  R y
  .
\end{split}
\end{equation}
Note that the solution $u$ of the PDE \eqref{eq:PDE_quasilinear_1} satisfies for all $t \in [0,T)$, $x=(x_1,x_2,\ldots,x_d)\in \R^d$
that
$
  u(T,x) =
  \min_{
    j \in \{ 1, 2, \dots, d \}
  }
    x_j 
$
and
\begin{multline}  
\label{eq:PDE_counterparty}
  ( \tfrac{ \partial }{ \partial t } u )( t, x )
  +
 \bar \mu \sum_{i=1}^d x_i \big(\tfrac{\partial}{\partial x_i}u\big)(t,x)
  +
  \tfrac{\bar \sigma^2}{ 2 }
  \sum_{i=1}^d
  |x_i|^2\big(\tfrac{\partial^2}{\partial x^2_i}u\big)(t,x)\\
  -
  \left( 1 - \delta \right)
  \min\!\left\{\gamma^h,
  \max\!\left\{\gamma^l,
   \tfrac{ ( \gamma^h - \gamma^l ) }{ ( v^h - v^l ) }
      \left( u(t,x) - v^h \right)
      +
      \gamma^h
      \right\}
      \right\}
      u(t,x)
      -R u(t,x)
  = 0
  .
\end{multline}
In \eqref{eq:PDE_counterparty} the function $u$ models the price of an European financial derivative with payoff $g$ at 
maturity $T$ whose issuer
may default. The number $u(t,x)\in \R$ describes the price of the financial derivative at time $t\in [0,T]$ in dependence on the prices $x=(x_1,\ldots,x_d)\in \R^d$ of the $d$ underlyings of the model given that no default has occurred before time $t$. In the model the arrival intensity of the default and the default-payoff depend on the price of the derivative itself.
 In the case $d=1$ we choose
  $v^h=50$ and $v^l=120$ and in the case $d=100$ we choose $v^h=47$ and $v^l=65$. 
 The thresholds $v^h,v^l \in (0,\infty)$ are adjusted to the dimension 
$d$ since the expectation $\E[g(\mathcal{X}^{ 0, 1, 0 }_{ x_0, 0, T })]$ and the variance
$\Var(g(\mathcal{X}^{ 0, 1, 0 }_{ x_0, 0, T }))$
depend on the dimension (Bender, Schweizer, \& Zhuo~\cite[Subsection 5.3]{BenderSchweizerZhuo2014} choose
  $ d = 5 $, $ v^h = 54 $, and $ v^l = 90 $; all parameters in this subsection except of the parameters $d$, $v^h$, and $v^l$ agree with
Bender, Schweizer, \& Zhuo~\cite[Subsection 5.3]{BenderSchweizerZhuo2014}).
{\sc Matlab} code \ref{code:datacounterd100} presents the 
parameter values in the case $d=100$. In the case $d=1$ line 3
of  {\sc Matlab} code \ref{code:datacounterd100}
is replaced by \texttt{dim=1;} and lines 8 and 9 are changed to \texttt{vh=50;} and \texttt{vl=120;},
respectively. 
The simulation results are shown in Figure \ref{fig:counter}, the left-hand side of Figure \ref{fig:countermultidim}, and Tables \ref{tab:counterd1} 
and \ref{tab:counterd100}.
Figure \ref{fig:counter} suggests an empirical convergence rate close to $\nicefrac{1}{3}$.

\begin{lstlisting}[
	caption={A {\sc Matlab} function that returns the parameter values for the 
	recursive pricing with default risk example.},
	label=code:datacounterd100]
function [T,dim,f,g,eta,mu,sigma] = modelparameters()
    T = 1;
    dim=100;
    sigma=0.2;
    mu=0.02;
    delta=2/3;
    R=0.02;
    vh=47;
    vl=65;
    gammah=0.2;
    gammal=0.02;
    f = @(t,x,y,z) -(1-delta)*y.*((y<vh).*gammah...
        +((vh<=y).*(y<vl)).*((gammah-gammal)/(vh-vl)*(y-vh)+gammah)+(vl<=y)*gammal)-R*y;
    g = @(x) min(x,[],1);
    eta=@(x) x;
end
\end{lstlisting}

\begin{table}[htb]
\begin{center}
\begin{tabular}{|r|rrrrrrr|}\hline 
$\rho$ & 1 & 2 & 3 & 4 & 5 & 6 & 7  \\ \hline 
average runtime in seconds &    0.002 &    0.008 &    0.108 &    1.386 &   11.339 &  119.388 & 5777.365\\ 
$ \overline {{\bf U}}^{[1]}_{ \rho ,  \rho  }(0, x_0 )=\frac{ 1 }{ 10 }\sum_{ i = 1 }^{ 10 }{\bf U}^{i,[1]}_{ \rho ,  \rho  }(0, x_0) $ &   90.807 &   94.345 &   98.138 &   98.697 &   97.712 &   97.749 &   97.703\\ 
$ \sqrt{\frac{1}{9}\sum_{ i = 1 }^{ 10 }|{\bf U}^{i,[1]}_{ \rho ,  \rho  }(0, x_0 )-\overline {{\bf U}}^{[1]}_{ \rho ,  \rho  }(0, x_0 )|^2}$ &   23.618 &    8.818 &    2.966 &    1.474 &    0.386 &    0.158 &    0.035\\ 
 $\frac{ 1 }{ 10 }\sum_{ i = 1 }^{ 10 }\frac{|{\bf U}^{i,[1]}_{ \rho ,  \rho  }(0, x_0 )-\texttt v|}{ |\texttt v|} $ &   0.1912 &   0.0723 &   0.0254 &   0.0148 &   0.0030 &   0.0011 &   0.0003 \\ \hline 
\end{tabular}

\caption{Average runtime, empirical mean, empirical standard deviation, and relative approximation error in the case $d=1$ for the recursive pricing with default risk example in Subsection \ref{subsec:counterparty}. The approximation by the finite difference approximation scheme in {\sc Matlab} code \ref{code:counterfinitediffgbm} yields $\texttt{v=approximateUfinitediffgbm(0,x0,2\^{}11)} \approx 97.705$.}\label{tab:counterd1}
\end{center}
\end{table}

\begin{table}[htb]
\begin{center}
\begin{tabular}{|r|rrrrrrr|}\hline 
$\rho$ & 1 & 2 & 3 & 4 & 5 & 6 & 7  \\ \hline 
average runtime in seconds &    0.002 &    0.010 &    0.115 &    1.085 &   13.209 &  151.868 & 8453.915\\ 
$ \overline {{\bf U}}^{[1]}_{ \rho ,  \rho  }(0, x_0 )=\frac{ 1 }{ 10 }\sum_{ i = 1 }^{ 10 }{\bf U}^{i,[1]}_{ \rho ,  \rho  }(0, x_0) $ &   61.302 &   57.494 &   57.816 &   57.876 &   58.145 &   58.085 &   58.113\\ 
$ \sqrt{\frac{1}{9}\sum_{ i = 1 }^{ 10 }|{\bf U}^{i,[1]}_{ \rho ,  \rho  }(0, x_0 )-\overline {{\bf U}}^{[1]}_{ \rho ,  \rho  }(0, x_0 )|^2}$ &    5.180 &    2.821 &    0.875 &    0.388 &    0.112 &    0.041 &    0.035\\ 
 ${\frac{1}{10} \sum_{ i = 1 }^{ 10 }\frac{|{\bf U}^{i,[1]}_{ \rho+1 ,  \rho+1  }(0, x_0 )-{\bf U}^{i,[1]}_{ \rho ,  \rho  }(0, x_0 )|  }{  |\overline {{\bf U}}^{[1]}_{ 7 ,  7  }(0, x_0 )| }}$ &   0.0782 &   0.0420 &   0.0110 &   0.0061 &   0.0022 &   0.0010 & \\ \hline 
\end{tabular}

\caption{Average runtime, empirical mean, empirical standard deviation, and relative approximation increments in the case $d=100$ for the recursive pricing with default risk example in Subsection \ref{subsec:counterparty}.}\label{tab:counterd100}
\end{center}
\end{table}

\subsection{Pricing with counterparty credit risk}\label{subsec:cva}
In this subsection we present a numerical simulation of a semilinear PDE that 
arises in the valuation of derivative contracts with counterparty credit risk.
The PDE is a special case of the PDEs that are, e.g., derived in Henry-Labord\`ere~\cite{Henry-Labordere2012} and 
Burgard \& Kjaer~\cite{BurgardKjaer2011}.

\begin{figure}[htb]
\begin{center}
\includegraphics[width=0.5\textwidth]{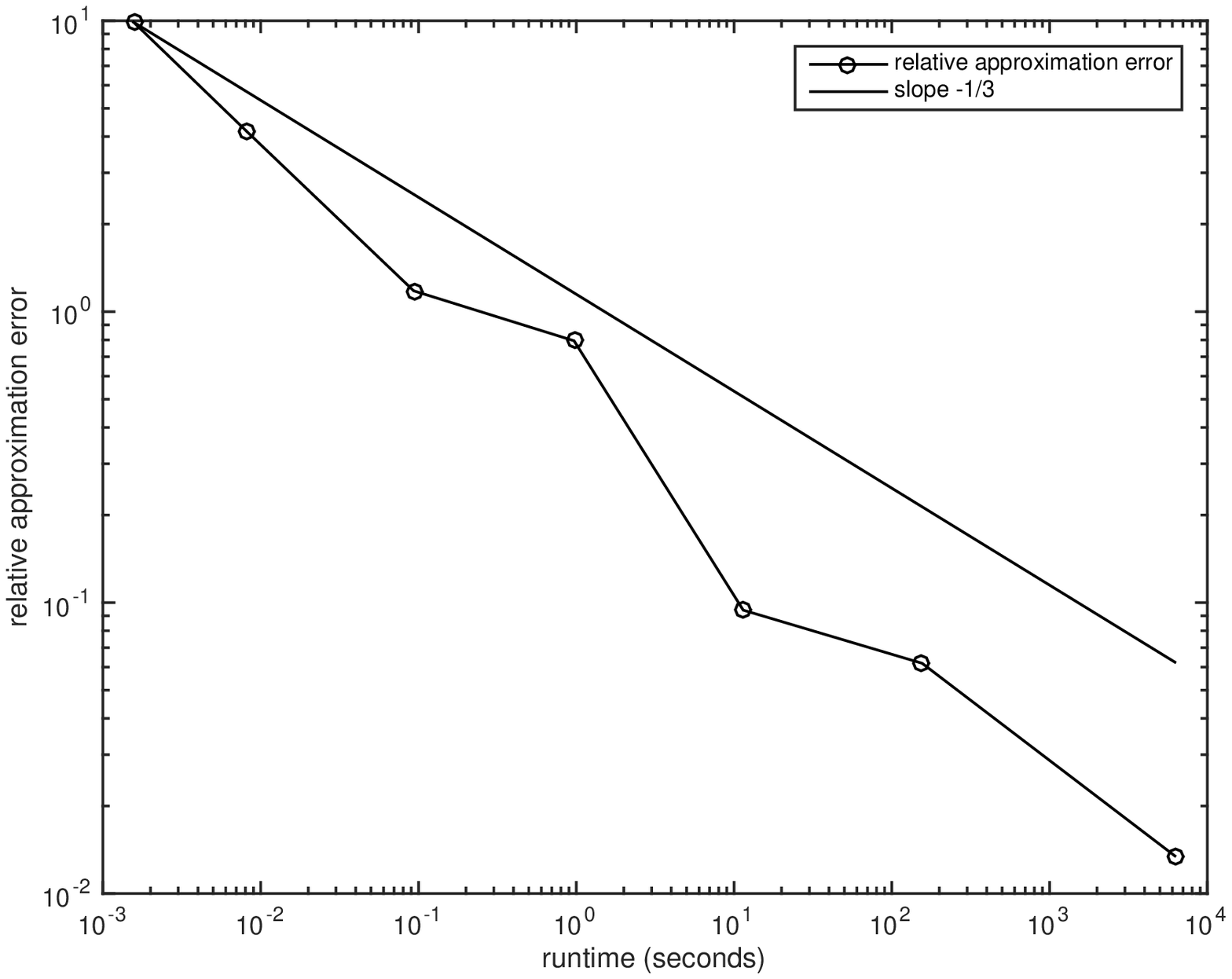}\includegraphics[width=0.5\textwidth]{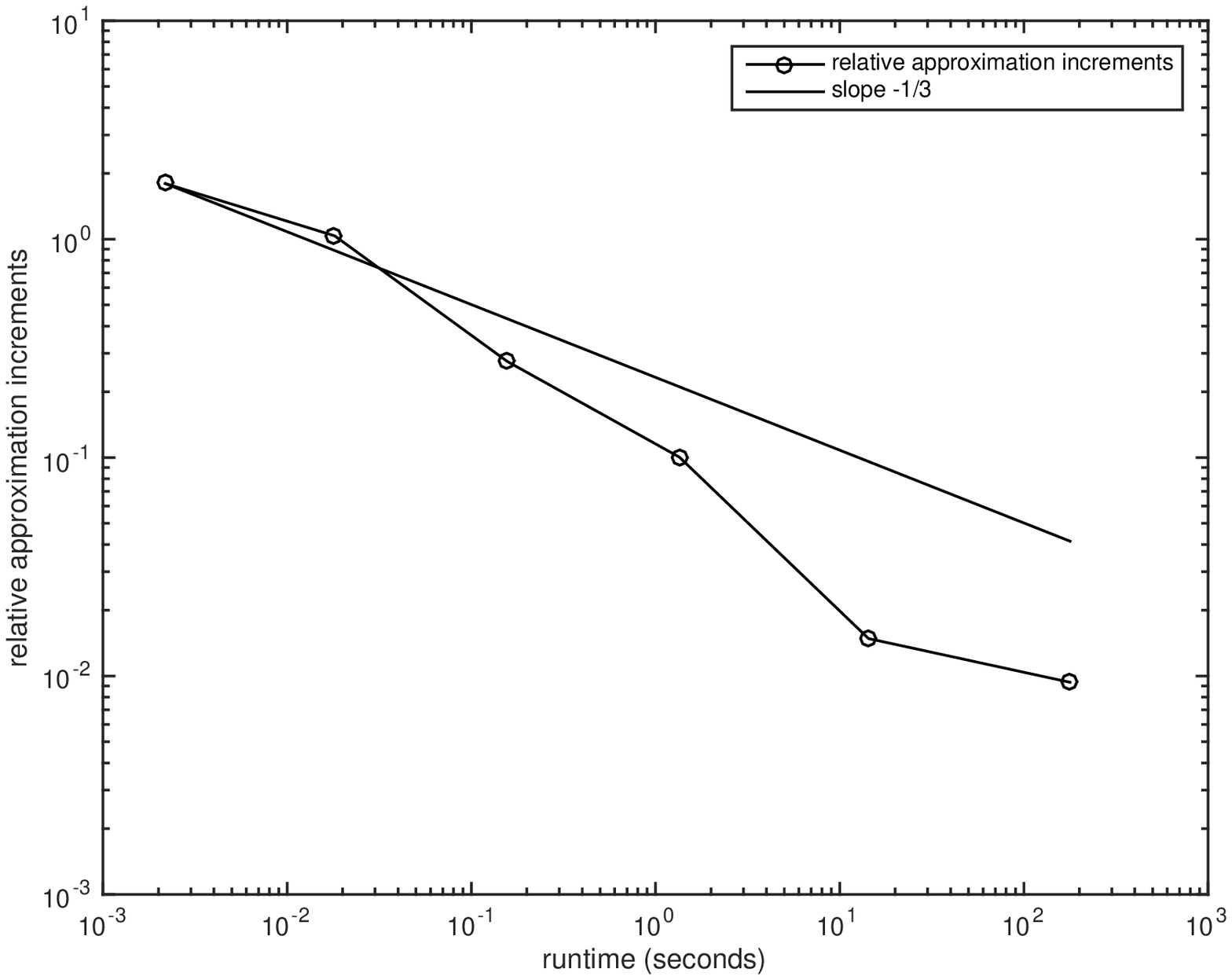}
\end{center}
\caption{Empirical convergence of the scheme \eqref{eq:scheme_Laplace2} for the pricing with counterparty credit risk example in Subsection \ref{subsec:cva}. 
  Left: Relative 
approximation errors 
$
\scriptstyle{
    \frac{ 1 }{ 10  |\texttt v |}
    \sum_{ i = 1 }^{ 10 }
    |
      {\bf U}^{i,[1]}_{ \rho ,  \rho  }(0, x_0 )
      -
      \texttt v
    |
   }
$
for $\rho \in \{1,2,\ldots, 7\}$
against the average 
runtime in the case $d=1$.  Right: Relative 
approximation increments
$\scriptstyle{
  \left(
    \frac{ 1 }{ 10 }
    \sum_{ i = 1 }^{ 10 }
    |
      {\bf U}^{i,[1]}_{ \rho + 1 , \rho + 1 }(0, x_0 )
      -
      {\bf U}^{i,[1]}_{ \rho, \rho }(0, x_0 ) 
    |
  \right)
  \big /
  \left(
  \frac{1}{10}
    |\sum_{i=1}^{10}{\bf U}^{i,[1]}_{7,7}(0,x_0)|
   \right)}
$
for $\rho \in \{1,2,\ldots,6\}$
against the average 
runtime in the case $d=100$.}
\label{fig:cva}
\end{figure}

\begin{figure}[htb]
\begin{center}
\includegraphics[width=0.5\textwidth]{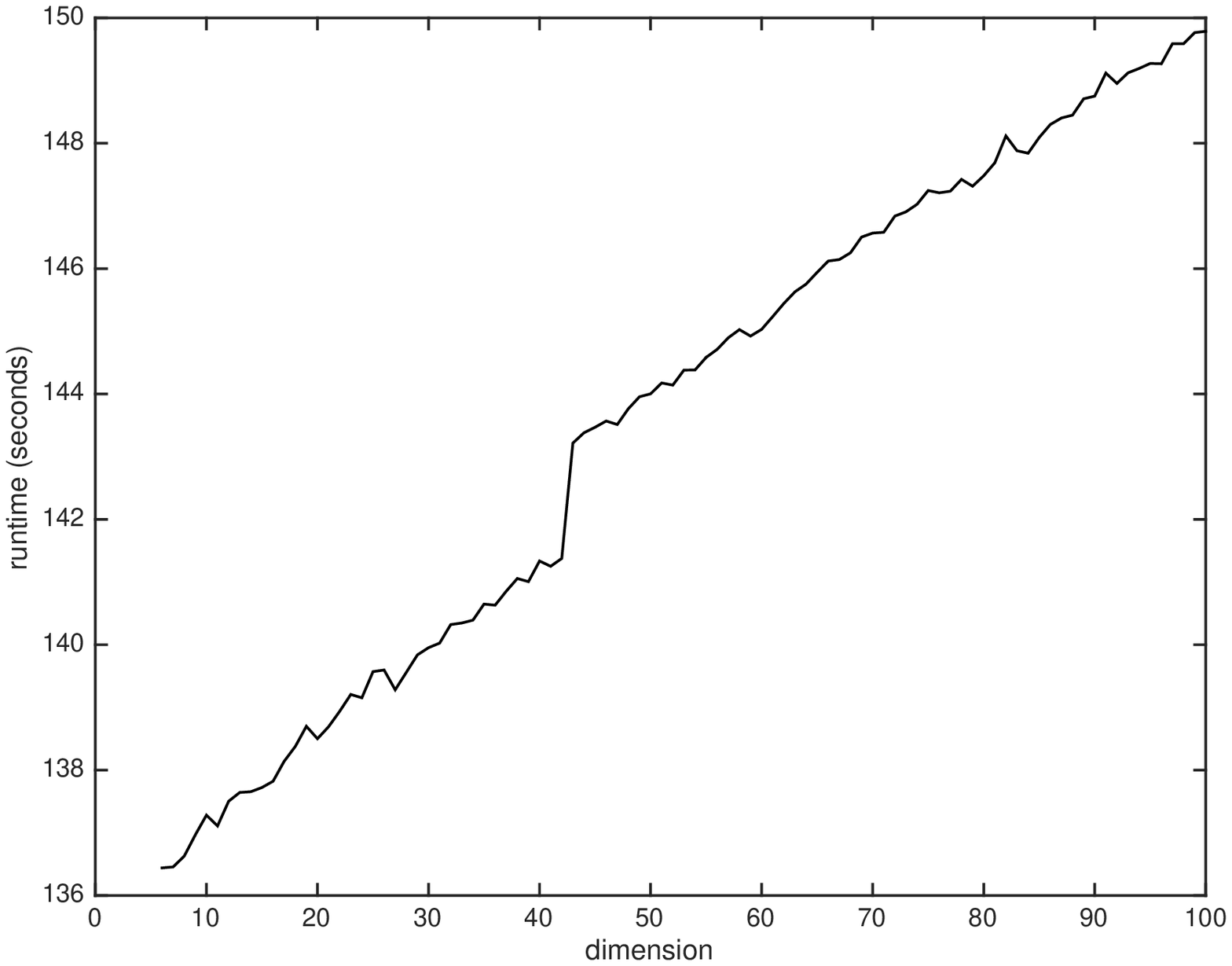}\includegraphics[width=0.5\textwidth]{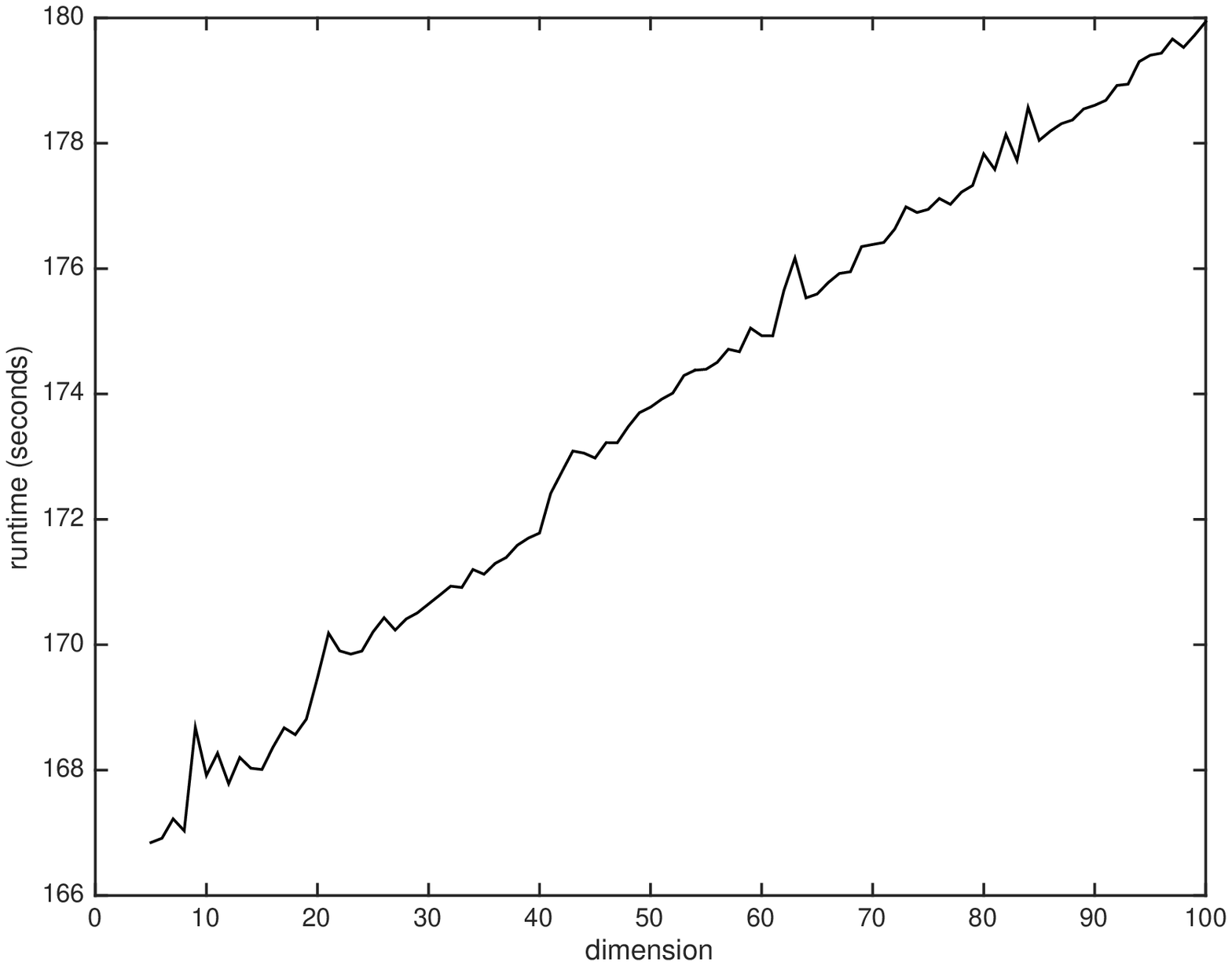}
\end{center}
\caption{Left: Runtime needed to compute 
one realization of ${\bf U}^1_{6,6}(0,x_0)$ against dimension $d\in \{5,6,\ldots,100\}$ for the recursive pricing with default risk example in Subsection \ref{subsec:counterparty}. Right: Runtime needed to compute 
one realization of ${\bf U}^1_{6,6}(0,x_0)$ against dimension $d\in \{5,6,\ldots,100\}$ for the pricing
with counterparty credit risk example in Subsection \ref{subsec:cva}.}
\label{fig:countermultidim}
\end{figure}

Throughout this subsection assume the setting in the beginning of Section~\ref{sec:numerics}, let
$
  \bar \sigma =0.2
$, $\beta=0.03$, $K_1,L\in \R$, $K_2\in (K_1,\infty)$ and assume
for all
$
  s \in [0,T]$, $t \in [s,T]$, $x = ( x_1, \dots, x_d ) \in \R^d$,  
$y \in \R$, $z\in \R^d$, 
$k\in \N_0$, 
  $\rho \in \N$, 
  $\theta \in \Theta
$
that $ T =2 $,
$
  \eta( x ) = x
$, 
$
  \mu( s,x ) = 0
$, 
$
  \sigma(s, x) = \bar \sigma \operatorname{diag}(x)
$, 
$
  x_0 = 100 \cdot ( 1, 1, \dots, 1 ) \in \R^{d }
$, 
$
  \mathcal{D}_{ x, s, t }^{ k, \rho, \theta }
  =
  \exp(-\tfrac{\bar{\sigma}^2}{2}(t-s))\exp\!\left(\bar{\sigma}\operatorname{diag}(\Delta W^{\theta}_{s,t})\right)
$, 
$
  \mathcal{X}^{ k, \rho, \theta }_{ x, s, t } =
  \mathcal{D}_{ x, s, t }^{ k, \rho, \theta } x
$, 
$
f(s,x,y,z) =
  \beta([y]^{+}-y)
 $,
 and
\begin{equation}
\begin{split}
g(x)=\left[\min_{j\in \{1,2,\ldots,d\}}x_j-K_1\right]^{+}-\left[\min_{j\in \{1,2,\ldots,d\}}x_j-K_2\right]^{+}-L.
\end{split}
\end{equation}
Note that the solution $u$ of the PDE \eqref{eq:PDE_quasilinear_1} satisfies for all $t \in [0,T)$, $x=(x_1,x_2,\ldots,x_d)\in \R^d$
that
$
  u(T,x) = \min \{ \max\{\min_{j\in \{1,2,\ldots,d\}}x_j,K_1\},K_2\}-K_1-L
$
and
\begin{equation}  
\label{eq:PDE_cva}
\begin{split}
&
  ( \tfrac{ \partial }{ \partial t } u )( t, x )
  - \beta \min\{u(t,x),0\}
  +
  \tfrac{\bar \sigma^2}{ 2 }
  \sum_{i=1}^d
  |x_i|^2\big(\tfrac{\partial^2}{\partial x^2_i}u\big)(t,x)
  = 0
  .
\end{split}     
\end{equation}
In \eqref{eq:PDE_cva} the function $u$
models the price of an European financial derivative with possibly negative payoff $g$ at 
maturity $T$ whose buyer
may default. The number $u(t,x)\in \R$ describes the price of the financial derivative at time $t\in [0,T]$ in dependence on the prices $x=(x_1,\ldots,x_d)\in \R^d$ of the $d$ underlyings of the model given that no default has occurred before time $t$. In the model the default-payoff depends on the price of the derivative itself.
The choice of the parameters is based on the choice of parameters in
Henry-Labord\`ere~\cite[Subsection 5.3]{Henry-Labordere2012}. 
In the case $d=1$ we choose that 
$K_1=90$, $K_2=110$, and $L=10$. 
In the case $d=100$ we choose that 
$K_1=30$, $K_2=60$, and $L=15$. 
{\sc Matlab} code \ref{code:datacvad100} presents the 
parameter values in the case $d=100$. In the case $d=1$ line 3
of {\sc Matlab} code \ref{code:datacvad100}
is replaced by \texttt{dim=1;} and lines 7, 8, and 9 are changed to \texttt{K1=90;}, \texttt{K2=110;}, and
\texttt{L=10;},
respectively. 
The simulation results are shown in Figure \ref{fig:cva}, the right-hand side of Figure \ref{fig:countermultidim}, and Tables \ref{tab:cvad1} 
and \ref{tab:cvad100}.
Figure \ref{fig:cva} suggests an empirical convergence rate close to $\nicefrac{1}{3}$.
\begin{lstlisting}[
	caption={A {\sc Matlab} function that returns the parameter values for the 
	pricing with counterparty credit risk example.},
	label=code:datacvad100]
function [T,dim,f,g,eta,mu,sigma] = modelparameters()
    T=2;
    dim=100;
    sigma=0.2;
    mu=0;
    beta=0.03;
    K1=30;
    K2=60;
    L=15;
    f = @(t,x,y,z) beta*(subplus(y)-y);
    g = @(x) subplus(min(x,[],1)-K1)-subplus(min(x,[],1)-K2)-L;
    eta=@(x) x;
end
\end{lstlisting}

\begin{table}[htb]
\begin{center}
\begin{tabular}{|r|rrrrrrr|}\hline 
$\rho$ & 1 & 2 & 3 & 4 & 5 & 6 & 7  \\ \hline 
average runtime in seconds &    0.002 &    0.008 &    0.093 &    0.968 &   11.430 &  155.607 & 6226.101\\ 
$ \overline {{\bf U}}^{[1]}_{ \rho ,  \rho  }(0, x_0 )=\frac{ 1 }{ 10 }\sum_{ i = 1 }^{ 10 }{\bf U}^{i,[1]}_{ \rho ,  \rho  }(0, x_0) $ &   -0.582 &   -3.614 &   -0.767 &   -0.433 &   -0.916 &   -0.866 &   -0.884\\ 
$ \sqrt{\frac{1}{9}\sum_{ i = 1 }^{ 10 }|{\bf U}^{i,[1]}_{ \rho ,  \rho  }(0, x_0 )-\overline {{\bf U}}^{[1]}_{ \rho ,  \rho  }(0, x_0 )|^2}$ &    9.346 &    3.964 &    1.279 &    0.782 &    0.105 &    0.067 &    0.015\\ 
 $\frac{ 1 }{ 10 }\sum_{ i = 1 }^{ 10 }\frac{|{\bf U}^{i,[1]}_{ \rho ,  \rho  }(0, x_0 )-\texttt v|}{ |\texttt v|} $ &   9.8923 &   4.1417 &   1.1794 &   0.7941 &   0.0943 &   0.0617 &   0.0135 \\ \hline 
\end{tabular}

\caption{Average runtime, empirical mean, empirical standard deviation, and relative approximation error in the case $d=1$ for the pricing with counterparty credit risk example in Subsection \ref{subsec:cva}. The approximation by the finite difference approximation scheme in {\sc Matlab} code \ref{code:counterfinitediffgbm} yields $\texttt{v=approximateUfinitediffgbm(0,x0,2\^{}11)} \approx -0.883$.}\label{tab:cvad1}
\end{center}
\end{table}

\begin{table}[htb]
\begin{center}
\begin{tabular}{|r|rrrrrrr|}\hline 
$\rho$ & 1 & 2 & 3 & 4 & 5 & 6 & 7  \\ \hline 
average runtime in seconds &    0.002 &    0.018 &    0.156 &    1.348 &   14.332 &  177.989 & 8935.848\\ 
$ \overline {{\bf U}}^{[1]}_{ \rho ,  \rho  }(0, x_0 )=\frac{ 1 }{ 10 }\sum_{ i = 1 }^{ 10 }{\bf U}^{i,[1]}_{ \rho ,  \rho  }(0, x_0) $ &    5.823 &    1.878 &    2.376 &    2.450 &    2.607 &    2.617 &    2.626\\ 
$ \sqrt{\frac{1}{9}\sum_{ i = 1 }^{ 10 }|{\bf U}^{i,[1]}_{ \rho ,  \rho  }(0, x_0 )-\overline {{\bf U}}^{[1]}_{ \rho ,  \rho  }(0, x_0 )|^2}$ &    5.741 &    3.051 &    1.041 &    0.335 &    0.053 &    0.027 &    0.008\\ 
 ${\frac{1}{10} \sum_{ i = 1 }^{ 10 }\frac{|{\bf U}^{i,[1]}_{ \rho+1 ,  \rho+1  }(0, x_0 )-{\bf U}^{i,[1]}_{ \rho ,  \rho  }(0, x_0 )|  }{  |\overline {{\bf U}}^{[1]}_{ 7 ,  7  }(0, x_0 )| }}$ &   1.7971 &   1.0354 &   0.2758 &   0.1004 &   0.0148 &   0.0094 & \\ \hline 
\end{tabular}

\caption{Average runtime, empirical mean, empirical standard deviation, and relative approximation increments in the case $d=100$ for the pricing with counterparty credit risk example in Subsection \ref{subsec:cva}.}\label{tab:cvad100}
\end{center}
\end{table}

\subsection{Pricing with different interest rates for borrowing and lending}\label{subsec:borrowlend}
We consider a pricing problem of an European option in a financial market with different interest rates for borrowing and lending. The model goes back to Bergman~\cite{Bergman1995} and serves as a standard example in the literature on numerical methods for BSDEs (see, e.g., \cite{GobetLemorWarin2005,BenderDenk2007,BenderSchweizerZhuo2014,BriandLabart2014,CrisanManolarakis2012}).

\begin{figure}[htb]
\begin{center}
\includegraphics[width=0.5\textwidth]{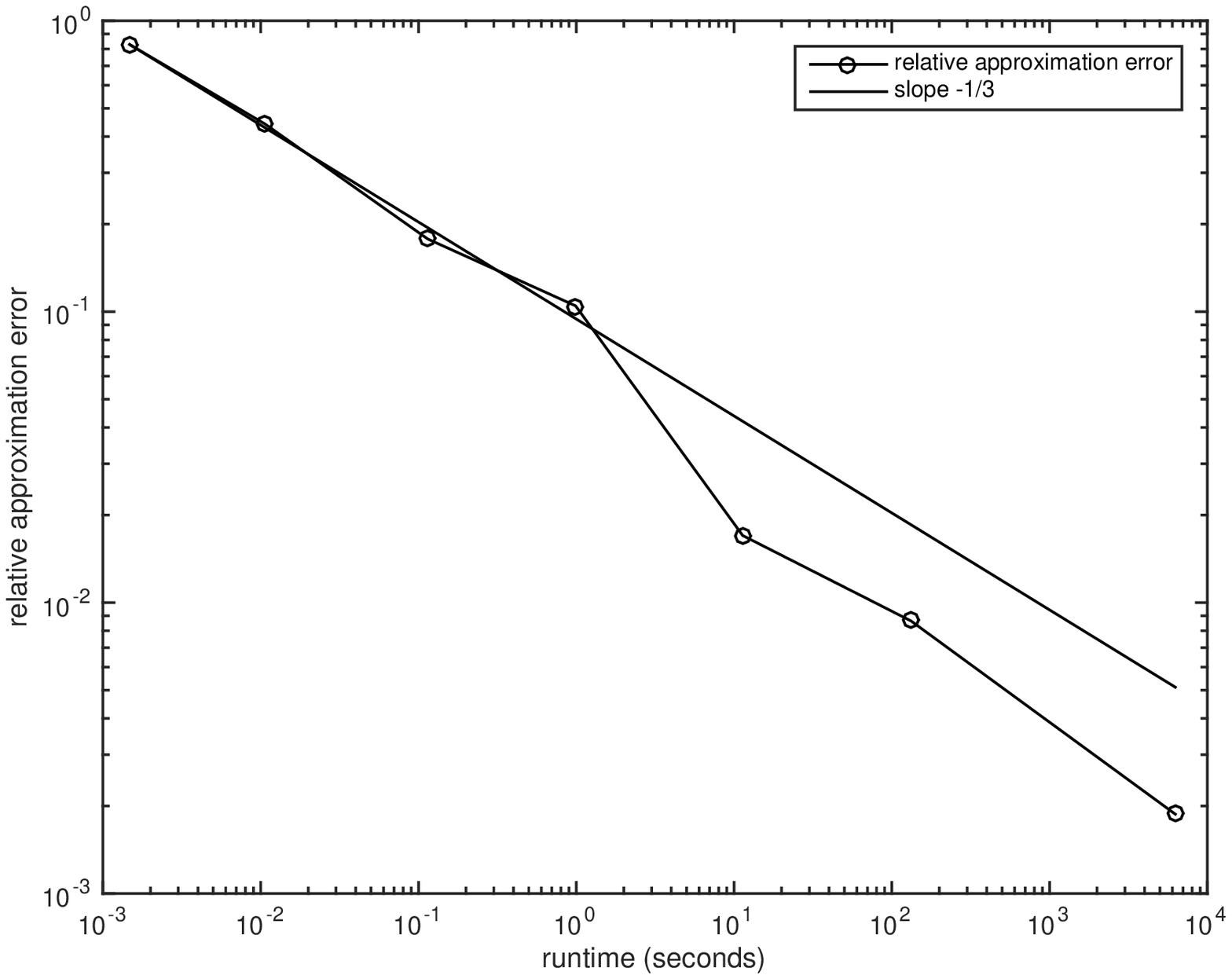}\includegraphics[width=0.5\textwidth]{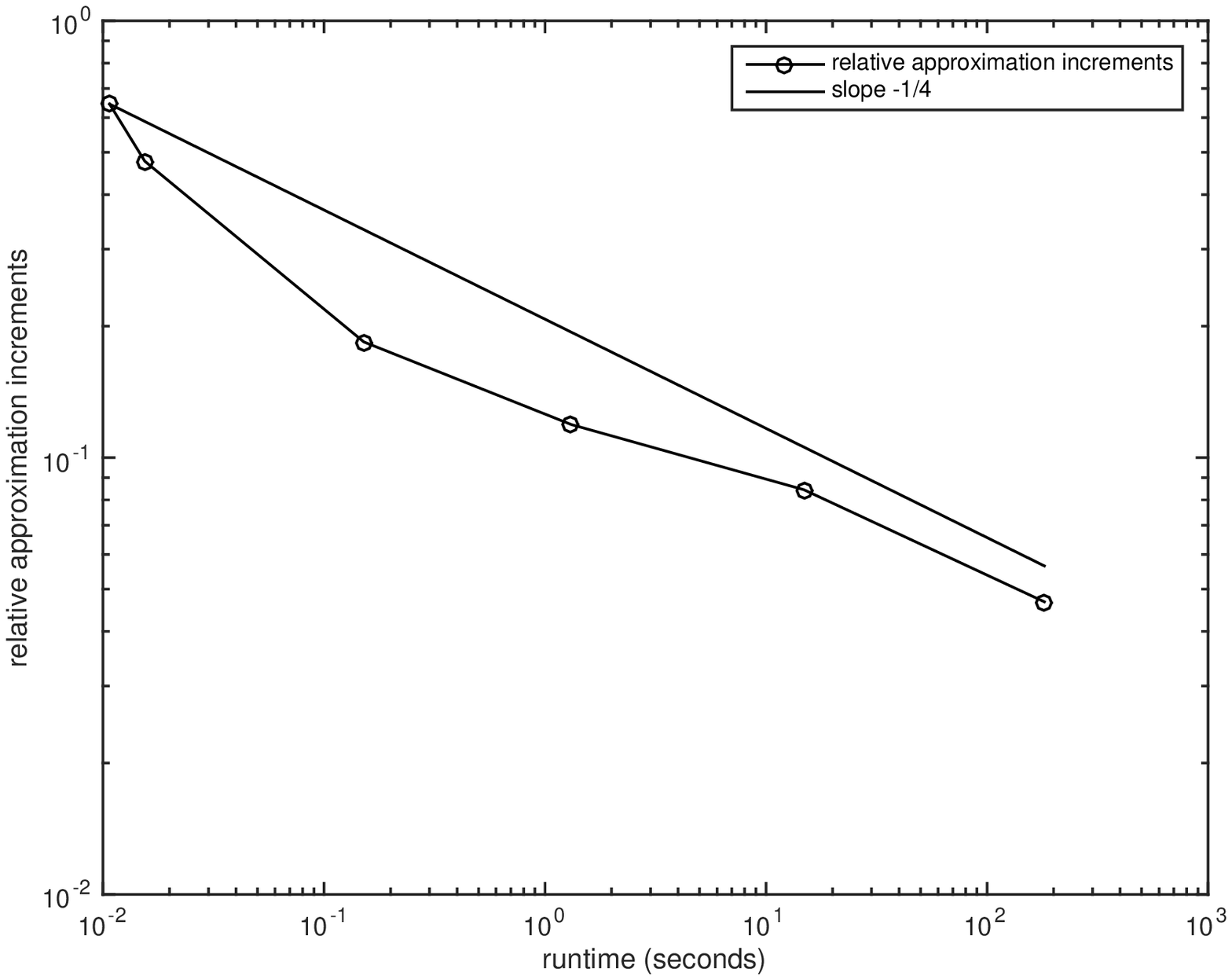}
\end{center}
\caption{Empirical convergence of the scheme \eqref{eq:scheme_Laplace2} for the pricing with different interest rates example in Subsection \ref{subsec:borrowlend}. 
  Left: Relative 
approximation errors  
$
\scriptstyle{
    \frac{ 1 }{ 10  |\texttt v |}
    \sum_{ i = 1 }^{ 10 }
    |
      {\bf U}^{i,[1]}_{ \rho ,  \rho  }(0, x_0 )
      -
      \texttt v
    |
   }
$
for $\rho \in \{1,2,\ldots, 7\}$
against the average 
runtime in the case $d=1$.  Right: Relative 
approximation increments
$\scriptstyle{
  \left(
    \frac{ 1 }{ 10 }
    \sum_{ i = 1 }^{ 10 }
    |
      {\bf U}^{i,[1]}_{ \rho + 1 , \rho + 1 }(0, x_0 )
      -
      {\bf U}^{i,[1]}_{ \rho, \rho }(0, x_0 ) 
    |
  \right)
  \big /
  \left(
  \frac{1}{10}
    |\sum_{i=1}^{10}{\bf U}^{i,[1]}_{7,7}(0,x_0)|
   \right)}
$
for $\rho \in \{1,2,\ldots,6\}$
against the average 
runtime in the case $d=100$.}
\label{fig:borrowlend}
\end{figure}

Throughout this subsection assume the setting in the beginning of Section~\ref{sec:numerics}, 
let $\bar \mu =0.06$, $\bar \sigma=0.2$, $R^l=0.04$, $R^b=0.06$, and assume 
for all
$
  s \in [0,T]$, $t \in [s,T]$, $x = ( x_1, \dots, x_d ) \in \R^d$,  
  $y\in \R$, $z=( z_1, \dots, z_d )\in \R^d$, 
  $k \in \N_0$, 
  $\rho \in \N$, 
  $\theta \in \Theta
$
that $T=0.5$,
$
  \eta( x ) = x
$, 
$
  \mu( s,x ) = \bar \mu x
$, 
$
  \sigma(s,x) = \bar \sigma \operatorname{diag}(x)
$, 
$
  x_0 = 100 \cdot ( 1, 1, \dots, 1 ) \in \R^{d }
$, 
$
  \mathcal{D}_{ x, s, t }^{ k, \rho, \theta }
  =
  \exp((\bar{\mu}-\tfrac{\bar{\sigma}^2}{2})(t-s))\exp\!\left(\bar{\sigma}\operatorname{diag}(\Delta W^{\theta}_{s,t})\right)
$, 
$
  \mathcal{X}^{ k, \rho, \theta }_{ x, s, t } =
  \mathcal{D}_{ x, s, t }^{ k, \rho, \theta } x
$, and
\begin{equation}
f(s,x,y,z)=-R^ly-\frac{(\bar \mu-R^l)}{\bar \sigma}\left(\sum_{i=1}^dz_i\right)+(R^b-R^l)\left[\frac{1}{\bar \sigma}\left(\sum_{i=1}^dz_i\right)-y\right]^{+}.
\end{equation}
Note that the solution $u$ of the PDE \eqref{eq:PDE_quasilinear_1} satisfies for all $t \in [0,T)$, $x=(x_1,x_2,\ldots,x_d)\in \R^d$
that
$
  u(T,x) = g(x) 
$
and
\begin{multline}  
\label{eq:PDE_borrowlend}
  ( \tfrac{ \partial }{ \partial t } u )( t, x )
  +
  \tfrac{\bar \sigma^2}{ 2 }
  \sum_{i=1}^d
  |x_i|^2\big(\tfrac{\partial^2}{\partial x^2_i}u\big)(t,x)
  \\
  -\min\!\bigg\{R^b\bigg(u(t,x)-\sum_{i=1}^d x_i \big(\tfrac{\partial}{\partial x_i}u\big)(t,x)\bigg),
  R^l\bigg(u(t,x)-\sum_{i=1}^d x_i \big(\tfrac{\partial}{\partial x_i}u\big)(t,x)\bigg)\bigg\}
  = 0
  .  
\end{multline}
In \eqref{eq:PDE_borrowlend} the function $u$ 
 models the price of an European financial derivative with payoff $g$ at 
maturity $T$ in a financial market with a higher interest rate for borrowing than for lending. The number $u(t,x)\in \R$ describes the price of the financial derivative at time $t\in [0,T]$ in dependence on the prices $x=(x_1,\ldots,x_d)\in \R^d$ of the $d$ underlyings of the model.
In the case $d=1$, assume that
for all $x \in \R$ it holds that $g(x)=[x-100]^{+}$.
This setting agrees with the setting in Gobet, Lemor, \& Warin~\cite[Subsection 6.3.1]{GobetLemorWarin2005}, where it
is also noted that the PDE \eqref{eq:PDE_borrowlend} is actually linear. More precisely, 
$u$ also satisfies \eqref{eq:PDE_quasilinear_1} for all $t\in [0,T)$, $x\in \R$ with $f$ being replaced by 
$\bar f\colon [0,T]\times \R \times \R \times \R \to \R$ satisfying for all
$t\in [0,T]$, $x,y,z\in \R$ that $\bar f(t,x,y,z)=-R^by-\frac{(\bar \mu-R^b)}{\bar \sigma}z$.
In the case $d=100$, assume that 
for all $ x = ( x_1, \dots, x_d ) \in \R^d$ it holds that $g(x)=[\max_{i\in \{1,\ldots, 100\}}x_i-120]^+-2[\max_{i\in \{1,\ldots, 100\}}x_i-150]^+$.
This choice of $g$ is based on the choice of $g$ in Bender, Schweizer, \& Zhuo~\cite[Subsection 4.2]{BenderSchweizerZhuo2014}. We note that with this choice of $g$ the PDE \eqref{eq:PDE_borrowlend} can not be reduced to a linear PDE.
{\sc Matlab} code \ref{code:databorrowlendd100} presents the 
parameter values in the case $d=100$. In the case $d=1$ line 3
of {\sc Matlab} code \ref{code:databorrowlendd100}
is replaced by \texttt{dim=1;} and line 9 of {\sc Matlab} code \ref{code:databorrowlendd100}
 is replaced by 
\texttt{ g=@(x) subplus(x-100);}.
The simulation results are shown in Figure \ref{fig:borrowlend}, the left-hand side of Figure \ref{fig:borrowlendmultidim}, and Tables \ref{tab:borrowlendd1} 
and \ref{tab:borrowlendd100}.
The left-hand side of Figure \ref{fig:borrowlend} suggests an empirical convergence rate close to $\nicefrac{1}{3}$ in the case $d=1$. Moreover, the right-hand side of Figure \ref{fig:borrowlend} suggests an empirical convergence rate close to
$\nicefrac{1}{4}$ in the case $d=100$.
\begin{lstlisting}[
	caption={A {\sc Matlab} function that returns the parameter values for the 
	pricing with different interest rates example.},
	label=code:databorrowlendd100]
function [T,dim,f,g,eta,mu,sigma] = modelparameters()
    T = 0.5;
    dim=100;
    Rl=0.04;
    Rb=0.06;
    mu=0.06;
    sigma=0.2;
    f = @(t,x,y,z) -Rl.*y-(mu-Rl)/sigma*sum(z,1)+(Rb-Rl)*subplus(1/sigma*sum(z,1)-y); 
    g= @(x) subplus(max(x,[],1)-120)-2*subplus(max(x,[],1)-150); 
    eta=@(x) x;
end
\end{lstlisting}

\begin{table}[htb]
\begin{center}
\begin{tabular}{|r|rrrrrrr|}\hline 
$\rho$ & 1 & 2 & 3 & 4 & 5 & 6 & 7  \\ \hline 
average runtime in seconds &    0.001 &    0.011 &    0.114 &    0.990 &   11.347 &  132.484 & 6264.475\\ 
$ \overline {{\bf U}}^{[1]}_{ \rho ,  \rho  }(0, x_0 )=\frac{ 1 }{ 10 }\sum_{ i = 1 }^{ 10 }{\bf U}^{i,[1]}_{ \rho ,  \rho  }(0, x_0) $ &    5.695 &    5.947 &    7.085 &    7.631 &    7.156 &    7.162 &    7.166\\ 
$ \sqrt{\frac{1}{9}\sum_{ i = 1 }^{ 10 }|{\bf U}^{i,[1]}_{ \rho ,  \rho  }(0, x_0 )-\overline {{\bf U}}^{[1]}_{ \rho ,  \rho  }(0, x_0 )|^2}$ &    7.780 &    4.080 &    1.612 &    0.811 &    0.151 &    0.071 &    0.016\\ 
 $\frac{ 1 }{ 10 }\sum_{ i = 1 }^{ 10 }\frac{|{\bf U}^{i,[1]}_{ \rho ,  \rho  }(0, x_0 )-\texttt v|}{|\texttt v|} $ &   0.8285 &   0.4417 &   0.1777 &   0.1047 &   0.0170 &   0.0086 &   0.0019 \\ \hline 
\end{tabular}

\caption{Average runtime, empirical mean, empirical standard deviation, and relative approximation error in the case $d=1$ for the pricing with different interest rates example in Subsection \ref{subsec:borrowlend}. The approximation by the finite difference approximation scheme in {\sc Matlab} code \ref{code:counterfinitediffgbm} yields $\texttt{v=approximateUfinitediffgbm(0,x0,2\^{}11)} \approx 7.156$.}\label{tab:borrowlendd1}
\end{center}
\end{table}

\begin{table}[htb]
\begin{center}
\begin{tabular}{|r|rrrrrrr|}\hline 
$\rho$ & 1 & 2 & 3 & 4 & 5 & 6 & 7  \\ \hline 
average runtime in seconds &    0.011 &    0.015 &    0.151 &    1.317 &   14.813 &  181.647 & 8825.390\\ 
$ \overline {{\bf U}}^{[1]}_{ \rho ,  \rho  }(0, x_0 )=\frac{ 1 }{ 10 }\sum_{ i = 1 }^{ 10 }{\bf U}^{i,[1]}_{ \rho ,  \rho  }(0, x_0) $ &   28.902 &   22.854 &   23.356 &   21.771 &   21.374 &   21.274 &   21.299\\ 
$ \sqrt{\frac{1}{9}\sum_{ i = 1 }^{ 10 }|{\bf U}^{i,[1]}_{ \rho ,  \rho  }(0, x_0 )-\overline {{\bf U}}^{[1]}_{ \rho ,  \rho  }(0, x_0 )|^2}$ &    8.798 &   11.317 &    4.492 &    2.953 &    1.449 &    1.376 &    0.467\\ 
 ${\frac{1}{10} \sum_{ i = 1 }^{ 10 }\frac{|{\bf U}^{i,[1]}_{ \rho+1 ,  \rho+1  }(0, x_0 )-{\bf U}^{i,[1]}_{ \rho ,  \rho  }(0, x_0 )|  }{  |\overline {{\bf U}}^{[1]}_{ 7 ,  7  }(0, x_0 )| }}$ &   0.6446 &   0.4770 &   0.1840 &   0.1190 &   0.0844 &   0.0467 & \\ \hline 
\end{tabular}

\caption{Average runtime, empirical mean, empirical standard deviation, and relative approximation increments in the case $d=100$ for the pricing with different interest rates example in Subsection \ref{subsec:borrowlend}.}\label{tab:borrowlendd100}
\end{center}
\end{table}

\subsection{Allen-Cahn equation}\label{subsec:allencahn}
In this subsection we consider the Allen-Cahn equation with a double well potential.

\begin{figure}[htb]
\begin{center}
\includegraphics[width=0.5\textwidth]{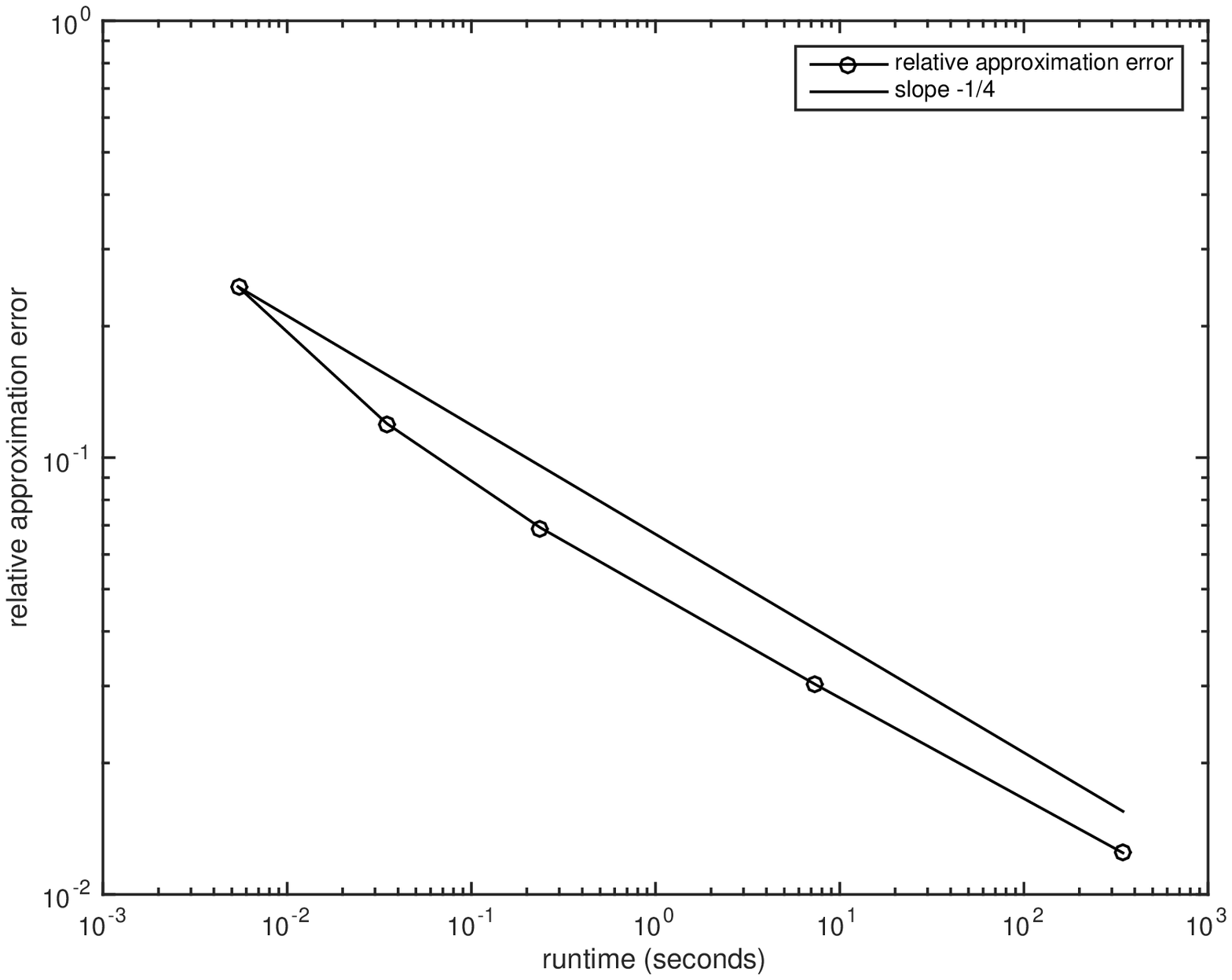}\includegraphics[width=0.5\textwidth]{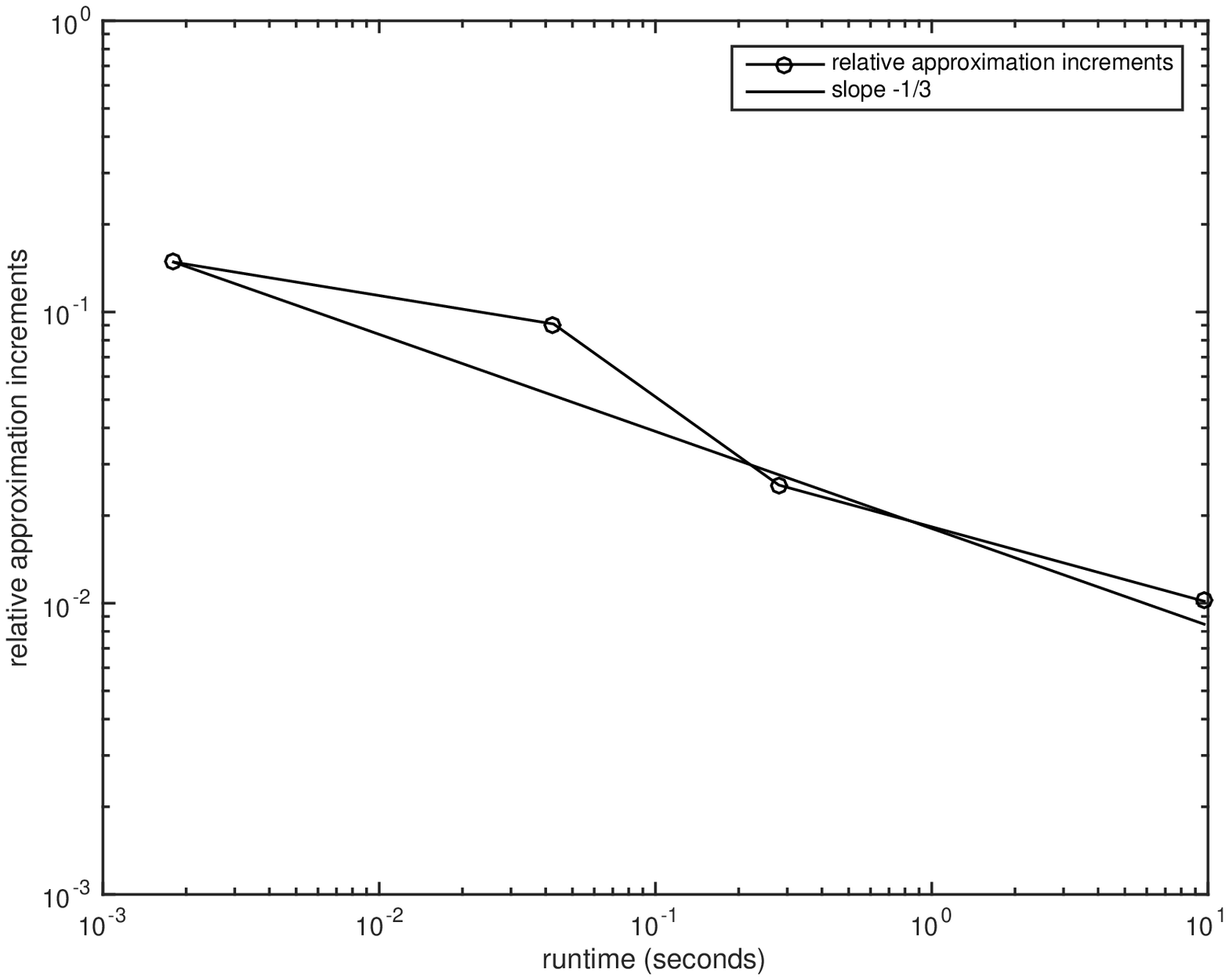}
\end{center}
\caption{Empirical convergence of the scheme \eqref{eq:scheme_Laplace} for the Allen-Cahn equation in Subsection \ref{subsec:allencahn}. 
  Left: Relative 
approximation errors 
$
\scriptstyle{
    \frac{ 1 }{ 10  |\texttt v |}
    \sum_{ i = 1 }^{ 10 }
    |
      {\bf U}^{i,[1]}_{ \rho ,  \rho  }(0, x_0 )
      -
      \texttt v
    |
   }
$
for $\rho \in \{1,2,\ldots, 5\}$
against the average 
runtime in the case $d=1$. Right: Relative 
approximation increments
$\scriptstyle{
  \left(
    \frac{ 1 }{ 10 }
    \sum_{ i = 1 }^{ 10 }
    |
      {\bf U}^{i,[1]}_{ \rho + 1 , \rho + 1 }(0, x_0 )
      -
      {\bf U}^{i,[1]}_{ \rho, \rho }(0, x_0 ) 
    |
  \right)
  \big /
  \left(
  \frac{1}{10}
    |\sum_{i=1}^{10}{\bf U}^{i,[1]}_{5,5}(0,x_0)|
   \right)}
$
for $\rho \in \{1,2,3,4\}$
against the average 
runtime in the case $d=100$.}
\label{fig:allencahn}
\end{figure}

\begin{figure}[htb]
\begin{center}
\includegraphics[width=0.5\textwidth]{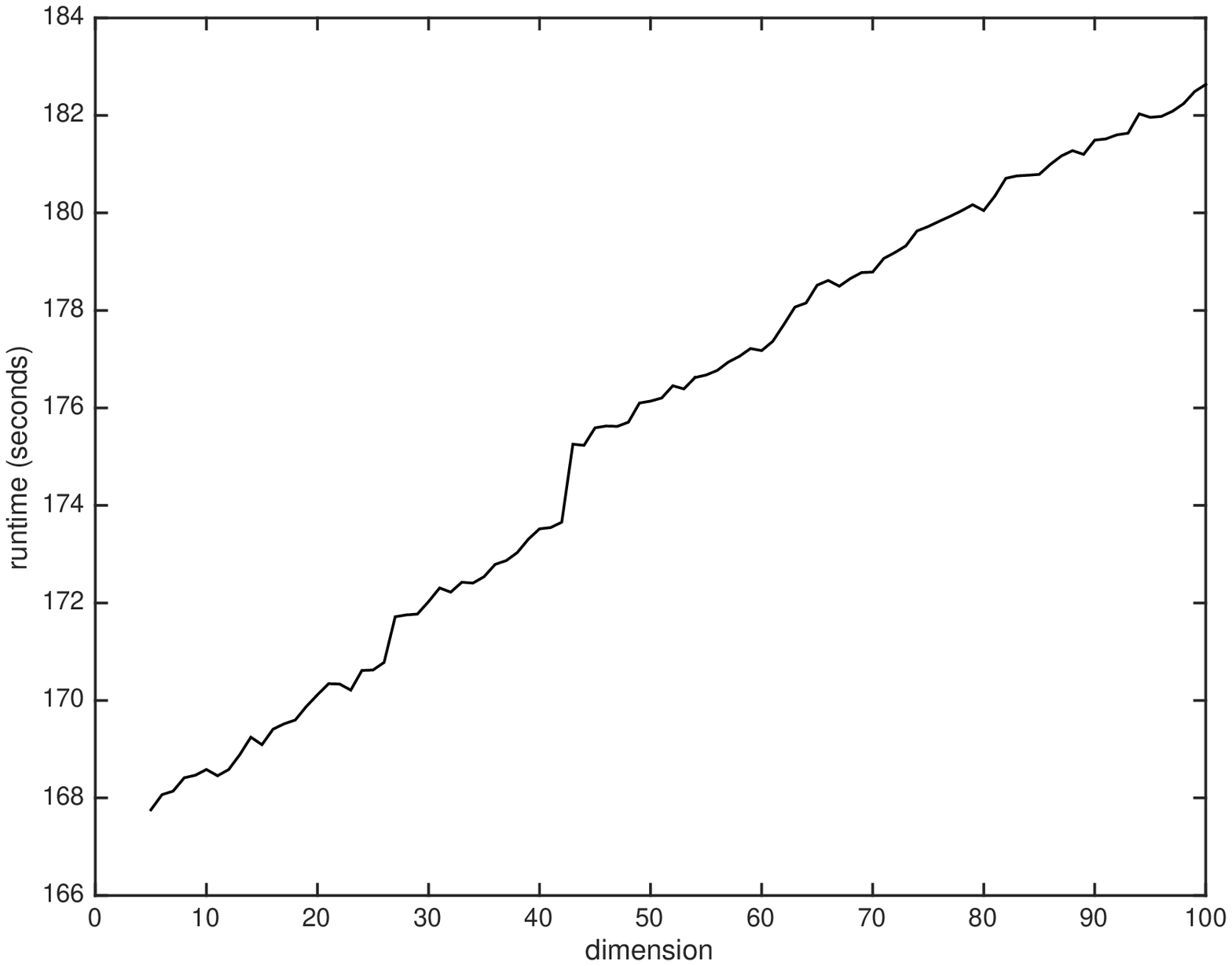}\includegraphics[width=0.5\textwidth]{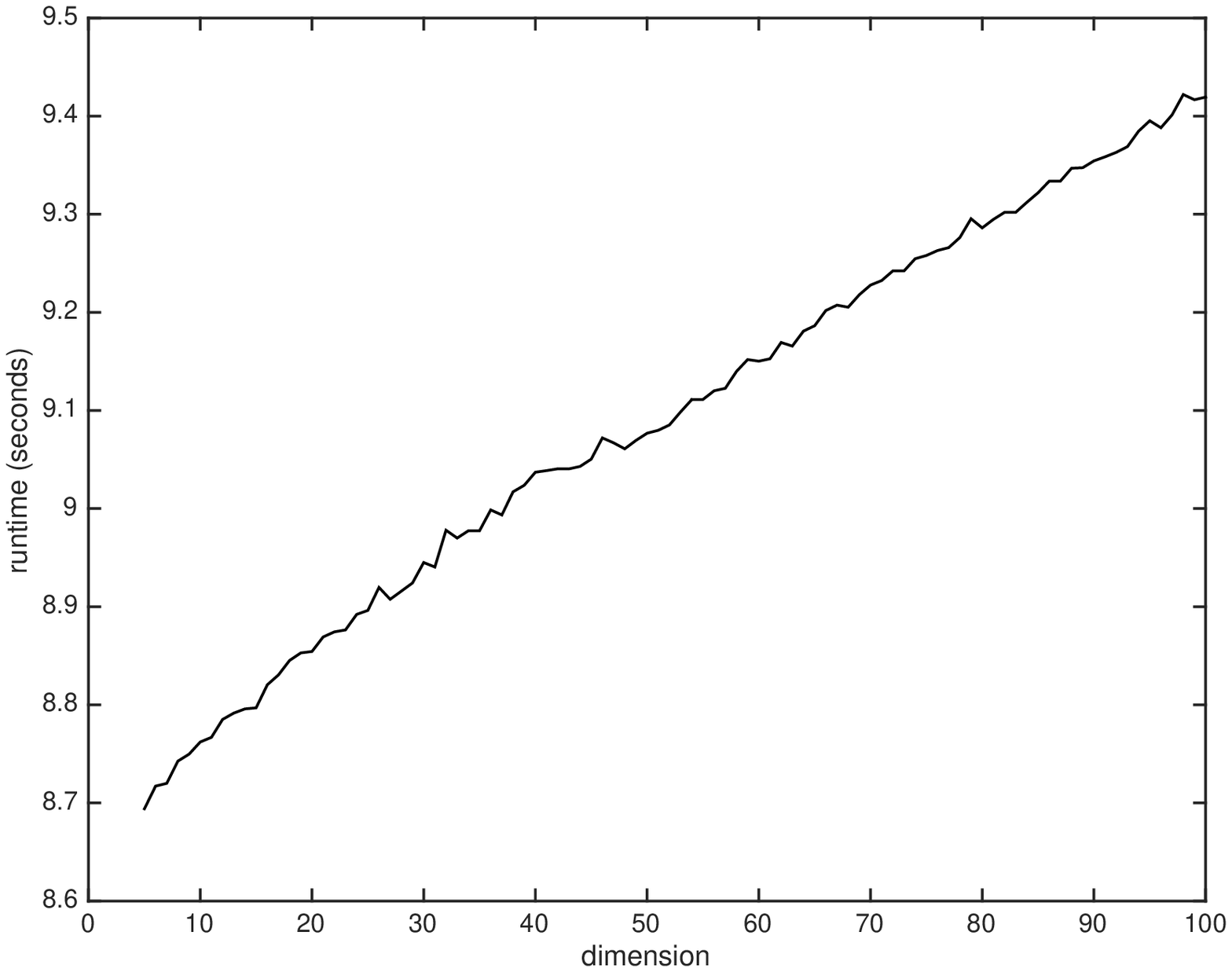}
\end{center}
\caption{Left: Runtime needed to compute 
one realization of ${\bf U}^1_{6,6}(0,x_0)$ against dimension $d\in \{5,6,\ldots,100\}$ for the pricing with different interest rates example in Subsection \ref{subsec:borrowlend}. Right: Average runtime needed to compute 
20 realizations of ${\bf U}^1_{4,4}(0,x_0)$ against dimension $d\in \{5,6,\ldots,100\}$ for the Allen-Cahn equation in Subsection \ref{subsec:allencahn}.}
\label{fig:borrowlendmultidim}
\end{figure}

Throughout this subsection assume the setting in the beginning of Section~\ref{sec:numerics} and assume
for all
$
  s \in [0,T]$, $t \in [s,T]$, $x = ( x_1, \dots, x_d ) \in \R^d$,  
  $y \in \R$, $z\in \R^d$, 
  $k \in \N_0$, 
  $\rho \in \N$, 
  $\theta \in \Theta
$
that
$T=1$, 
$
  \eta( x ) = x
$, 
$
  \mu( s,x ) = 0
$, 
$
  \sigma(s, x) =\operatorname{I}_{\R^{d\times d}} $, 
$
  x_0 =  ( 0, 0, \dots, 0 ) \in \R^{d }
$, 
$
  \mathcal{X}^{ k, \rho, \theta }_{ x, s, t } =x+W^{ \theta }_t - W^{ \theta }_s
$, 
$
\mathcal{D}^{ k, \rho, \theta }_{ x, s, t } 
  =\operatorname{I}_{\R^{d\times d}}
$, 
$
f(s,x,y,z)=y-y^3
$, 
and
$g(x) = \frac{ 1 }{ 1 + \max\{ | x_1 |^2, ..., | x_d |^2 \} } 
$.
Note that the solution $u$ of the PDE \eqref{eq:PDE_quasilinear_1} satisfies for all $t \in [0,T)$, $x\in \R^d$
that
$
  u(T,x) = \frac{ 1 }{ 1 + \max\{ | x_1 |^2, ..., | x_d |^2 \} }  
$
and
\begin{equation}  
\label{eq:PDE_allencahn}
\begin{split}
&
  ( \tfrac{ \partial }{ \partial t } u )( t, x )
  + 
  u(t,x)-\big[u(t,x)\big]^3
  +
  \tfrac{1}{ 2 }\big(\Delta_x u \big)(t,x)
  = 0
  .
\end{split}     
\end{equation}
{\sc Matlab} code \ref{code:dataallencahnd100} presents the 
parameter values in the case $d=100$. In the case $d=1$ line 3
of  {\sc Matlab} code \ref{code:dataallencahnd100}
is replaced by \texttt{dim=1;}.
The simulation results are shown in Figure \ref{fig:allencahn}, the right-hand side of Figure \ref{fig:borrowlendmultidim}, and Tables \ref{tab:allencahnd1} 
and \ref{tab:allencahnd100}.
The left-hand side of Figure \ref{fig:allencahn} suggests an empirical convergence rate close to $\nicefrac{1}{4}$ in the case $d=1$. Moreover, the right-hand side of Figure \ref{fig:allencahn} suggests an empirical convergence rate close
$\nicefrac{1}{3}$ in the case $d=100$. 
\begin{lstlisting}[
	caption={A {\sc Matlab} function that returns the parameter values for the 
	Allen-Cahn example.},
	label=code:dataallencahnd100]
function [T,dim,f,g,eta,mu,sigma] = modelparameters()
    T = 1;
    dim=100;
    mu=0;
    sigma=1;
    f = @(t,x,y,z) y-y.^3; 
    g = @(x) 1./(1+max(x.^2,[],1));
    eta=@(x) x;
end
\end{lstlisting}

\begin{table}[htb]
\begin{center}
\begin{tabular}{|r|rrrrr|}\hline 
$\rho$ & 1 & 2 & 3 & 4 & 5  \\ \hline 
average runtime in seconds  &    0.005 &    0.035 &    0.237 &    7.402 &  345.124\\ 
$ \overline {{\bf U}}^{[1]}_{ \rho ,  \rho  }(0, x_0 )=\frac{ 1 }{ 10 }\sum_{ i = 1 }^{ 10 }{\bf U}^{i,[1]}_{ \rho ,  \rho  }(0, x_0) $ &    1.027 &    0.866 &    0.918 &    0.894 &    0.897\\ 
$ \sqrt{\frac{1}{9}\sum_{ i = 1 }^{ 10 }|{\bf U}^{i,[1]}_{ \rho ,  \rho  }(0, x_0 )-\overline {{\bf U}}^{[1]}_{ \rho ,  \rho  }(0, x_0 )|^2}$ &    0.219 &    0.131 &    0.078 &    0.037 &    0.013\\ 
 $\frac{ 1 }{ 10 }\sum_{ i = 1 }^{ 10 }\frac{|{\bf U}^{i,[1]}_{ \rho ,  \rho  }(0, x_0 )-\texttt v|}{ |\texttt v|} $ &   0.2462 &   0.1197 &   0.0691 &   0.0302 &   0.0124 \\ \hline 
\end{tabular}

\caption{Average runtime, empirical mean, empirical standard deviation, and relative approximation error in the case $d=1$ for the Allen-Cahn equation in Subsection \ref{subsec:allencahn}. The approximation by the finite difference approximation scheme  in {\sc Matlab} code \ref{code:counterfinitediffabm} yields $\texttt{v=approximateUfinitediffgbm(0,x0,2\^{}11)} \approx 0.905$.}\label{tab:allencahnd1}
\end{center}
\end{table}

\begin{table}[htb]
\begin{center}
\begin{tabular}{|r|rrrrr|}\hline 
$\rho$ & 1 & 2 & 3 & 4 & 5  \\ \hline 
average runtime in seconds  &    0.002 &    0.043 &    0.280 &    9.687 &  453.418\\ 
$ \overline {{\bf U}}^{[1]}_{ \rho ,  \rho  }(0, x_0 )=\frac{ 1 }{ 10 }\sum_{ i = 1 }^{ 10 }{\bf U}^{i,[1]}_{ \rho ,  \rho  }(0, x_0) $ &    0.246 &    0.284 &    0.313 &    0.319 &    0.317\\ 
$ \sqrt{\frac{1}{9}\sum_{ i = 1 }^{ 10 }|{\bf U}^{i,[1]}_{ \rho ,  \rho  }(0, x_0 )-\overline {{\bf U}}^{[1]}_{ \rho ,  \rho  }(0, x_0 )|^2}$  &    0.043 &    0.013 &    0.007 &    0.004 &    0.002\\ 
 ${\frac{1}{10} \sum_{ i = 1 }^{ 10 }\frac{|{\bf U}^{i,[1]}_{ \rho+1 ,  \rho+1  }(0, x_0 )-{\bf U}^{i,[1]}_{ \rho ,  \rho  }(0, x_0 )|  }{  |\overline {{\bf U}}^{[1]}_{ 5 ,  5  }(0, x_0 )| }}$  &   0.1484 &   0.0909 &   0.0254 &   0.0102 & \\ \hline 
\end{tabular}

\caption{Average runtime, empirical mean, empirical standard deviation, and relative approximation increments in the case $d=100$ for the Allen-Cahn equation in Subsection \ref{subsec:allencahn}.}\label{tab:allencahnd100}
\end{center}
\end{table}

\subsection{An example with an explicit solution}\label{subsec:exp}

In this subsection we discuss an example with an explicit solution whose three-dimensional version has been considered in Chassagneux~\cite{Chassagneux2014}.

\begin{figure}[htb]
\begin{center}
\includegraphics[width=0.5\textwidth]{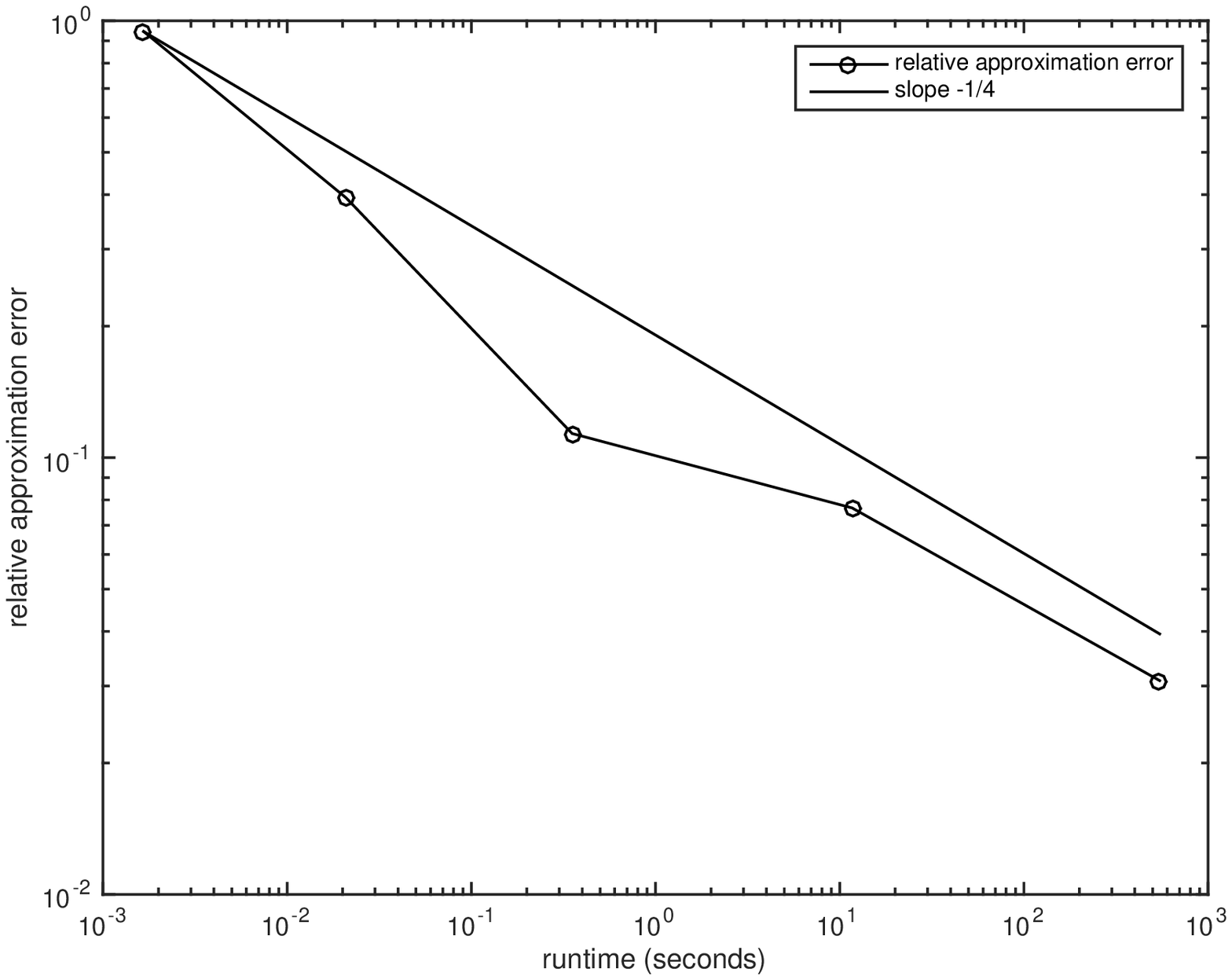}\includegraphics[width=0.5\textwidth]{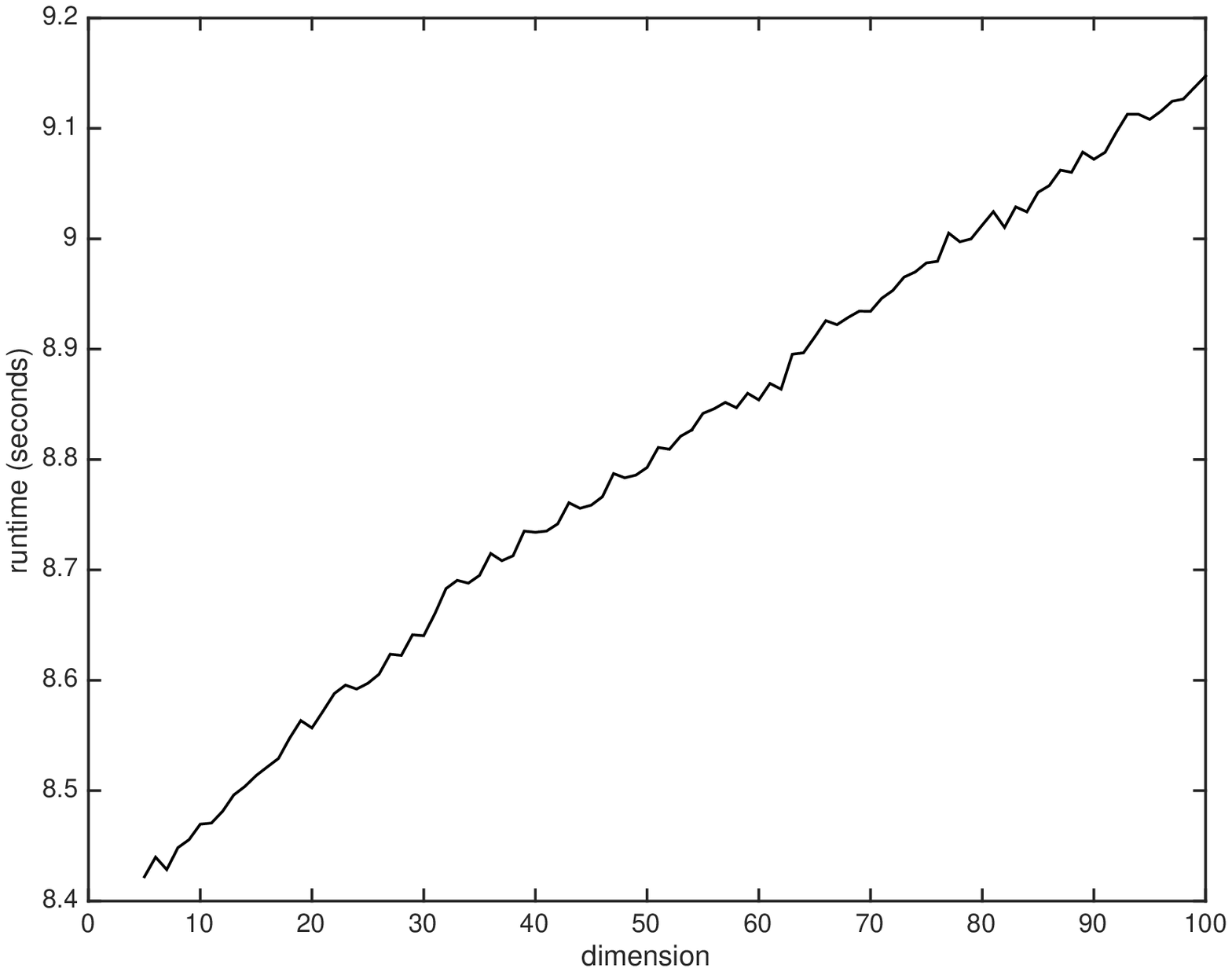}
\end{center}
\caption{ Performance of the scheme \eqref{eq:scheme_Laplace} for the example PDE of Subsection \ref{subsec:exp}.
  Left: Relative 
approximation errors 
$
\scriptstyle{
    \frac{ 1 }{ 10  |\texttt v |}
    \sum_{ i = 1 }^{ 10 }
    |
      {\bf U}^{i,[1]}_{ \rho ,  \rho  }(0, x_0 )
      -
      \texttt v
    |
   }
$
for $\rho \in \{1,2,\ldots, 5\}$
against the average 
runtime for the case $d=100$. Right: Average runtime needed to compute 
20 realizations of ${\bf U}^1_{4,4}(0,x_0)$ against dimension $d\in \{5,6,\ldots,100\}$.}
\label{fig:exp}
\end{figure}

Throughout this subsection assume the setting in the beginning of Section~\ref{sec:numerics}, let $\bar \sigma = 0.25$,
and assume
for all
$
  s \in [0,T]$, $t \in [s,T]$, $x = ( x_1, \dots, x_d ) \in \R^d$, 
  $y \in \R$, $z=(z_1,\ldots, z_d)\in \R^d$,  
  $k \in \N_0$, 
  $\rho \in \N$, 
  $\theta \in \Theta
$
that
$T=0.5$, $d=100$, 
$
  \eta( x ) = x
$, 
$
  \mu( s,x ) = 0
$, 
$
  \sigma(s, x) =  \bar \sigma \operatorname{I}_{\R^{d\times d}}
$, 
$
  x_0 =  ( 0, 0, \dots, 0 ) \in \R^{d }
$, 
$
  \mathcal{X}^{ k, \rho, \theta }_{ x, s, t } =x+W^{ \theta }_t - W^{ \theta }_s
$, 
$
\mathcal{D}^{ k, \rho, \theta }_{ x, s, t } 
  =\operatorname{I}_{\R^{d\times d}}
$, 
$
f(s,x,y,z)=\bar \sigma \big(y-\frac{2+\bar\sigma^2 d}{2\bar \sigma^2 d}\big)\big(\sum_{i=1}^d z_i\big)
$, 
and
$
g(x)=\frac {\exp(T+\sum_{i=1}^d x_i)}{1+\exp(T+\sum_{i=1}^d x_i)}
$.
Note that the solution $u$ of the PDE \eqref{eq:PDE_quasilinear_1} satisfies for all $t \in [0,T)$, $x=(x_1,x_2,\ldots,x_d)\in \R^d$
that
$
  u(T,x) = \frac {\exp(T+\sum_{i=1}^d x_i)}{1+\exp(T+\sum_{i=1}^d x_i)}
$
and
\begin{equation}  
\label{eq:PDE_exp}
\begin{split}
&
  ( \tfrac{ \partial }{ \partial t } u )( t, x )
  + 
    \left[\bar \sigma^2 u(t,x)-\tfrac{1}{ d}-\tfrac {\bar \sigma^2}{2} \right]\bigg[\sum_{i=1}^d  ( \tfrac{ \partial }{ \partial x_i } u )( t, x )\bigg]
  +
  \tfrac{\bar \sigma^2}{ 2 }\big(\Delta_x u\big)(t,x)
  = 0
  .
\end{split}     
\end{equation}
Next we observe
 that $u$ satisfies for all $s\in[0,T]$, $x= ( x_1, \dots, x_d )\in \R^d$ that
\begin{equation}\label{eq:exp_sol}
u(s,x)=\frac {\exp(s+\sum_{i=1}^d x_i)}{1+\exp(s+\sum_{i=1}^d x_i)}.
\end{equation}
{\sc Matlab} code \ref{code:dataexpd100} presents the 
parameter values. 
The simulation results are shown in Figure \ref{fig:exp} and Table \ref{tab:expd100}.
The left-hand side of Figure \ref{fig:exp} suggests an empirical convergence rate close to $\nicefrac{1}{4}$.
\begin{lstlisting}[
	caption={A {\sc Matlab} function that returns the parameter values for the 
	Allen-Cahn example.},
	label=code:dataexpd100]
function [T,dim,f,g,eta,mu,sigma] = modelparameters()
    T = 0.5;
    dim=100;
    mu=0;
    sigma=0.25;
    f = @(t,x,y,z) sigma*(y-(2+sigma^2*dim)/(2*sigma^2*dim)).*sum(z,1); 
    g = @(x) 1-1./(1+exp(T+sum(x,1)));
    eta=@(x) x;
end
\end{lstlisting}

\begin{table}[htb]
\begin{center}
\begin{tabular}{|r|rrrrr|}\hline 
$\rho$ & 1 & 2 & 3 & 4 & 5  \\ \hline 
average runtime in seconds  &    0.002 &    0.021 &    0.352 &   11.677 &  545.871\\ 
$ \overline {{\bf U}}^{[1]}_{ \rho ,  \rho  }(0, x_0 )=\frac{ 1 }{ 10 }\sum_{ i = 1 }^{ 10 }{\bf U}^{i,[1]}_{ \rho ,  \rho  }(0, x_0) $ &    0.751 &    0.522 &    0.523 &    0.520 &    0.495\\ 
$ \sqrt{\frac{1}{9}\sum_{ i = 1 }^{ 10 }|{\bf U}^{i,[1]}_{ \rho ,  \rho  }(0, x_0 )-\overline {{\bf U}}^{[1]}_{ \rho ,  \rho  }(0, x_0 )|^2}$ &    0.477 &    0.248 &    0.070 &    0.040 &    0.018\\ 
 $\frac{ 1 }{ 10 }\sum_{ i = 1 }^{ 10 }\frac{|{\bf U}^{i,[1]}_{ \rho ,  \rho  }(0, x_0 )-\texttt v|}{ |\texttt v|} $ &   0.9449 &   0.3931 &   0.1135 &   0.0767 &   0.0309 \\ \hline 
\end{tabular}

\caption{Average runtime, empirical mean, empirical standard deviation, and relative approximation error in the case $d=100$ for the example PDE of Subsection \ref{subsec:exp}. The exact value of the solution is $\texttt{v}=u(0,0)=0.5$.}\label{tab:expd100}
\end{center}
\end{table}

\section{Discussion of approximation methods from the literature}\label{sec:lit_discussion}

Deterministic methods for second-order parabolic PDEs  are known to have exponentially growing computational effort.
Since a program with $10^{80}$, say, floating point operations will never terminate (on a non-quantum computer),
deterministic methods
such as finite elements methods, finite difference methods, spectral Galerkin approximation methods, or sparse grid methods
are not suitable for solving high-dimensional nonlinear second-order parabolic PDEs no matter what
the convergence rate of the method is.
For this reason we discuss only stochastic approximation methods for nonlinear second-order parabolic PDEs.
In the literature we have found the following articles
\cite{
BriandDelyonMemin2001,
MaProtterSanMartinTorres2002,
BallyPages2003a,
BallyPages2003b,
BouchardTouzi2004,
Zhang2004,
GobetLemorWarin2005,
LemorGobetWarin2006,
DelarueMenozzi2006,
BenderDenk2007,
GobetLemor2008,
CrisanManolarakis2010,
CrisanManolarakisTouzi2010,
GobetLabart2010,
CrisanManolarakis2012,
BriandLabart2014,
ChassagneuxCrisan2014,
Chassagneux2014, 
CrisanManolarakis2014,
LionnetReisSzpruch2015,
Turkedjiev2015,
RuijterOosterlee2015,
RuijterOosterlee2016,
GobetTurkedjiev2016regression, 
GobetTurkedjiev2016Malliavin, 
Henry-Labordere2012,
Henry-LabordereTanTouzi2014,
Henry-LabordereEtAl2016,
GeissLabart2016,
ChangLiuXiong2016,
LeCavilOudjaneRusso2016%
}
which propose (possibly non-implementable) stochastic approximation methods for nonlinear second-order parabolic PDEs.
All of these methods except for
\cite{Henry-Labordere2012,Henry-LabordereTanTouzi2014,Henry-LabordereEtAl2016,
ChangLiuXiong2016,LeCavilOudjaneRusso2016}
exploit a stochastic representation with BSDEs due to Pardoux \& Peng~\cite{PardouxPeng1992}.
Moreover, all of these methods except for
\cite{GobetLabart2010,BriandLabart2014,GeissLabart2016,Henry-Labordere2012,Henry-LabordereTanTouzi2014,Henry-LabordereEtAl2016,
ChangLiuXiong2016,LeCavilOudjaneRusso2016}
can be described in two steps.
In the first step, time in the corresponding BSDE is discretized backwards in time via an explicit or an implicit Euler-type method
which was investigated in detail, e.g., in Bouchard \& Touzi~\cite{BouchardTouzi2004} and Zhang~\cite{Zhang2004}.
The resulting approximations involve nested conditional expectations and, therefore, are not implementable.
In the second step, these conditional expectations are approximated
by 'straight-forward' Monte Carlo simulations,
by the quantization tree method (proposed in~\cite{BallyPages2003a}),
by a regression method based on kernel-estimation or on Malliavin calculus (proposed in~\cite{BouchardTouzi2004}),
by projections on function spaces (proposed in~\cite{GobetLemorWarin2005}),
or
by the cubature method on Wiener space (developed in~\cite{LyonsVictoir2004} and proposed in~\cite{CrisanManolarakis2012}).
The first step does not cause problems in high dimensions in the sense that
the backward (explicit or implicit) Euler-type approximations converge
under suitable assumptions with rate at least $0.5$
(see Theorem 5.3 in Zhang~\cite{Zhang2004} and Theorem 3.1 in Bouchard \& Touzi~\cite{BouchardTouzi2004}
for the backward implicit Euler-type method)
and the computational effort (assuming the conditional expectations are known exactly) grows at most linearly
in the dimension for fixed accuracy.
For this reason, we discuss below in detail only the different methods for discretizing conditional expectations.
In addition, we discuss
the Wiener chaos decomposition method proposed in~\cite{BriandLabart2014},
the branching diffusion method proposed in~\cite{Henry-Labordere2012},
and methods based on density estimation proposed
in~\cite{ChangLiuXiong2016,LeCavilOudjaneRusso2016}.

A difficulty in our discussion below is that the discussed algorithms (except for the branching diffusion method)
depend on different parameters
and the optimal choice of these parameters is unknown since no lower estimates for the approximation errors are known.
For this reason we will choose parameters which are optimal with respect to the best known upper error bound.
For these parameter choices we will show below
for the discussed algorithms (except for the branching diffusion method)
that
the computational effort fails to grow at most polynomially 
both
in the dimension and in the
reciprocal of the best known upper error bound.

Throughout this section assume the setting in Subsection \ref{sec:algorithm_semilinear},
 let $u^{\infty}\in C^{1,2}([0,T]\times \R^d, \R)$ be a function which satisfies~\eqref{eq:PDE_quasilinear_1}
and we denote by $Y\colon[0,T]\times\Omega\to\R$ the stochastic process which satisfies for all $t\in[0,T]$ that
$Y_t=u^{\infty}(t,W_t^0)$.

\subsection{The 'straight-forward' Monte Carlo method}\label{ssec:MC}
The 'straight-forward' Monte Carlo method 
approximates the conditional expectations involved in backward Euler-type approximations by Monte Carlo simulations.
The resulting nesting of Monte Carlo averages is computational expensive in the following sense.
If for a set $\pi\subseteq[0,T]$ with $|\pi|\in\N$  many points
and for $M\in\N$ the random variable $Y^{\pi,M}\colon\Omega\to\R$ is
the 'straight-forward' Monte Carlo approximation of $Y$ with time grid $\pi$ and $M$ Monte Carlo averages
for each involved conditional expectation,
then
the number of realizations of scalar standard normal random variables
required to compute
$Y^{\pi,M}$
is $(Md)^{|\pi|}$
and
the $L^2$-error 
satisfies for a suitable constant $c\in\R$ independent of $\pi$ and $N$ that
 \begin{equation}\label{eq:error.brute.force}
  \max_{t\in\pi}\|Y_{t}- Y^{\pi,M}_{t}\|_{L^2(\P;\R)}
  \leq c\left(|\pi|^{-\frac{1}{2}}+|\pi|M^{-1/2}\right),
\end{equation}
see, e.g., 
Theorem 4.3 and Display (4.14) in Crisan \& Manolarakis~\cite{CrisanManolarakis2010}.
Thus the computational effort
$(Md)^{(\frac{1}{|\pi|^{-1/2}})^2}
\geq
(Md)^{(\frac{c}{c|\pi|^{-1/2}+c|\pi|M^{-1/2}})^2}
$
grows at least exponentially in the
reciprocal of the right-hand side of~\eqref{eq:error.brute.force}.
This suggests an at most logarithmic convergence rate
of the 'straight-forward' Monte Carlo method.
We are not aware
of a statement in the literature  claiming that the 'straight-forward' Monte Carlo method has a polynomial convergence rate.

\subsection{The quantization tree method}
The quantization tree method has been introduced in
Bally \& Pag\`es~\cite{BallyPages2003a,BallyPages2003b}.
In the proposed algorithm, time is discretized by the explicit backward Euler-type method.
Moreover, 
one chooses a space-time grid
and computes the transition probabilities for the underlying forward
diffusion projected to this grid by Monte Carlo simulation.
With these discrete-space transition probabilities one can then approximate
all involved conditional expectations.
If for a set $\pi\subseteq[0,T]$ with $|\pi|\in\N$  many points
and for $N\in\N$ the random variable $Y^{\pi,N}\colon\Omega\to\R$ is
the quantization tree approximation of $Y$ with time grid $\pi$, a specific space grid with a total number of $N$ nodes
and explicitly known transition probabilities of the forward diffusion
and if the coefficients are sufficiently regular,
then
the number of realizations of scalar standard normal random variables
required to compute
$Y^{\pi,N}$
is at least $Nd|\pi|$
and
Display (6) in Bally \& Pag\`es~\cite{BallyPages2003b} shows for optimal grids 
and a constant $c\in\R$ independent of $\pi$ and $N$ that
\begin{equation}
  \max_{t\in\pi}\|Y_{t}- Y^{\pi,N}_{t}\|_{L^2(\P;\R)}
  \leq c\left(\frac{1}{|\pi|}+\frac{|\pi|^{1+1/d}}{N^{1/d}}\right).
\end{equation}
To ensure that this upper bound does not explode as $|\pi| \to \infty$ it is thus necessary
to choose a space-time grid with at least $N=|\pi|^{d+1}$ many nodes when there are $|\pi|\in \N$ many time steps. With this choice
the computational effort
of this algorithm grows exponentially fast in the dimension.
We have not found a statement in the literature on the quantization tree method
claiming that there exists a choice of parameters such that
the computational effort grows at most polynomially both in the dimension and in the reciprocal
of the prescribed accuracy.

\subsection{The Malliavin calculus based regression method}
The Malliavin calculus based regression method has been introduced
in Section 6
in Bouchard \& Touzi~\cite{BouchardTouzi2004} and is based on the implicit backward Euler-type method.
The algorithm involves iterated Skorohod integrals
which by Display (3.2) in
Crisan, Manolarakis, \& Touzi~\cite{CrisanManolarakisTouzi2010}
can be numerically computed with $2^d$ many independent
standard normally distributed random variables. 
In that case the computational effort grows exponentially fast in the dimension. We are not aware of an approximation method of the involved iterated Skorohod integrals whose computational effort does not grow
exponentially fast in the dimension.
Example 4.1 
in Bouchard \& Touzi~\cite{BouchardTouzi2004} also mentions a method for approximating all involved conditional expectations
using kernel estimation. For this method we have not found an upper error estimate in the literature so that we do not
known how to choose the bandwidth matrix of the kernel estimation given the number of time grid points.

\subsection{The projection on function spaces method}
The projection on function spaces method has been proposed in Gobet, Lemor, \& Warin~\cite{GobetLemorWarin2005}.
The algorithm is based on estimating the involved conditional
expectations by considering the projections of the random variables
on a finite-dimensional function space and then estimating these projections by
Monte Carlo simulation.
In general the projection error 
and the computational effort depend on the choice of the basis functions.
In the literature we have found the following two choices of basis functions.
In Gobet, Lemor, \& Warin~\cite{GobetLemorWarin2005} (see also Gobet \& Lemor~\cite{GobetLemor2008}
and Lemor, Gobet, \& Warin~\cite{LemorGobetWarin2006}) indicator functions of hypercubes are employed as basis functions.
In this case there exists $c\in \R$ such that a projection error $\eps \in (0,\infty)$
can be achieved by simulating $\lfloor c\eps^{-(3+2d)}|\log(\eps)| \rfloor$ paths of the forward diffusion. 
With this choice, the computational effort of the algorithm grows exponentially fast
in the dimension for fixed accuracy $\eps \in (0,1)$.
Ruijter \& Oosterlee~\cite{RuijterOosterlee2015}
 use
certain cosine functions as basis functions
and motivate this with
a
Fourier cosine series expansion.
The resulting approximation method has only been specified in a
one-dimensional setting so that the computational effort
in a high-dimensional setting remained unclear.
We have not found a statement in the literature on the projection on function spaces method claiming that
there exists a choice of function spaces and other algorithm parameters such that
the computational effort of the method grows at most polynomially both in the dimension of the PDE and in the reciprocal of the prescribed
accuracy.

\subsection{The cubature on Wiener space method}
The cubature on Wiener space method  for approximating solutions of PDEs has been introduced
in Crisan \& Manolarakis \cite{CrisanManolarakis2012}.
This method combines
the implicit backward Euler-type scheme with the cubature method developed in Lyons \& Victoir~\cite{LyonsVictoir2004}
for constructing finitely supported measures that approximate the distribution of the solution of
a stochastic differential equation.
This method has a parameter $m\in\N$
and constructs for every finite time grid $\pi\subseteq[0,T]$ (with $|\pi|\in\N$ points)
a sequence
$w_1,\ldots,w_{(N_{m,d})^{|\pi|}}\in C^0([0,1],\R^d)$ of paths with bounded variation
where $N_{m,d}\in\N$ is the number of nodes needed for a cubature formula of degree $m$ with
respect to the $d$-dimensional Gaussian measure.
We note that this construction is independent of $f$ and $g$ and can be computed once and then tabularized.
Using these paths,
Corollary 4.2 in
Crisan \& Manolarakis \cite{CrisanManolarakis2012}
shows in the case $m\geq 3$ that there exists a constant $c\in[0,\infty)$,
a sequence $\pi_n\subseteq[0,T]$, $n\in\N$, of finite time grids 
and
there exist implementable approximations
$ Y ^{n}\colon\pi_n\times\Omega\to\R$, $n\in\N$, of the exact solution $Y$
such that for all $n\in\N$ it holds that $0\in\pi_n$, $\pi_n$ has $n+1$ elements
and 
 \begin{equation}\label{eq:upper.bound.cubature}
  \|Y_{0}- Y^{\pi_n}_{0}\|_{L^2(\P;\R)}\leq \tfrac{c}{n}.
\end{equation}
In this form of the algorithm, the computational effort for calculating $Y^{\pi_n}$, which is at least the number
$(N_{m,d})^n$ of paths to be used,
grows exponentially in the reciprocal of the right-hand side of~\eqref{eq:upper.bound.cubature}.
To avoid this exponential growth of the computational effort in the number of cubature paths,
Crisan \& Manolarakis~\cite{CrisanManolarakis2012}
specify two methods (a tree based branching algorithm of Crisan \& Lyons~\cite{CrisanLyons2002}
and the recombination method of Litterer \& Lyons~\cite{LittererLyons2012})
which reduce the number of nodes  and which result
in approximations which converge with polynomial rate; cf.\ Theorem 5.4 in\cite{CrisanManolarakis2012}.
The constant in the upper error estimate in Theorem 5.4 in~\cite{CrisanManolarakis2012}
may depend on the dimension (cf.\ also (5.16) and the proof of Lemma 3.1 in~\cite{CrisanManolarakis2012}).
Simulations in the literature on the cubature method were performed in dimension 1
(see Figures 1--4
 in~\cite{CrisanManolarakis2012} and Figure 1 in~\cite{CrisanManolarakis2014})
 or dimension 5
(see Figures 5--6
 in~\cite{CrisanManolarakis2012}).
 To the best of our knowledge, there exist no statement in the literature on the cubature method which asserts that the computational effort
of the cubature method together with a suitable complexity reduction method
grows at most polynomially both in the dimension of the PDE and in the reciprocal of the prescribed
accuracy.

\subsection{The Wiener chaos decomposition method}\label{ssec:WienerChaos}
The Wiener chaos decomposition method has been introduced in
Briand \& Labart~\cite{BriandLabart2014} and has been extended to the case of BSDEs with jumps in Geiss \& Labart~\cite{GeissLabart2016}.
The algorithm is based on Picard iterations of the associated BSDE and evaluates the involved nested conditional expectations using Wiener chaos decomposition formulas. 
This algorithm does not need to discretize time since Wiener integrals over integrands with explicitly known antiderivative
can be simulated exactly.
The computational complexity of approximating the solution of a BSDE of dimension $d$
using a Wiener chaos decomposition of order $p\in \N$, $K\in \N$ Picard iterations, $M\in \N$ Monte Carlo samples for each
conditional expectation, and $N\in \N$ many time steps is of order $O(K\times M \times p \times (N\times d)^{p+1})$; see Section 3.2.2 in Briand \& Labart~\cite{BriandLabart2014}.
To ensure that the approximation error converges to $0$ requires the order $p$ of the chaos decomposition to increase to $\infty$. 
This implies that the computational effort fails to grow at most polynomially both in the dimension of the PDE and in the reciprocal of the prescribed
accuracy.
We also mention that we are not aware of a result in the literature that establishes a polynomial rate of convergence for the Wiener chaos decomposition method (see, e.g., Remark 4.8 in Briand \& Labart~\cite{BriandLabart2014}).

\subsection{The branching diffusion method}\label{ssec:branchingdiffusions}

The branching diffusion method has been proposed in Henry-Labord\`ere~\cite{Henry-Labordere2012};
see also the extensions to the non-Markovian case in~\cite{Henry-LabordereTanTouzi2014}
and to nonlinearities depending on derivatives in~\cite{Henry-LabordereEtAl2016}.
This method approximates the nonlinearity $f$ by polynomials and then exploits that the solution of a semilinear PDE with polynomial nonlinearity
(KPP-type equations) can be represented as an expectation of a functional of a branching diffusion process
due to Skorohod~\cite{Skorohod1964}.
This expectation can then be numerically approximated with the standard Monte Carlo method
and pathwise approximations of the branching diffusion process.
The branching diffusion method does not suffer from the 'curse of dimensionality by construction' 
and works in all dimensions.
It's convergence rate is $0.5$ if the forward diffusion can be simulated exactly
and, in general, its rate is $0.5-$ using a pathwise approximation of the forward diffusion and 
the multilevel Monte Carlo method proposed in Giles~\cite{g08a}.

The major drawback of the branching diffusion method is its insufficient applicability.
This method replaces potentially 'nice' nonlinearities by potentially 'non-nice' polynomials.
Semilinear PDEs with certain polynomial nonlinearities, however, can 'blow up' in finite time;
see, e.g., Fujita~\cite{Fujita1966}, Escobedo \& Herrero~\cite{EscobedoHerrero1991} for analytical proofs
and, e.g.,
Nagasawa \& Sirao~\cite{NagasawaSirao1969} and Lopez-Mimbela \& Wakolbinger~\cite{Lopez-MimbelaWakolbinger1998}
for probabilistic proofs.
If the approximating polynomial nonlinearity, the time horizon, and the terminal condition 
satisfy a certain condition,
then the PDE does not 'blow up' until time $T$ and the branching diffusion method is known to work well.
More specifically, if
there exist 
$\beta\in(0,\infty)$ and functions $(a_k)_{k\in\N_0}\colon[0,T]\times\R^d\to\R$
such that
for all $(t,x,y,z)\in[0,T]\times\R^d\times\R\times\R^d$
it holds
that $f(t,x,y,z)=\beta\sum_{k=0}^{\infty}a_k(t,x)y^k-\beta y$,
if the functions $\mu$ and $\sigma$ are bounded, continuous and Lipschitz in the second argument,
and if $\forall x\in\R^d\colon\eta(x)=x$,
then Theorem 2.13 in~\cite{Henry-LabordereTanTouzi2014}
(see also Proposition 4.2 in~\cite{Henry-Labordere2012} or Theorem 3.12 in~\cite{Henry-LabordereEtAl2016})
shows that a sufficient condition for a stochastic representation with a branching diffusion to hold is that
\begin{equation}  \begin{split}\label{eq:sufficient.branching.representation}
  \int_{\sup_{x\in\R^d}|g(x)|}^{\infty}\tfrac{1}{\beta\max\left\{0,\sum_{k=0}^{\infty}\sup_{(t,x)\in[0,T]\times\R^d}|a_k(t,x)|y^k-y\right\}}\,dy> T.
\end{split}     \end{equation}
For the branching diffusion method to converge with rate $0.5$ the random variables in the stochastic representation
need to have finite second moments which leads to a more restrictive condition than~\eqref{eq:sufficient.branching.representation};
see Remark 2.14 in~\cite{Henry-LabordereTanTouzi2014}.
However,
condition~\eqref{eq:sufficient.branching.representation} is also necessary
for the stochastic representation in~\cite{Henry-LabordereTanTouzi2014} to hold
if the functions $g$ and $(a_k)_{k\in\N_0}$ are constant and positive
and
if $\mu$ and $\sigma$ are constant
(then the PDE~\eqref{eq:PDE_quasilinear_1} reduces to an ODE for which the 'blow-up'-behavior is well-known);
see, e.g.,
Lemma 2.5 in~\cite{Henry-LabordereTanTouzi2014}.

The branching diffusion method also seems to have problems with polynomial nonlinearities
where the exact solution does not 'blow up' in finite time.
Since no theoretical results are available in this direction,
we illustrate this with simulations for 
an Allen-Cahn equation
(a simplified version of the Ginzburg-Landau equation).
More precisely, for the rest of this subsection assume that
$T$, $\mu$, $\sigma$, $f$, and $g$ satisfy
for all $(t,x,y,z)\in[0,T]\times\R^d\times\R\times\R^d$
that $T=1$, $\mu(x)=\mu(0)$, $\sigma(x)=\sigma(0)$, $f(t,x,y,z)=y-y^3$, $g(x)=g(0)$, and $g(0)\ge 0$.
Then~\eqref{eq:PDE_quasilinear_1} is an ODE and the solution $u^{\infty}$ satisfies
for all $(t,x)\in[0,T]\times\R^d$ that $u^{\infty}(t,x)=\left(1-(1-(g(0))^{-2})e^{2t-2}\right)^{-\frac{1}{2}}$.
In this situation
the sufficient condition~\eqref{eq:sufficient.branching.representation}
(choose for all $(t,x)\in[0,T]\times\R^d$ and $k\in\N_0\setminus\{1,3\}$
that $\beta=1$, $a_1(t,x)=2$, $a_3(t,x)=-1$, and $a_k(t,x)=0$)
from
Theorem 2.13 in~\cite{Henry-LabordereTanTouzi2014}
is equivalent to
\begin{equation}  \begin{split}
  1<\int_{|g(0)|}^{\infty}\tfrac{1}{y+y^3}\,dy
  =\left[\log(y)-\tfrac{1}{2}\log(1+y^2)\right]_{|g(0)|}^{\infty}=\tfrac{1}{2}\log\left(1+\tfrac{1}{|g(0)|^2}\right)
\end{split}     \end{equation}
which is equivalent to $|g(0)|<(e^2-1)^{-\frac{1}{2}}=0.395623\ldots$
We simulated the branching diffusion approximations of $u^{\infty}(0,0)$ for different
values of $g(0)$. Each approximation is an average over $M=10^5$ independent copies $(\Psi_{0,0}^{i})_{i\in\{1,2,\ldots,10^5\}}$ of
the random variable $\Psi_{0,0}$ defined in (2.12)
in Henry-Labord\`ere, Tan, \& Touzi~\cite{Henry-LabordereTanTouzi2014}
where we choose for all $k\in\N_0$ that $p_k=0.5\1_{\{1,3\}}(k)$.
We also report the estimated standard deviation
$\big(10^{-5}\sum_{i=1}^{10^5}(\Psi_{0,0}^i)^2-\big[10^{-5}\sum_{i=1}^{10^5}\Psi_{0,0}^i\big]^2\big)/\sqrt{10^5}$
of the branching diffusion approximation $10^{-5}\sum_{i=1}^{10^5}\Psi_{0,0}^i$.
\begin{table}[htb]
\begin{center}
\begin{tabular}{|c||c|c|c|c|c|} \hline
$g(0)$ & exact value $u^{\infty}(0,0)$ & approximation & estimated standard deviation \\
\hline\hline
$0.1$   &  0.263540   & 0.271007        &   0.0169791  \\
$0.2$   &  0.485183   & 0.499103        &   0.0361975  \\
$0.3$   &  0.649791   & 0.848879        &   0.2211004  \\
$0.4$   &  0.764605   & 3.495457        &   2.8179089  \\
$0.5$   &  0.843347   & 21.68436        &   20.978325  \\
$0.6$   &  0.897811   & 136.6667        &   110.02696  \\
$0.7$   &  0.936233   & 7321.326        &   5404.5849  \\
 \hline
\end{tabular}
\caption{\footnotesize Approximation of the PDE $\tfrac{\partial}{\partial t}u+\tfrac{1}{2}\Delta_xu+u-u^3=0$ defined on $[0,1]\times\R$
with terminal condition $u(1,\cdot)=g(0)$
with the branching diffusion method from Theorem 2.13 in
Henry-Labord\`ere, Tan, \& Touzi~\cite{Henry-LabordereTanTouzi2014}.}
\label{tab1}
\end{center}
\end{table}
Table~\ref{tab1} shows that the branching diffusion approximations of $u^{\infty}(0,0)$ become poor as $g(0)$ increases from $0.1$ to $0.7$.
Thus the branching diffusion method fails to produce good approximations for $u^{\infty}(0,0)$ in our example as soon as
condition~\eqref{eq:sufficient.branching.representation} is not satisfied.

A further minor drawback of the branching diffusion method is that it requires a suitable
approximation of the nonlinearity with polynomials and this might not be available.
In addition, certain functions (e.g.\ for $\R\ni x\mapsto \max\{0,x\}\in\R$)
can only be approximated by polynomials on finite intervals
so that choosing suitable approximating polynomials might require appropriate a priori bounds
on the exact solution of the PDE.

\subsection{Approximations based on density representations}\label{ssec:FokkerPlanck}
Recently two approximation methods were proposed
in Chang, Liu, \& Xiong~\cite{ChangLiuXiong2016}
and in Le Cavil, Oudjane, \& Russo~\cite{LeCavilOudjaneRusso2016}.
Both methods are based
on a stochastic representation with a McKean-Vlasov-type SDE where
$u^{\infty}$ is the density of a certain measure.
As a consequence both methods proposed
in~\cite{ChangLiuXiong2016,LeCavilOudjaneRusso2016} encounter the difficulty
of density estimation in high dimensions.
More precisely if $\bar{u}^{\eps,N,n}$ is the approximation of $u^{\infty}$
defined in (5.33)
in~\cite{LeCavilOudjaneRusso2016} with a uniform time grid with $n\in\N$ time points,
$N\in\N$ Monte Carlo averages, and bandwidth matrix $\eps \operatorname{I}_{\R^{d\times d}}$ where $\eps \in(0,1]$,
then the proof of Theorem 5.6, the proof of Corollary 5.4
in~\cite{LeCavilOudjaneRusso2016} imply under suitable assumptions
existence of constants $C,\bar{C}\in(0,\infty)$
and a function $c\colon(0,1]\to(0,\infty)$ (which are independent of $n,N$ and $\eps$)
such that
\begin{equation}  \begin{split}
  &\sup_{t\in[0,T]}\E\left[\int_{\R^d}\left|\bar{u}^{\eps,N,n}(t,x)-u^{\infty}(t,x)\right|\,dx
  +\int_{\R^d}\left|(\nabla_x(\bar{u}^{\eps,N,n}-u^{\infty}))(t,x)\right|\,dx
  \right]
 \\& \leq \left(\tfrac{\bar{C}}{\eps^{d+3}\sqrt{n}}+\tfrac{C}{\sqrt{\eps^{d+4}N}}+
       \right)e^{\frac{C}{\eps^{d+1}}}+c_{\eps}.
\end{split}     \end{equation}
This upper bound becomes only small if we choose the bandwidth $\eps$ small and if $n$ and $N$
grow exponentially in the dimension.
The upper bounds established
in~\cite{ChangLiuXiong2016}
are less explicit in the dimension.
However, following the estimates 
in the proofs in~\cite{ChangLiuXiong2016},
it becomes apparent that the
number of initial particles in
branching particle system approximations defined on pages 30 and 18
in~\cite{ChangLiuXiong2016} need to grow exponentially in the dimension.

\subsubsection*{Acknowledgement}
This project has been partially supported through the research grants ONR N00014-13-1-0338 and DOE DE-SC0009248
and through the German Research Foundation via RTG 2131 \textit{High-dimensional Phenomena in Probability -- Fluctuations and Discontinuity}
and via research grant HU 1889/6-1.

{\small
\bibliographystyle{acm}
\bibliography{Bib/bibfile}
}

\end{document}